\documentclass[12pt]{article}
\usepackage{amssymb}
\usepackage{amsmath}
\usepackage{amsthm}
\usepackage{mathrsfs}
\usepackage{ifthen}
\usepackage{tikz}
\usepackage{enumerate}
\usepackage{verbatim}
\usepackage{varwidth}

\usepackage{caption}
\usepackage[utf8]{inputenc}

\usepackage[letterpaper]{geometry}
\usetikzlibrary{calc,arrows,decorations.pathmorphing,decorations.pathreplacing,shapes.symbols,backgrounds,fit}

\usepackage{hyperref}
\hypersetup{colorlinks, citecolor=blue}

\newtheorem{thm}{Theorem}[section]
\newtheorem{lem}[thm]{Lemma}
\newtheorem{prop}[thm]{Proposition}
\newtheorem{cor}[thm]{Corollary}
\newtheorem{mainthm}{Theorem}

\theoremstyle{definition}
\newtheorem{defn}[thm]{Definition}

\newtheorem{conv}[thm]{Convention}

\newcommand{\drawpetal}[3]{\draw [shift=#3,scale=#2,rotate=#1](0,0) .. controls (-0.25,0.75) and (-0.2,1) .. (0,1) .. controls (0.2,1) and (0.25,0.75) .. (0,0);}

\newcommand{\drawbouquet}[3]{[scale=#2]
  \foreach \x in {120,0,-40,-80,-120}
  {
    \pgfmathparse{\x+#1}
    \let\spam\pgfmathresult
    \drawpetal{\spam}{#2}{#3}
  }
  \filldraw [shift=#3,scale=#2,rotate=#1] (0,0) circle (0.1)
  (120:0.5) circle (0.03)
  (150:0.5) circle (0.03)
  (180:0.5) circle (0.03);}

\newcommand{\base}[4]{\draw[thick] (#2+0.05,#1) -- node{$\blacktriangleright$}
  node[above]{\ensuremath{#4}} (#2+#3-0.05,#1);}

\newcommand{\opdualbase}[4]{\draw[thick] (#2+0.05,#1) -- node
  {$\blacktriangleright$}
  node[above]{\ensuremath{\overline{#4}}}(#2+#3-0.05,#1);}

\newcommand{\opolpair}[5]{\base{#1}{#2}{#3}{#5}
  \opdualbase{#1+0.5}{#2+#4}{#3}{#5}
}

\newcommand{\opband}[4]{
\begin{scope}

\filldraw[black!10!white,draw=black]
(#1+#4/2,#2+#4/2) arc (90:0:#4/2) -- (#1+#4,#2) -- (#1+#4+#3,#2) arc
(0:90:#4/2) -- (#1+#4/2,#2+#4/2);

\filldraw[black!20!white,draw=black]
(#1,#2) arc (180:90:#4/2) -- (#1+#3+#4/2,#2+#4/2) arc(90:180:#4/2) --
(#1,#2);

\draw[dashed] (#1+#3+#4/2,#2+#4/2) arc(90:180:#4/2);
\end{scope}
}

\newcommand{\opside}[3]{
  \draw[ultra thick](#1,#2) arc (180:0:#3/2);
}

\newcommand{\olband}[4]{
\begin{scope}
\filldraw[black!10!white,draw=black] (#1,#2) .. controls (#1+1,#2+0.5) and (#1+#4/2+1,#2+2)
.. (#1+#4/2,#2+2) .. controls (#1+#4/2-1,#2+2) and (#1+#4-1,#2+0.5)
.. (#1+#4,#2) -- (#1+#3+#4,#2) .. controls (#1+#3+#4-1,#2+0.5) and
(#1+#3+#4/2-1,#2+2) .. (#1+#3+#4/2,#2+2) ..controls (#1+#4/2+#3+1,#2+2)
and (#1+#3+1,#2+0.5) .. (#1+#3,#2);
\filldraw[black!20!white,draw=black] (#1,#2) .. controls (#1+1,#2+0.5) and (#1+#4/2+1,#2+2)
.. (#1+#4/2,#2+2) -- (#1+#4/2+#3,#2+2) ..controls (#1+#4/2+#3+1,#2+2)
and (#1+#3+1,#2+0.5) .. (#1+#3,#2) --cycle;
\draw[dashed](#1+#4/2,#2+2) .. controls (#1+#4/2-1,#2+2) and (#1+#4-1,#2+0.5)
.. (#1+#4,#2)
(#1+#3+#4,#2) .. controls (#1+#3+#4-1,#2+0.5) and
(#1+#3+#4/2-1,#2+2) .. (#1+#3+#4/2,#2+2)
;
\end{scope}
}

\newcommand{\olside}[4]{ 
  \draw[ultra thick] (#1,#2) .. controls (#1+1,#2+0.5) and
  (#1+#4/2+1,#2+2) .. (#1+#4/2,#2+2).. controls (#1+#4/2-1,#2+2) and
  (#1+#4-1,#2+0.5) .. (#1+#4,#2);}

\newcommand{\tubeside}[3]{ 
  \draw[ultra thick] (#1,#2) .. controls (#1+1,#2+0.5) and
  (#1+#3/2+1,#2+2) .. (#1+#3/2,#2+2);
  \draw[ultra thick,dashed] (#1+#3/2,#2+2) .. controls (#1+#3/2-1,#2+2) and
  (#1+#3-1,#2+0.5) .. (#1+#3,#2);}

\newcommand{\squareofnodes}[4]
{\node (1) at (0,0) {\ensuremath{#1}};
  \node (2) at (2,0) {\ensuremath{#2}};
  \node (3) at (0,-1) {\ensuremath{#3}};
  \node (4) at (2,-1) {\ensuremath{#4}};}

\newcommand{\define}[1]{\emph{#1}}

\newcommand{\calH}{\mathcal{H}}

\newcommand{\mo}{{-1}}
\newcommand{\pmo}{{\pm 1}}
\renewcommand{\tilde}{\widetilde}
\newcommand{\Z}{\mathbb{Z}}
\newcommand{\bra}{\langle}
\newcommand{\kett}{\rangle}
\newcommand{\ol}{\overline}
\newcommand{\ncl}[1]{\bk{\bk{#1}}}

\newcommand{\C}{\mathscr{C}}

\newcommand{\N}{\mathbb{N}}
\newcommand{\closure}[1]{\ensuremath{\mathbf{closure}\left(#1\right)}}

\newcommand{\axis}{\mathrm{axis}}
\newcommand{\bk}[1]{{\bra #1 \kett}}

\newcommand{\gp}{\mathrm{Gp}}
\newcommand{\bdy}{\partial}

\newcommand{\complex}[1]{{#1}}
\newcommand{\fungrp}[1]{\pi_{1}(#1)}

\newcommand{\elliptics}{\ensuremath{\mathcal{H}}}
\newcommand{\Band}{\ensuremath{\mathbb{B}}}
\newcommand{\band}[1]{\ensuremath{\mathbb{B}(#1)}}

\newcommand{\mbc}{\ensuremath{\mathscr{C}}}

\newcommand{\Ann}{\ensuremath{\mathcal{A}}}
\newcommand{\ann}[1]{{\ensuremath{\mathcal{A}(#1)}}}
\newcommand{\bandeq}{\doteqdot} 
\newcommand{\bouquet}[1]{\ensuremath{B_{#1}}}

\newcommand{\tube}[1]{\ensuremath{\tau_{#1}}}
\newcommand{\relloop}[2]{\ensuremath{\rho_{\inter{#1}}^{#2}}}
\newcommand{\tr}[1]{\ensuremath{\mathrm{tr}(#1)}}
\newcommand{\inter}[1]{\ensuremath{\sigma(#1)}}
\newcommand{\reltr}[2]{\ensuremath{\mathrm{tr}_{\inter{#1}}(#2)}}
\newcommand{\excess}[2]{\ensuremath{\psi_{#1}\left(#2\right)}}

\newcommand{\frep}{\ensuremath{f_{\mathrm{rep}}}}

\newcommand{\period}[2]{\ensuremath{\mathrm{period}_{#1}(#2)}}
\newcommand{\tpbf}[2]{\ensuremath{\mathfrak{T}_{\mathrm{PBF}}(#1,#2)}}

\newcommand{\orelloop}[2]{\ensuremath{\tau_{#2}^{#1}}}

\newcommand{\aut}{\ensuremath{\mathrm{Aut}}}

\newcommand{\resmap}{\rho}
\newcommand{\quotmap}{\pi}

\newcommand{\into}{\ensuremath{\hookrightarrow}}
\newcommand{\interior}[1]{\ensuremath{{\mathbf{interior}\left(#1\right)}}}

\newcommand{\basepair}[1]{\ensuremath{(#1,\ol{#1})}}
\newcommand{\size}[1]{\ensuremath{\mathrm{Size}(#1)}}
\newcommand{\onto}{\ensuremath{\twoheadrightarrow}}
\newcommand{\immersion}{\looparrowright}
\newcommand{\dual}[1]{\overline{#1}}
\newcommand{\dualtree}[2]{\ensuremath{{T(#1,#2)}}}

\newcommand{\ETone}[1]{{\ensuremath{\mathfrak{T}}(#1)}}
\newcommand{\ADone}[1]{{\ensuremath{\mathfrak{A}}(#1)}}
\newcommand{\ETtwo}[1]{\ensuremath{\mathfrak{T}_+(#1)}}
\newcommand{\ADtwo}[1]{\ensuremath{\mathfrak{A}_+(#1)}}

\newcommand{\trnhd}[2]{\ensuremath{\mathrm{N}_{#1}{\left(#2\right)}}}
\newcommand{\effTracks}[1]{\ensuremath{\mathbf{tracks}{\left(#1\right)}}}
\newcommand{\kappaTracks}[1]{\ensuremath{ \mathbf{tracks}_{\kappa}
    {\left(#1\right)}}} 

\newcommand{\mob}[1]{\ensuremath{{\mathcal{M}(#1)}}}
\newcommand{\branch}[1]{\ensuremath{\mathfrak{b}_{#1}}}

\newcommand{\ETP}[1]{\ensuremath{{\mathfrak{T}^P}_+(#1)}}

\newcommand{\Taux}[1]{\ensuremath{\mathfrak{T}^{\mathrm{aux}}\left(#1\right)}}

\newcommand{\graphofgroups}[1]{\ensuremath{\mathbb{#1}}}

\author{Nicholas W.M. Touikan\\Department of Mathematical Sciences\\
  Stevens Institute of Technology
  \\email:~\texttt{nicholas.touikan@gmail.com}} \title{Detecting
  geometric splittings in finitely presented groups.}
\begin{document}
\maketitle
\abstract{We present an algorithm which given a presentation of a
  group $G$ without 2-torsion, a solution to the word problem with
  respect to this presentation, and an acylindricity constant
  $\kappa$, outputs a collection of tracks in an appropriate
  presentation complex. We give two applications: the first is an
  algorithm which decides if $G$ admits an essential free
  decomposition, the second is an algorithm which, if $G$ is
  relatively hyperbolic, decides if it admits an essential elementary
  splitting.}

\section{Introduction}\label{sec:introduction}

A important group invariant is whether or not it splits as a certain
type of graph of groups. In this paper we prove an algorithmic
analogue of Sela's $\kappa$-acylindrical super accessibility
\cite[Theorem 4.3]{Sela-acylindrical} for the class of one edged
$\kappa$-acylindrical geometric splittings. In particular the main
result, Theorem \ref{thm:main}, gives an algorithm that produces a
list that contains a representative of every one edged
$\kappa$-acylindrical geometric splitting of $\fungrp{\complex C}$, up
to equivalence in $\aut\left(\fungrp{\complex C}\right)$. We give some
corollaries of this theorem.

\begin{mainthm}\label{thm:grushko}
  There is an algorithm that takes as input a finite presentation
  $\bk{X \mid R}$ of a group $G$ without 2-torsion and a solution to
  the word problem with respect to this presentation and decides whether or not
  the group $G$ admits an \emph{essential} free decomposition, i.e. a
  free decomposition \[G=H_1*H_2\] with $H_1 \neq \{1\} \neq H_2$.
\end{mainthm}

This theorem is proved in Section \ref{sec:grushko}. As a consequence
we have the following corollary whose proof we leave as an
exercise. (Hint: if we can solve the word problem, then we can decide if a
finitely generated group is abelian and we can decide, given a finite
presentation, if an abelian group is cyclic.)

\begin{cor}\label{cor:grushko}
  Let $G=\bk{X \mid R}$ be as in the statement of Theorem
  \ref{thm:grushko}, then we can find a Grushko decomposition for $G$.
\end{cor}

In a sense, aside from the no 2-torsion assumption, this is the
strongest result of this type possible: the restrictions on the input
are as minimal as can be reasonably expected. This result also extends
all previously known results (at least in the case without 2-torsion),
which we now briefly survey.

Diao and Feighn in \cite{Diao-Feighn-free-split} showed how to find a
Grushko decomposition of a fundamental group of a graph of free
groups. Their techniques rely on Whitehead methods refined by Gersten
and group actions on square complexes. Kharlampovich and Miasnikov in
\cite{KM-JSJ} showed how to find a Grushko decomposition of a fully
residually free group by running their Elimination Process: the free
decomposition becomes apparent by ``separating the variables'' in the
defining equations. 

Even in the presence of 2-torsion, Dahmani and Groves in
\cite{DG-free-split} are able to detect free splittings of certain
relatively hyperbolic groups, by generalizing an unpublished algorithm
for hyperbolic groups due to Gerasimov. Their approach is to decide
some connectivity criterion of the boundary of toral relatively
hyperbolic groups. Our work implies this result in the 2-torsion-free
case.  Another result \cite{Wilton-Groves-Enumerate} due to Groves and
Wilton, which works in the presence of 2-torsion, is that given a
finite presentation of a group $G$ and a solution to the word problem
with respect to that presentation we can decide if $G$ is free. If $G$
is without 2-torsion, then this is an easy consequence of Corollary
\ref{cor:grushko}. At the end of Section \ref{sec:grushko} we will
explain later how to obtain this result in complete generality from
the work in this paper.

It is also worth noting that Casals-Ruiz and Kazachkov used methods
related to ours to describe solutions to equations over free products
\cite{C-r-K-FP}.

The algorithm given in this paper is also well-suited to relatively
hyperbolic groups we have the following immediate corollary to Theorem
\ref{thm:elementary-splittings}, which is proved in Section \ref{sec:rel-hyp}.

\begin{cor}
  We can decide if a torsion-free relatively hyperbolic group with
  polycyclic parabolics has a trivial elementary JSJ decomposition, in
  the sense of \cite{bowditch1998cut, Bow_periph}or
  \cite[Theorem 4]{GL-trees}.
\end{cor}

This generalizes a result of Dahmani and Groves in
\cite{Dahmani-Groves-iso} for toral relatively hyperbolic groups in
two ways. Firstly, the present approach works for a larger wider class
of groups. Secondly, it can detect splittings that are not in the
class $\mathcal{Z}_{\mathrm{max}}$. Both of these earlier limitations
arise from fact that all previous algorithms to detect splittings in
relatively hyperbolic groups depend on ``equational'' methods. In
particular they will not work with nilpotent parabolics since we can't
solve equations over nilpotent groups \cite{Romankov-undecidable}
and they can't detect non-$\mathcal{Z}_{\mathrm{max}}$ splittings
since Dehn twists around such groups give trivial automorphisms.

In an earlier preprint \emph{Effective Grushko decomposition}
(\url{http://arxiv.org/abs/0906.3902v1}) the author claimed Theorem
\ref{thm:grushko} \emph{without} the no 2-torsion hypothesis. There is
a gap in that proof: the argument is incomplete because the author did
not take Möbius strips into consideration. In the second version of
this paper we modified the argument so that it can handle
$\kappa$-acylindricity and the existence of Möbius strips, at the
cost of having to exclude 2-torsion. Otherwise, there was no gap in
the second version of this paper, but it was horribly written. This
third version attempts to rectify the issue and has more pictures.

\subsection{Acknowledgements}
I wish to thank Olga Kharlampovich for her answers to some very
technical questions, enabling me to adapt the ideas in
\cite{KM-IrredII} to obtain the results of sections
\ref{sec:periodicity-reduction} and \ref{sec:bound-periodicity}. I am
grateful to Alexei Miasnikov for his numerous (ultimately successful)
attempts to explain the elimination process in advanced courses and
for drawing my attention to tracks. I am also grateful to Ilya
Kazachkov, Montserrat Casals-Ruiz, Martin Dunwoody, Gilbert Levitt,
François Dahmani, Vincent Guirardel, Daniel Groves, and Henry Wilton
for discussion and encouragement. The second version of this paper was
written while the author was an NSERC postdoctoral fellow at CIRGET in
UQAM and at the Oxford Mathematical Institute. Although much of the
third version paper was written on the agonizingly slow Amtrak
Adirondack Train that connects New York and Montreal, I feel no
gratitude towards Amtrak. Finally I thank the anonymous referee for
actually reading this paper and giving precise and relevant
feedback. The paper is much better now.

\subsection{Outline of the paper}
First we will give the basic definitions and results needed to make
sense of the statement of Theorem \ref{thm:main}. After stating it and discussing some of its limitation,
we apply it to detect free decompositions of finitely presented groups
and elementary splittings of relatively hyperbolic groups. These
applications also serve the role of providing a ``tutorial'' on how to
use Theorem \ref{thm:main}. 

Next we will define band complexes which are similar, but not
identical, to the band complexes in
\cite{Bestvina-Feighn-1995}. Instead of using them to study minimal
foliations, we will treat them as combinatorial objects to study
Dunwoody patterns (see Definition \ref{defn:pattern}) . Next we will
define transformations done to band complexes and tracks they
carry. This constitutes the Rips machine of
\cite{Bestvina-Feighn-1995}.

The Rips machine is designed to study a single lamination in a cell
complex. The elimination process, inspired from works of Makanin and
Razborov \cite{Makanin-1982,Razborov-1987} as read from
\cite{KM-IrredII}, is a branching search algorithm that constructs a
finite rooted directed tree that decides the existence of certain
types of laminations, in our case, tracks. We will show how to
construct this tree one level at a time and give an analogy with
splitting sequences for surface train tracks.  We will then define
various inadmissibility criteria which will forbid the elimination
tree from growing at certain nodes.

As usual in this business, it will be relatively easy to handle the
thin/Levitt/7-10 case as well as the surface/quadratic/12 case. The
real difficulty is in handling the superquadratic/axial/15 case and this
is where most of the new ideas in this paper reside.

Eventually we will have given sufficiently many inadmissibility
criteria, including a periodicity bound, to force the elimination tree
to be finite. The leaves of this tree will give us the output of the
main algorithm.

\subsection{Patterns, tracks, and geometric splittings}

We take it for granted that the reader is comfortable with Bass-Serre
theory. The best reference, especially for this paper, would be
\cite{Scott-Wall}. Another standard reference is
\cite{Serre-arbres}. We also assume the reader is well acquainted with
polygonal 2-complexes, their fundamental groups, and the actions of fundamental groups on universal covers by deck transformations.

The graphs of groups $\graphofgroups{X}$
has underlying graph $X$. We will write
\define{$G$ splits as a graph of groups $\graphofgroups{X}$} or even
\define{$\graphofgroups{X}$ is a splitting of $G$} instead writing
``$G$ is the fundamental group of the graph of groups
$\graphofgroups{X}$.''  We will also use the action of a group $G$
on a simplicial tree $T$ and the corresponding splitting
$\graphofgroups{X}$, where $X = G\backslash T$,
interchangeably. Finally all trees are assumed to be minimal.

\begin{conv}
  In order to be sure to avoid any pathologies, we will restrict
  ourselves to the piecewise linear category of topological spaces.
\end{conv}

Throughout this paper $\complex{C}$ will be a polygonal 2-complex. If
$f:X \to Y$ is a continuous map we denote its functorial image
$f_{\sharp}: \fungrp X \to \fungrp{Y}$, which is well-defined up to
conjugacy.  Let $Y \subset X$ be connected cell complexes. Consider
the natural map:
\[
\begin{tikzpicture}
  \squareofnodes{Y}{X}{\pi_1(Y)}{\pi_1(X)} \draw[right hook->] (1) --
  node[above]{$i$} (2); \draw[->] (3) -- node[above]{$i_\sharp$} (4);
\end{tikzpicture}
\]
where $i$ denotes the inclusion map. We denote
\[\gp(Y) = i_\sharp(\pi_1(Y)),\]
which gives a well defined conjugacy class in $\pi_1(X)$.

\begin{defn}\label{defn:pattern}Let $\complex{C}$ be a 
  polygonal 2-complex. A \emph{pattern} $P \subset \complex{C}$ is an
  embedded 1-complex such that: \begin{enumerate}[(i)]
  \item for every 2-cell $D \subset \complex{C}$, $P\cap D$ is a
    (possibly empty) \emph{finite} collection of closed arcs joining
    \emph{distinct} sides of $D$;
  \item $P$ does not meet $\complex{C}^{(0)}$;
  \item $P$ has a regular neighbourhood $N(P)\subset \complex{C}$
    homeomorphic to ${P\times [-1,1]}$.
  \end{enumerate}
\end{defn}

\begin{defn}\label{defn:track} A connected component of a pattern is called a
  \emph{track}.
\end{defn}

This definition of a pattern is slightly non-standard in that the last
condition implies that our pattern is \emph{2-sided}. 2-sidedness
however should be standard because it implies that the pattern is
locally separating into two components. The Seifert-van Kampen Theorem
immediately implies that the decomposition (which is essentially a graph
of spaces)
\begin{equation}\label{eqn:cut-along-P} \complex{C} =
  N(P) \cup \complex{C}\setminus P
\end{equation} 
splits $\pi_1(\complex{C})$ as a graph of groups $\graphofgroups{X}^P$
where the vertex groups are given by $\gp(\complex{C}_i)$ and the
edge groups are given by $\gp(t_j)$, where the $\complex{C}_i$ denote the
connected components of $\complex{C} \setminus P$ and the $t_j$ denote
the tracks in $P$ respectively.

\begin{prop}\label{prop:tree}
  Let $P \subset \complex{C}$ be a pattern and let $\tilde{P}$ be the
  lift of $P$ in the universal cover $\tilde{\complex{C}}$ of
  $\complex{C}$. Each connected component of $\tilde{P}$ separates
  $\tilde{\complex{C}}$ into two components. This gives rise to the
  $\pi_1(\complex{C})$-tree $\dualtree{P}{\complex{C}}$ whose vertices are
  connected components of $\tilde{\complex{C}} \setminus \tilde{P}$, whose
  edges are connected components of $\tilde{P}$ and such that the
  edge $\tilde{t}$ is adjacent to the vertex $\tilde{\complex{C}_i}$ if
  $\tilde{t}$ is contained in the closure of of
  $\tilde{\complex{C}_i}$. Thus $\dualtree{P}{\complex C}$ can be
  obtained by a $\fungrp{\complex C}$-equivariant identification map
  \begin{equation}
    \label{eq:phi-map}
    \quotmap: \tilde{\complex C} \onto \dualtree{P}{\complex C}
  \end{equation}

\end{prop}
\begin{proof}
  The lift $\tilde{P} \subset \tilde{\complex{C}}$ of $P$ is again a
  pattern in $\tilde{\complex{C}}$. This gives a decomposition of
  $\tilde{\complex{C}}$ as in (\ref{eqn:cut-along-P}) which expresses
  $\pi_1(\tilde{C})$ as a graph of groups $\graphofgroups{Y}$. If some
  component of $\tilde{P}$ is not separating the underlying graph $Y$
  of $\graphofgroups{Y}$ contains a cycle contradicting the fact that
  $\pi_1(\tilde{C})=1$. It therefore follows that the graph $Y$ is a
  tree $\dualtree{P}{\complex{C}}$, which is easily seen to be a
  $\pi_1(\complex{C})$-tree.

  $\quotmap$ is obtained collapsing each track neighbourhood
  $N(\tilde t) = \tilde t \times[-1,1] \onto [-1,1]$ and collapsing every
  connected component of
  $\closure{\tilde{\complex{C}} \setminus N(\tilde P)}$ to a point.
\end{proof}

This next proposition follows immediately by thinking about the action
of $\pi_1(\complex{C})$ on $T(P,\complex{C})$ induced by deck
transformations and the meaning of the Seifert-van Kampen Theorem, or
simply by thinking of (\ref{eqn:cut-along-P}) as a \emph{graph of
  spaces decomposition} \`a la \cite{Scott-Wall}.

\begin{prop}\label{prop:pattern-tree}
  Let $P \subset \complex{C}$ be a pattern. Then the action of
  $\pi_1(\complex{C})$ on $T(P,\complex{C})$ gives the splitting of
  $\pi_1(\complex{C})$ as the graph of groups $\graphofgroups{X}^P$
  induced by the decomposition (\ref{eqn:cut-along-P}).
\end{prop}

All that being sorted, we can now make sense of the second and third
words of the title of the paper.

\begin{defn}\label{defn:dual-bass-serre}
  For a pattern $P \subset \complex C$ the tree $\dualtree{P}{\complex
  C}$ obtained in Proposition \ref{prop:tree} is called the
\define{Bass-Serre tree dual to $P$}, or simply the \define{dual
  Bass-Serre tree}.
\end{defn}

\begin{defn}
  A splitting of $\pi_1(\complex{C})$ is \define{geometric} if it is
  represented by a pattern $P$, i.e. the Bass-Serre tree of the
  splitting is given by the action of $\pi_1(\complex{C})$ on
  $T(P,\complex{C})$. The pattern $P$ is said to be \define{essential}
  if $T(P,\complex{C})$ is infinite.
\end{defn}

The following fact is important since it implies that the class of
geometric splittings is significant.

\begin{thm}[Restatement of {\cite[Lemma
    2.2]{D-S-JSJ}}]\label{thm:resolution}
  Let $\pi_1(\complex{C})$ act minimally on a tree $T$. Then there
  exists a pattern $P \subset \complex{C}$ such that there is a
  $\pi_1(\C)$-equivariant surjective \emph{simplicial} map called a
  \define{resolution} \[\resmap:
  \dualtree{P}{\complex{C}} \rightarrow T.\] In particular, the edge
  stabilizers of $\dualtree{P}{\complex{C}}$ are conjugate to subgroups of the
  edge stabilizers of $T$. Moreover if the action of
  $\pi_1(\complex{C})$ on $T$ is non-trivial then some track in $P$
  will be essential.
\end{thm}

Immediately we get:

\begin{cor}\label{cor:track-free-decomp}
  If $\pi_1(\complex{C})$ is freely decomposable, then some essential
  free decomposition is geometric.
\end{cor}

\begin{cor}\label{cor:track-jsj}
  If $\complex{C}$ is a finite complex, then any
  Guirardel-Levitt JSJ deformation space \cite[Definition 4]{GL-JSJ-I}
  of $\fungrp{\complex{C}}$ contains a geometric splitting.
\end{cor}
\begin{proof}
  We refer the reader to the introduction of \cite{GL-JSJ-I} for the terminology in this
  proof.  Let $T$ be some JSJ tree for $\fungrp{\complex C}$ over some
  class of groups $\mathcal{A}$, i.e. a domination-maximal universally
  elliptic $\mathcal{A}$-tree. Theorem \ref{thm:resolution} implies
  the existence of a geometric tree $\dualtree{P}{\complex{C}}\to T$ that
  dominates $T$. Since the edge groups of $\dualtree{P}{\complex{C}}$ are
  contained in edge groups of $T$ and $\mathcal{A}$ is assumed to be
  closed under taking subgroups, $T(P,\complex C)$ is also an
  $\mathcal{A}$-tree. Since the edge groups of $T$ are
  $\mathcal{A}$-universally elliptic, so must the edge groups of
  $\dualtree{P}{\complex{C}}$. It follows that
  $\dualtree{P}{\complex{C}}$ is also a domination-maximal universally
  elliptic $\mathcal{A}$-tree so the result follows.
\end{proof}

\subsection{Relative  splittings}\label{sec:relative-splittings}
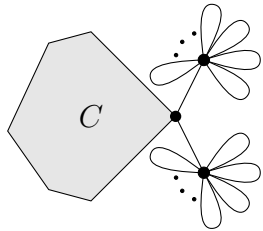
\begin{figure}[htb]
\centering
\begin{tikzpicture}[scale =1.5]
  \draw[fill=black!10!white] (0:0.75) -- (90:0.75) -- (120:0.75) --
  (190:0.75) -- (240:0.75) -- (270:0.75)-- cycle;
  \draw (0.75,0) -- (1,0.5);
  \draw (0.75,0) -- (1,-0.5);
  \drawbouquet{-10}{0.5}{{(1cm,0.5cm)}}
  \begin{scope}[yscale=-1]
    \drawbouquet{-10}{0.5}{{(1cm,0.5cm)}}
  \end{scope}
  \draw (0,0) node {$\complex{C}$};
  \fill (0.75,0) circle (0.05);
\end{tikzpicture}

\caption{Attaching bouquets of circles to a presentation complex. We
  must then attach 2-cells to preserve the fundamental group.}
\label{fig:attaching-bouquets}
\end{figure}

Suppose we are given a finite collection of finite generating
\[ S=\big\{ \{h_i\}_{i\in I_n} \mid n = 1,\ldots, m \big\}\]
of subgroups of $\pi_1(\complex{C})$, and that we want to study the
geometric splittings of $\fungrp{\complex C}$ in which the subgroups
$\bk{h_i}$ are elliptic. Then we can make a new 2-complex
$\complex{C}_S \supset \complex{C}$ with $\pi_1(\complex{C}) \approx \pi_1(\complex{C}_S)$
as follows (see Figure \ref{fig:attaching-bouquets}.) For each
$S_n = \{h_i\}_{i\in I_n}$
\begin{enumerate}[(1)]
\item\label{it:e_h} Make a bouquet of circles $\bouquet{n}$ such that for each $h
  \in S_n$ there is a directed edge $e_h$ in $\bouquet{n}.$
\item\label{it:a_n} Attach the vertex $v_n$ of $\bouquet{n}$ to the vertex $v$ of
  $\complex{C}$ by an arc $\alpha_n$.
\item Attach a 2-cell so that the loop $\alpha_n*e_h*\alpha_n^\mo$ is
  now homotopic to $h \in \pi_1(\complex{C},v)$.
\end{enumerate}
We call the resulting 2-complex $\complex{C}_S$. We note that
$\pi_1(\complex{C}) \approx \pi_1(\complex{C}_s)$ because of the
obvious deformation retraction $\complex{C}_S \onto \complex{C}$.
We now employ the following trick to restrict to relative geometric
splittings.

\begin{prop}\label{prop:relative-splittings}
  Let $\elliptics$ be the set of subgroups generated by the elements
  of $S$. Then every track $t$ dual to a geometric splitting of
  $\pi_1(\complex C)$ relative to $\elliptics$ can be extended to a
  track $t' \subset \complex{C}_S$ so that $t'$ is disjoint from the
  edges $e_h$ in item (\ref{it:e_h}) of the construction of $\complex{C}_S$.
\end{prop}
\begin{proof}[sketch]
  Let $\pi_1(\complex{C}) \times T \rightarrow T$ be a geometric
  action with the subgroups $\elliptics$ acting elliptically then we
  can extend the pattern $P \subset \complex{C} \subset \complex{C}_S$
  to a pattern $P' \subset \complex{C}_S$ such that $P'$ has empty
  intersection with the edges $e_h$ and such that we have a
  $\pi_1(\complex{C})$-equivariant
  isomorphism\begin{equation}\label{eqn:relative-extension}
    T(P',\complex{C}_S) \rightarrow T(P,\complex{C}).\end{equation} We
  do this by taking a resolution \[\resmap': \widetilde{\complex{C}_S}
  \rightarrow T(P,\complex{C})\] which extends
  $\resmap:\widetilde{\complex{C}} \rightarrow T(P,\complex{C})$ such
  that the lifts of the vertices of $v_i$ are mapped to vertices
  stabilized by appropriate conjugates of $\bk{S_i}$. We refer the
  reader to Section 2 of \cite{D-S-JSJ} for details on the resolution
  construction. It therefore follows geometric splittings of
  $\pi_1(\complex{C})$ relative to $\elliptics$ are given exactly by
  patterns in $\pi_1(\complex{C}_S)$ that do not intersect the new
  edges $e_h$.
\end{proof}

\subsection{Equivalence under automorphisms}\label{sec:auto-equiv}

Let 
\begin{eqnarray*} \varphi: G \times T &\to &T \\
  (g,x)& \mapsto& g\cdot x\\
\end{eqnarray*} 
be an action of a group $G$ on a tree $T$. Then for any $\alpha \in
\aut(G)$ we may \define{twist} $\varphi$ by $\alpha$ to get a new
action \begin{eqnarray*}
  \varphi^\alpha:G \times T & \to & T \\
  (g,x) & \mapsto & \alpha(g)\cdot x
\end{eqnarray*}

\begin{defn}\label{defn:aut-equiv}
  Let $\varphi:G\times T \mapsto T$ and $\psi:G\times S \mapsto S$ be
  two actions of the group $G$ on simplicial trees. Let $\alpha \in
  \aut(G)$ . We say the actions $\varphi$ and $\psi$ are
  \define{$\aut(G)$-equivalent}, written $\varphi \sim_{\aut(G)} \psi$, if
  there is a simplicial isomorphism $f: T \to S$ and an $\alpha \in
  \aut(G)$ that makes the the following diagram commutative: \[
  \begin{tikzpicture}
    \squareofnodes{G\times T}{T}{G \times S}{S}
    \draw[->] (1) -- node[above]{$\varphi^\alpha$} (2);
    \draw[->] (3) -- node[above]{$\psi$} (4);
    \draw[->] (1) -- node[left]{$1 \times f$} (3);
    \draw[->] (2) -- node[right]{$f$} (4);
  \end{tikzpicture}
  \]
\end{defn}

\begin{defn}\label{defn:pattern-equiv}
  Let $P,P'$ be patterns in $\complex{C}$. We say that \define{ the
    patterns $P$ and $P'$ are $\aut(\pi_1(\complex{C}))$-equivalent},
  written $P\sim_{\aut(\pi_1(\complex{C}))} P'$ if the natural actions
  $\pi_1(\complex{C}) \times T(P,\complex{C}) \to T(P,\complex{C})$
  and $\pi_1(\complex{C}) \times T(P',\complex{C}) \to
  T(P',\complex{C})$ are $\aut(\pi_1(\complex{C}))$-equivalent.
\end{defn}

In the case of closed surfaces patterns arise as multicurves and
automorphic equivalence of multicurves coincides with equivalence
under homeomorphisms. For general 2-complexes (which may have trivial
mapping class groups) these equivalences may not coincide.

\subsection{The main result}\label{sec:main-result}
A solution to the word problem in $\pi_1(\complex{C})$ is a procedure
that decides if a loop (given as a sequence of directed edges in
$\complex{C}^{(1)}$) is nullhomotopic in $\complex C$. In the case
where $\complex C$ is a presentation complex such a loop uniquely
defines a word in the prescribed generating set.

\begin{mainthm}\label{thm:main}
  There is an algorithm which takes as input a finite 2-complex
  $\complex{C}$ such that $\pi_1(\complex{C})$ has no 2-torsion, a
  solution to the word problem in $\pi_1(\complex{C})$, some positive
  integer $\kappa$, a finite collection
  \[ S=\big\{ \{h_i\}_{i\in I_n} \mid n = 1,\ldots, m \big\}\]
  of finite generating sets of subgroups
  $\elliptics = \big\{ \bk{h_i}_{i\in I_n}\big\}_{n=1}^m$ of
  $\pi_1(\complex{C})$ and outputs a finite collection of tracks
  $t_1,\ldots,t_{n(\complex{C},\kappa,S)}$ which lie in a complex
  $\complex{C}_S \supseteq \complex{C}$ (with equality if
  $S = \emptyset$) such that the isomorphism
  $\pi_1(\complex{C}_S) \approx \pi_1(\complex{C})$ is given
  explicitly. These tracks give splittings of $\pi_1(\complex{C})$
  relative to $\elliptics$ with the following property: if
  $\pi_1(\complex{C})$ admits a geometric $\kappa$-acylindrical
  splitting relative to $\elliptics$ represented by a track $t$ then
  there is some $i \in \{1,\ldots,n(\complex{C},\kappa,S)\}$ such that
  $t \sim_{\aut(\pi_1(\complex{C}))} t_i$.
\end{mainthm}

This theorem is proved in Section \ref{sec:main-proof}, where the main
algorithm is given. The
2-complex $\complex{C}_s$ was defined in Section
\ref{sec:relative-splittings}. 

The result, as stated, is about splittings that are geometric,
$\kappa$-acylindrical and with one edge group. This result is not the
strongest possible, but it gives us all the applications we need.

It could be strengthened as follows: by Theorem \ref{thm:resolution} every $\kappa$-acylindrical tree $T$
is resolved by a geometric tree
$\resmap: \dualtree{P}{\complex C} \onto T$. The resolving tree
$\dualtree{P}{\complex C}$, however, may not itself be
$\kappa$-acylindrical. The analysis of the relationship between the
trees $\dualtree{P}{\complex C}$ and $T$ in
\cite{delzant1999accessibilite} combined with the arguments of this
paper actually give a finite collection of tracks that \emph{resolve}
every $\kappa$-acylindrical tree. The geometric resolving splittings
themselves may not be $\kappa$-acylindrical, but they are ``locally''
$\kappa$-acylindrical in a way that is good enough for our
arguments. Although this would give us a full algorithmic version of
Sela's super accessibility \cite[Theorem 4.3]{Sela-acylindrical} for
the class of $\kappa$-acylindrical one edged splittings, we have opted
for a simpler formulation; thus removing a layer of notation. We hope
the reader will agree that this is for the best.

For the sake of simplicity we have also restricted ourselves to
one-edged splitting or tracks, instead of general patterns. This does
not weaken the result because we allow relative splittings which
enables us to produce refinements. Results such as
\cite{Weidmann-2002,delzant1999accessibilite} then give explicit
bounds on the number of components of the pattern.

It should also be noted that Theorem \ref{thm:main} does not
necessarily enable us to detect whether $\fungrp{\complex\C}$ actually
has a geometric $\kappa$-acylindrical splitting. To reach such a
conclusion we must be able to further analyze the collection of tracks
produced by the main algorithm. This means we must be able to solve
more delicate algorithmic problems in the ambient group. The next two
applications, especially the proof of Theorem
\ref{thm:elementary-splittings}, will illustrate the necessary extra
requirements.

Finally there is the issue of torsion. The current algorithm cannot
handle actions on trees with arbitrarily long arcs with non-trivial
pointwise stabilizers, even if these stabilizers are
finite. Forbidding 2-torsion, for example, controls a problem that occur
with Möbius bands by bounding their width, which gives terminating
conditions. If the algorithm were to run in the presence of 2-torsion,
then it would still produce a (possibly empty) list of tracks and
terminate. However, because the algorithm will have stopped
prematurely this list may be missing some tracks.

There is no reason these torsion issues cannot be overcome. For
example, \cite{Dahmani-Guirardel-foliations} deals with laminations in
band complexes with torsion. The author suspects to solve this problem
one would have to generalize band complexes to some version for
2-orbihedra with finite cell stabilizers.

\subsection{Computing Grushko decompositions}\label{sec:grushko}

This is an application of Theorem \ref{thm:main} with $\elliptics = \emptyset$ and $\kappa =
0$. In this case we only need to be able to solve the word problem.

\begin{proof}[Proof of Theorem \ref{thm:grushko}]
  Let $\complex{C}$ be a presentation 2-complex form $\bk{X \mid R}$
  and consider a maximal splitting of $\bk{X \mid R}$ over finite
  groups. If $G$ admits an essential free decomposition then by
  Corollary \ref{cor:track-free-decomp} there is a track $t \subset
  \complex{C}$ that represents this splitting. 

  Free decompositions correspond exactly to 0-acylindrical
  actions on trees. We now apply the algorithm of Theorem
  \ref{thm:main} to get a finite collection of
  tracks \[t_1,\ldots,t_n.\] If there is a track $t$ that
  represents an essential free splitting of $G$ then $t$ is
  $\aut(G)$-equivalent to some $t_i$ in our finite
  collection. So $G$ admits an essential free decomposition if and
  only if some $t_i$ represents an essential free decomposition.

  With our solution to the word problem we are able to check for each
  $t_i$ if \begin{itemize}
  \item $\gp(t_i) = \{1\}$,
  \item
    \begin{itemize}
    \item if $\complex{C} \setminus t_i$  is not connected, then both
      components of must have non trivial image in $\fungrp C$
      via the inclusion map, or
    \item if $\complex{C} \setminus t_i$ is connected, there is
      nothing to show,
    \end{itemize}
  \end{itemize}
  and thus decide if $t_i$ represents an essential free
  decomposition.
\end{proof}

We now give another method to decide if a finitely presented group $G$
with decidable word problem is free (see
\cite{Wilton-Groves-Enumerate}.) First note that if we can solve the
word problem, we can decide if a finitely presented group is abelian
(check if the generators commute) and then, by linear algebra, we can
compute its isomorphism type. In particular we can decide if $G$ is
isomorphic to $\mathbb Z$. We repeatedly apply the algorithm for
Theorem \ref{thm:grushko} to attempt to compute the Grushko
decomposition of $G$. If $G$ is 2-torsion-free then this will be the
correct Grushko decomposition. Otherwise we still will obtain some
(possibly trivial) free decomposition of $G$. We can then decide if
each factor of this decomposition is isomorphic to $\mathbb Z$ or
$\{1\}$. This will be the case if and only if $G$ is free.

\subsection{Detecting splittings of relatively hyperbolic groups}
\label{sec:rel-hyp}
For this section we assume that the reader is familiar with relatively
hyperbolic groups. The reader can consult \cite{Farb-rel-hyp} or
\cite{Hruska} for definitions.  Let $G$ be a finitely presented
\emph{torsion-free} group that is hyperbolic relative to the finitely
generated subgroups $\elliptics = \{H_1,\ldots,H_m\}$. We assume that
the groups $H_i \in \elliptics$ are pairwise distinct and
non-conjugate. Before continuing we need to give some definitions.

\begin{defn}
  An element of $g \in G$ (respectively a subgroup $K \leq G$) is
  \emph{parabolic} if there exists some $h \in G$ such that $h^\mo gh
  \in H_i$ (respectively $h^\mo K h \leq H_i$) for some $i \in
  \elliptics$.
\end{defn}

\begin{conv}
  We will assume in this section that all algorithms in a group are
  with respect to a presentation and that the (tuples of) elements of
  the input are given as (tuples of) words in the symmetrized
  generating set.
\end{conv}

\begin{defn}\label{defn:elementary}
  A splitting of $G$ is \define{elementary} if all parabolic subgroups
  are elliptic and the edge groups are either trivial, infinite
  cyclic, or parabolic.
\end{defn}

\begin{defn}\label{defn:tractable-triple}
  A triple $\big(\bk{S \mid R}, \texttt{CP}, \texttt{Gen} \big)$
  where:\begin{enumerate}[(i)]
  \item $\bk{S \mid R}$ is a finite group presentation,
  \item $\texttt{CP}$ is an algorithm which solves the the conjugacy
  problem with respect to $\bk{S \mid R}$, and 
\item\label{it:gen-pb}$\texttt{Gen}$ is an algorithm
  which decides whether or not a finite tuple generates $\bk{S \mid
    R}$
\end{enumerate} is called an \define{algorithmically tractable
  triple}.
\end{defn}

\begin{defn}\label{defn:tractable-class}
  A class $\mathcal{C}$ of finitely presented groups is called
  an \emph{algorithmically tractable class of parabolics} if there is
  an algorithm which enumerates algorithmically tractable triples
  corresponding to the groups in $\mathcal{C}$.
\end{defn}

It is worth pointing out that by \cite{BCRS-PF}  the class of \emph{polycyclic-by-finite}
groups is algorithmically tractable.
We now collect some well known facts about torsion-free relatively hyperbolic
groups.

\begin{prop}[{\cite[Example 1 p.819]{Farb-rel-hyp}}]\label{prop:rel-hyp-malnormal}
  For all $h \in G$, $\left(h^\mo H_i h\right) \cap H_j \neq \{1\}$ if and only if
  $i=j$ and $h \in H_i$.
\end{prop}

This next result about elements of $G$ follows from the work in
\cite{Bumagin-rel-hyp-cp} but is stated explicitly in 
\cite{Osin-rel-hyp}. The generalization to explicitly
given subgroups of $G$ follows applying Proposition \ref{prop:rel-hyp-malnormal}.

\begin{thm}[c.f. {\cite[Theorem 5.6]{Osin-rel-hyp}}]\label{thm:parabolicity-problem}
  Given $g\in G$ (respectively $K=\bk{k_1,\ldots,k_n} \leq G$) if we are given a solution to the
  conjugacy problem for each $H_i, i=1,\ldots,n$ then we can decide
  whether there is some $h \in G$ such that $h^\mo g h \in H_i$
  (respectively $h^\mo Kh \leq H_i$) for some $i \in \{1,\ldots,n\}$, and
  find $h$ if it exists.
\end{thm}

The following two facts are well known, however I couldn't find any
precise references, they are stated here and proved.

\begin{prop}\label{prop:maximal-cyclics}
  Let $g \in G$ be a non-parabolic element. Then its centralizer
  $C(g)$ is infinite cyclic and malnormal.
\end{prop}
\begin{proof}
  By Theorems 4.16 and Theorem 4.19, the centralizer $C(g)$ of $g$ is
  a word hyperbolic group. By \cite[Corollary
  3.6]{Notes-on-hyp-groups}, since $G$ is torsion free $C(g)$ is
  infinite cyclic. Assume now for simplicity that $C(g) = \bk{g}$.

  Suppose there was some $h \in G$ such that
  $h^\mo \bk{g} h \cap \bk{g} \neq \{1\}$ then by \cite[Corollary
  4.26]{Osin-rel-hyp} there is some $l \in \mathbb{Z}_{\neq 0}$ such
  that $h^\mo g^l h = g^{\pm l}$. This means that $h^2 \in C(g^l)$ and
  that $\bk{g^l}$ is normal in $\bk{h,g^l}$. Now as explained before
  $\bk{k}=C(\bk{g^l}) \geq C(g) = \bk{g}$, which implies that
  $k \in C(g)$, so $k \in \bk{g}$; thus $h^2 \in \bk{g}$. From this we
  get that $[\bk{h,g^l}:\bk{g^l}] \leq 2l$.

  Now \cite[Lemma 11.4]{Hempel-3-manifolds} states that if a group $Q$
  contains a an infinite cyclic subgroup of finite index, then $Q$
  contains a finite subgroup $K$, such that $Q/K$ is either isomorphic
  to $\mathbb Z$ or $\mathbb Z_2 * \mathbb Z_2$. If $Q \leq G$, then
  $Q$ is torsion-free; so $K$ must be trivial and $Q$ must be infinite
  cyclic. It follows that $\bk{h,g^l}$ is infinite cyclic, say
  $\bk{h,g^l} = \bk{z}$. Then $z \in C(g)$, so in particular
  $h \in C(g)$.
\end{proof}

\begin{cor}\label{cor:2-acylindrical}
  A one edged elementary splitting of a torsion free relatively
  hyperbolic group $G$ is 2-acylindrical.
\end{cor}
\begin{proof}
  We first prove the following. \emph{Claim: let $e$ be an edge in
    Bass-Serre $T$ tree connecting the vertices $u,v$. Then at least
    one of the images of $G_e \leq G_u$ or $G_e \leq G_v$ is
    malnormal.}

  Indeed, if the splitting in question is free, then the result holds.
  Suppose now that the edge group is $\bk{g}$ for some non parabolic
  $g \in G$. Then $\bk{g}$ must be maximal cyclic in at least one of
  its images in the vertex groups, otherwise its centralizer is not
  cyclic because in the amalgam $\bk{x}*_{x^r=y^s}\bk{y}$ if $r,s \neq
  1$ then $\bk{xy,y^s} \approx \Z \oplus \Z \leq C(y^s)$. Also the
  images of $\bk{g}$ in the associated subgroups cannot intersect
  since by \cite[Corollary 4.27]{Osin-rel-hyp} any Baumslag-Solitar
  group must be parabolic.
  
  Suppose now that the edge group $G_e$ is parabolic, but not maximal
  parabolic, and hence malnormal, in either $G_u,G_v$. Then we have
  parabolic proper overgroups $G_e < P_u \leq G_u$ and $G_e < P_v \leq
  G_v$. One one hand since $|P_u \cap P_v| = |G_e| = \infty$ they must
  lie in a common maximal parabolic subgroup $P$. On the other hand
  $P$ does not act elliptically on $T$ (it has a non-trivial induced
  splitting) contradicting the fact that the splitting is
  elementary. \emph{This proves our claim.}

  We now prove 2-acylindricity. Let $u$ be a vertex in the Bass-Serre
  tree and let $e,f$ be edges such that $e \cap f = \{u\}$. Since
  $e,f$ are in the same $G$-orbit we have that $G_e,G_f$ are
  conjugate. On the other hand if $e \neq f$ then if $G_v \geq G_e
  \cap G_f \neq \{1\}$ then $G_e$ is not malnormal in $G_v$.  Suppose
  towards a contradiction that there is some $g \in G_e$ such for some
  edge $h \subset T$ such that $e \cap h =\emptyset$ we have
  $g\cdot h = h$, and suppose moreover that there is some edge $f
  \subset T$ such that $e \cap f =\{u\}$ and $f \cap h = \{v\}$. Then
  we must have that $g \in G_e \cap G_f \cap G_h$. $g \in G_e \cap
  G_f$ implies that $G_f$ is not malnormal in $G_u$ which means by our
  earlier claim that $G_f$ must be malnormal in $G_v$ so $G_f
  \cap G_h =\{1\}$ contradiction. Therefore no element of $G \setminus
  \{1\}$ fixes an arc of $T$ of length more than 2.
\end{proof}

We finally need the following.

\begin{thm}[Theorem 3 of
  \cite{Dahmani-Guirardel-Presenting-parabolics}]\label{thm:find-parabolics}
  There exists an algorithm as follows. It takes an input of a finite
  presentation of a group $G$, a solution to its word problem, and a
  recursive class of finitely presented groups $\mathcal{C}$ (given by
  a Turing machine enumerating presentations of these groups).

  It terminates if and only if $G$ is properly hyperbolic relative to
  subgroups that are in the class $\mathcal{C}$.

  In this case, the algorithm outputs an [relative linear]
  isoperimetry constant $K$ [in the sense of \cite[Definition
  2.30]{Osin-rel-hyp}], a generating set and a finite presentation for
  each of the parabolic subgroups.
\end{thm}

Now we have our second application.

\begin{mainthm}\label{thm:elementary-splittings}
  Suppose we are given a finite presentation $\bk{X\mid R}$ of a
  torsion free group $G$ that is relatively hyperbolic with finitely
  many parabolics that lie in an algorithmically tractable class of
  parabolics $\mathcal{C}$. Suppose also that we are given a solution
  to the word problem with respect to $\bk{X\mid R}$ and a finite collection
  $S$ of finite generating sets for a set of subgroups
  $\elliptics'$. Then we can decide if $G$ admits an elementary
  splittings relative to $\elliptics'$.
\end{mainthm}
\begin{proof}
  We first note that by Theorem \ref{thm:resolution} if $\complex{C}$
  is the presentation 2-complex associated to $\bk{X \mid R}$, then
  $G$ admits an essential elementary splitting if and only if
  $\pi_1(\complex{C})$ admits an essential elementary \emph{geometric}
  splitting. Any elementary splitting is, by Corollary
  \ref{cor:2-acylindrical}, 2-acylindrical.

  We first run the algorithm of Theorem \ref{thm:find-parabolics} to
  find the finite collection $\{H_1,\ldots,H_n\}$ (given by generating
  sets in $\bk{X\mid R}$) of parabolic subgroups we then apply the
  algorithm of Theorem \ref{thm:main} with $\kappa = 2$ and
  \[\elliptics = \elliptics' \cup \{H_1,\ldots,H_n\}\] with the
  collection of generating sets $S$. This gives us a finite collection
  of tracks $t_1,\ldots,t_{n(\complex{C},\kappa,S)}$ that lie in
  $\complex{C}_S$. It is now enough to check for each of these tracks
  if they represent an essential elementary splitting. Let $t$ be one
  of these tracks. By Theorem \ref{thm:parabolicity-problem} we can
  decide if $\gp(t)$ is parabolic.

  If $\gp(t)$ is trivial then as in the proof of Theorem
  \ref{thm:grushko} we can decide if it gives an essential
  splitting. 
  
  Suppose now that $\gp(t)$ isn't parabolic. Since we can solve the word problem we can
  check whether $\gp(t)$ is abelian, if it isn't then it certainly
  cannot represent an elementary splitting. Otherwise $\gp(t)$ is
  abelian. By Proposition \ref{prop:maximal-cyclics}, $\gp(t)$ is
  contained in the centralizer of some non-parabolic element and is
  therefore contained in a non-parabolic maximal cyclic group. If $t$
  is a non-separating track then it gives an essential elementary
  splitting. Otherwise $t$ separates $\complex{C}_S$, and $G$ splits
  as a free product with amalgamation over $\gp(t)$. To check if the
  splitting is essential it suffices to check, using the word problem,
  whether the generators of the vertex groups commute with
  $\gp(t)$. Indeed, since we are assuming that $\gp(t)$ is
  non-parabolic, we can assume that the vertex groups are
  non-parabolic; so by Proposition \ref{prop:maximal-cyclics} if one
  of the vertex groups commutes with $\gp(t)$ then it is at most a
  finite index cyclic overgroup of $\gp(t)$. Deciding if the vertex
  group coincides with $\gp(t)$ can now be solved using item 3) of
  Theorem 1.16 of \cite{Osin-rel-hyp}.

  Suppose finally that $\gp(t)$ is parabolic. Again, if $t$ is
  non-separating the splitting is essential. Otherwise the splitting
  is essential if and only if $\gp(t)$ doesn't equal one of the vertex
  groups. If neither of the vertex groups are parabolic then the
  splitting is essential. Otherwise at most one of the vertex groups
  is parabolic and we can decide if it is generated by $\gp(t)$ using
  Theorem \ref{thm:parabolicity-problem} and our solution to the
  generation problem given by the algorithmic tractability assumption
  (Definition \ref{defn:tractable-triple} (\ref{it:gen-pb})).
\end{proof}

\section{Band complexes}\label{sec:band-complexes}

The algorithm of Theorem \ref{thm:main} is a procedure that will
produce a rooted directed tree (i.e. a branching sequence) of band
complexes.  Band complexes first appeared in
\cite{Bestvina-Feighn-1995} to classify stable actions of finitely
presented groups on $\mathbb{R}$-trees. Our version of band complexes
differ in that they are combinatorial objects: they do not come with
laminations, instead we will allow a band complex to \emph{carry}
multiple laminations or, in our case, tracks.

As combinatorial objects, our band complexes will contain the same
amount of information as Makanin's generalized equations
(c.f. \cite{KM-IrredII}.)

\subsection{Definitions and terminology}

\begin{defn}\label{defn:measured-band}
  A \define{band} $\Band$ is a Cartesian product ${J_\Band \times
    [-1,1]}$ where $J_\Band$ homeomorphic to a closed interval. The
  subsets $J_\Band \times \{\pm 1\}$ are called \define{bases}. 
  If $\mu = J_\Band \times {\{\pm 1\}}$ is the base of a band then we
  call the base $\ol{\mu} = J_\Band \times {\{\mp 1\}}$ the
  \define{dual} of $\mu$.
\end{defn}

A band is therefore a rectangle with well-defined bases and a vertical
direction.

\begin{conv}
  The letters $\lambda, \mu,\eta,\nu$ shall be used to denote bases
  and $\ol{\lambda}$ will always denote the dual of $\lambda$. We
  shall denote by $\band{\lambda}$ the band that contains $\lambda$.
\end{conv}

\begin{defn}\label{defn:MBC}
  A \define{band complex} $\mbc$ is a 2-complex that is constructed
  in the following way:
  \begin{enumerate}[(1)]
  \item\label{it:underlying graph} Start with a simplicial graph $\Gamma$.
  \item Attach the bases of the bands $\Band_1,\ldots
    \Band_m$ to the interiors of edges of $\Gamma$ via
    embeddings \[g_{\mu_i}:\mu_i \into \Gamma \setminus \Gamma^{(0)}\]
    where the $\{\mu_i\}$ is the set of bases of the bands.
  \item\label{2-cell-requirements} Let \[U = \big(\Gamma \cup \Band_1
    \cup \ldots \cup \Band_m\big)/\sim\] be the resulting
    identification space. We finally obtain $\mbc$ by attaching discs
    $D_1,\ldots,D_l$ via immersions $f_i: \bdy D_i \immersion U$ with
    the following requirement:
    \begin{enumerate}
    \item\label{it:conn-a} For all $i,j$, $f_i(\bdy D_i) \cap \Band_j$
      can be expressed as a finite union of embedded arcs $\alpha_i$
      that travel from one base of $\Band_j$ to the other. Such arcs
      are called \define{connections}.
    \item Connections are pairwise disjoint.
    \item\label{it:conn-c} If a connection has non-trivial
      intersection with a side of a band $\Band_j$, then it coincides
      with that side.
    \end{enumerate}
  \end{enumerate}
\end{defn}

In the case of measured band complexes \cite{Bestvina-Feighn-1995} the
complicated requirement \ref{2-cell-requirements} above on the 2-cell
attaching map is ensured if the 2-cell attaching maps intersect
measured bands in \emph{vertical} subsets.

Connections (as described in items (\ref{it:conn-a})-(\ref{it:conn-c})
of the definition above) will occur exactly where \emph{boundary
  connections} occur when working with generalized equations (see for
example \cite{KM-IrredII}.) Controlling their cardinality is a key
step in the repetition argument which deals with the thinning and
superquadratic cases of the elimination process (Section
\ref{sec:repetition}.)

\begin{conv}
  Although a band complex is a 2-complex, whenever we mention a 2-cell
  we \emph{really mean a 2-cell $D_i$ that gets attached in step 3. of
  Definition \ref{defn:MBC}.}
\end{conv}

Thus, it is possible that a 2-cell $D_i$ in a band complex never
intersect any bands, in which case the image of the attaching map
$f_i(\bdy D_i)$ lies entirely in the underlying graph $\Gamma$ in item
(\ref{it:underlying graph}) of Definition \ref{defn:MBC}. 

\begin{conv}\label{conv:bases}
  Formally speaking, a base $\mu$ isn't a subset of the band complex
  $\mbc$. That being said we will still write $x \in \mu$ for some
  point $x \in \mbc$ such that $x \in g_\mu(\mu)$. We will also write
  $\lambda \subset \mu$ if $g_\lambda(\lambda) \subset g_\mu(\mu)$. In
  the case where $g_\lambda(\lambda) = g_\mu(\mu)$ we will use the
  evocative symbol $\lambda \bandeq \mu$ to avoid confusion. We will
  also treat the bands $\Band_i$ as subsets of $\mbc$ when it is convenient.
\end{conv}

\begin{defn}\label{defn:matched-base}
  We say that bases $\mu, \ol{\mu}$ are \define{matched bases} if
  $\mu \bandeq \ol{\mu}$ and $\band\mu$ forms an annulus in
  $\mbc$. Otherwise a base is called \emph{unmatched}.
\end{defn}

\begin{defn}\label{defn:maximal-section}
  A union of unmatched bases $U = \bigcup_{\mu \in S} \mu$ is called
  \emph{strongly connected} if the union of the interior of the bases
  $U' = \bigcup_{\mu \in S} \interior{\mu}$ is also connected (and
  therefore an interval.) A maximal (with respect to inclusion)
  strongly connected union of unmatched bases is called a
  \define{maximal section}.
\end{defn}

Maximal sections are almost the blocks in \cite{Bestvina-Feighn-1995}
and the closed sections in \cite{KM-IrredII}.

\begin{defn}[Carrying a track]\label{def:carry-track}
  Let $\mbc$ be a band complex and let $t \subset \mbc$ be a
  track. $\mbc$ \define{carries} $t$ if $t$ is contained in the union
  of the bands in $\mbc$ and furthermore,
  \begin{enumerate}[(i)]
  \item For each band $\Band_j$, $t\cap\Band_j$ consists of a union of
    pairwise disjoint embedded arcs travelling from one base of
    $\Band_j$ to the other, and
  \item $\Band_j \cap t$ is disjoint from the connections in $\Band_j$
    as well as from its sides.
  \end{enumerate}
\end{defn}
\begin{defn}[Efficiently carrying]\label{def:efficiently-carry}
  For a base $\mu$, let $S_\mu \subset \mu$ be the finite set
  containing the points the form $x=\mu \cap c$, where $c$ is a
  connection, and the points $x$ that are the endpoints of bases.
  $\mbc$ \define{carries $t$ efficiently} if for every base $\mu$ and
  every distinct $x,y \in S_\mu$ there is some point in $t\cap\mu$
  that separates them.
 \end{defn}

 Thus, if a $\mbc$ carries a track $t$, then $t$ is confined to the
 interior of the bands.  The notion of carrying naturally generalizes
 to arbitrary measured laminations, but since we will only be focusing
 on one leaf laminations we only need to deal with the hitting
 measure.

 \begin{defn}[Measure from a track]\label{defn:measure}
   Let $\mbc$ be a band complex, let $t$ be a track carried by $\mbc$, and
   let $S \subset \mbc$ be a union of bases. We define the
   \define{hitting measure on $S$ with respect to $t$}, denoted $|S|_t$, to be
   the cardinality of the intersection\[ |S|_t=|S\cap t|.
   \] If $\mu$ is a base of $\mbc$ then we will sometimes call
   $|\mu|_t$ the \define{length of $\mu$ with respect to $t$}.
\end{defn}

\subsection{Constructing (measured) band complexes from
  tracks}\label{sec:correspondences}

Let $\complex C$ be a standard CW 2-complex and let
$t \subset \complex C$ be a track. We obtain a band complex $\mbc$
from $\complex C$ as follows. 

\begin{figure}[htp]
  \centering
    \begin{tikzpicture}
  \begin{scope}
    \draw[clip] (-30:2) -- (90:2) -- (210:2) --cycle;
    \draw[thick,dashed] (-2,1.15) -- (2,1.15) (-2,1) -- (2,1);
    \draw[thick,dashed](-0.2,-1.73) -- (1.8,1.73);
    \draw[thick,dashed] (-2.4,1.73) -- (-0.4,-1.73) (-2,1.73) -- (0,-1.73) (-2.2,1.73) -- (-0.2,-1.73);

  \end{scope}
  \draw[very thick] (-30:2) -- (90:2) -- (210:2) --cycle;
  \begin{scope}[shift={(4,0)}]
    \begin{scope}
    \draw[clip] (-30:2) -- (90:2) -- (210:2) --cycle;
    \fill[black!30!white] (0,-1) -- ++(120:3) -- ++(-1.5,0) --
    ++(-60:3) -- cycle;
    \draw (0,-1) -- ++(120:3) -- ++(-1.5,0) --
    ++(-60:3) -- cycle;
    \fill[black!30!white] (0,-1) -- ++(60:3) -- ++(0.5,0) --
    ++(-120:3) --cycle;
    \draw (0,-1) -- ++(60:3) -- ++(0.5,0) --
    ++(-120:3) --cycle;
    \fill[black!30!white] (-2,0.5) -- (2,0.5) -- ++ (60:1) --
    ++(-4,0) --cycle;
    \draw (-2,0.5) -- (2,0.5) -- ++ (60:1) --
    ++(-4,0) --cycle;

    \draw[thick,dashed] (-2,1.15) -- (2,1.15) (-2,1) -- (2,1);
    \draw[thick,dashed](-0.2,-1.73) -- (1.8,1.73);
    \draw[thick,dashed] (-2.4,1.73) -- (-0.4,-1.73) (-2,1.73) -- (0,-1.73) (-2.2,1.73) -- (-0.2,-1.73);
      \end{scope}
    \draw[decorate,decoration={brace,amplitude=0.1cm}] (-1,0.5)
    ++(-120:1.5) -- node[shift={(-0.2,0.2)}]{3} (-1,0.5);
    \draw[decorate,decoration={brace,amplitude=0.1cm}] (-1,0.5)--
    node[shift={(-0.2,0.2)}]{2} ++(60:1);
    \draw[decorate,decoration={brace,amplitude=0.1cm}] (1,0.5)
    --node[shift={(0.2,0.2)}]{1}   ++(-60:0.5);
    \draw[very thick] (-30:2) -- (90:2) -- (210:2) --cycle;
  \end{scope} 
\end{tikzpicture}
 
\caption{On the left, the intersection of a track $t$ with a 2-cell
  $U$ in a 2-complex. On the right how to construct the corresponding
  measured band complex efficiently carrying $t$ by dividing $U$ into
  three bands and four 2-cells.}
  \label{fig:correspondence}
\end{figure}
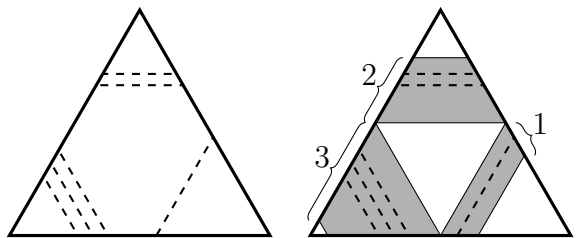

For each 2-cell $D \subset \complex{C}$, $t \cap D$ is a disjoint
union of arcs travelling from one edge of $\partial D$ to another edge
of $\partial D$. Metrize the 1-skeleton $\complex{C}^{(1)}$, giving
each edge $e_i$ a length of $|e_i \cap t|+2$. Subdivide each 2-cell
$D$ into a union of bands and 2-cells, such that a base $\mu$ has
length $|\mu \cap t|=|\mu|_t$. Explicitly parameterize each band
\[\band\mu = \left[0,|\mu|_t\right]\times[-1,1],\]
so that $t \cap \band\mu \cap t$ is a union of vertical sets as
follows:
\[ t \cap \band\mu = \bigcup_{i=1}^{|\mu|_t} \{i-\frac 1 2\} \times
[-1,1];
\]
this is illustrated in Figure \ref{fig:correspondence}.

This construction gives a band complex as in Definition
\ref{defn:MBC}. The underlying graph is the 1-skeleton
$\complex{C}^{(1)}$, the bases of the bands are embedded in this graph
and avoid the vertices $\complex{C}^{(0)}$. We further see that the
remaining 2-cells have embedded (thus, immersed) boundaries that only
intersect bands in their vertical sides.  Furthermore, since no 2-cell
attaching maps go through the interior of any of the bands, $t$ is
efficiently carried by the band complex $\mbc$.

\begin{defn}\label{defn:measured}
  A band complex $\mbc$ is \define{measured} if every base $\mu$ is
  metrized as a real closed interval $\left[0,|\mu|_t\right]$ (recall
  Definition \ref{defn:measure}). A subset $v$ of a band
  $\band\mu = \left[0,|\mu|_t\right] \times [-1,1]$ is called
  \define{vertical} if it is of the form \[ v = \{x_v\} \times [-1,1].
  \]
\end{defn}

We have proved:

\begin{prop}\label{prop:construct-bc}
  For any track $t$ contained in a 2-complex $\complex C$, we can
  subdivide $\complex C$ into a measured band complex $\mbc$ in which
  $t$ consists of a union of vertical sets. Moreover $\mbc$ carries
  $t$ efficiently.
\end{prop}

\begin{defn}[Combinatorial equivalence of band complexes]\label{defn:BC-equiv}
  Two band complexes $\mbc$ and $\mbc'$ are said to be
  \define{equivalent} if there is a homeomorphism $\mbc \to \mbc'$
  that sends the underlying simplicial graph (Definition
  \ref{defn:MBC} (\ref{it:underlying graph})) to the underlying
  simplicial graph, sends bands to bands, sends 2-cells to
  2-cells, and for each of these objects restricts to a homeomorphism.
\end{defn}

If we forget the measures on band complexes, we are left with only
finitely many possibilities; thus,

\begin{prop}\label{prop:finite-ubcs}
  Let $\complex{C}$ be a finite 2-complex and let $S$ be a finite
  collection of finite subsets of $\fungrp{\complex C}$. Then there
  are only finitely many possible band complexes (up to the
  combinatorial equivalence of Definition \ref{defn:BC-equiv}) that
  arise from the possibly infinite collection of tracks
  $t \subset \complex{C}_S$ (Proposition
  \ref{prop:relative-splittings}.)  Furthermore this list can be
  effectively constructed.
\end{prop}

\section{Moves on band complexes carrying
  tracks}\label{sec:Rips-machine}

We will present moves that transform a band complex carrying a track
into a new band complex carrying a new track.  These moves are
essentially the moves given in \cite[\S
6.1]{Bestvina-Feighn-1995}. Our treatment is slightly different since
we want to explicitly realize each move as a continuous map
$m:\mbc \to \mbc'$ sending the track $t$ to some track
$t' \subset \mbc'$, which we denote $(\mbc,t) \to (\mbc',t')$. This is
accomplished with the zipping moves. They will be of use in a later
section.

We will employ the convention of \cite{KM-IrredII} and reuse the names of bases, as is customary in
computer science.

\subsection{The basic moves on band complexes that carry a track}\label{sec:basic-moves}
Let $\mbc$ be a band complex efficiently carrying a track $t$. Suppose
furthermore that $\mbc$ is measured (Definition \ref{defn:measured}) so
that a base $\mu$ is metrized with length $|\mu|_t$ and $t \subset
\mbc$ is a union of vertical sets. 

We first define elementary moves $(\mbc,t) \to (\mbc',t')$ which
transform the underlying band complex and track, while preserving the
fundamental group and dual Bass-Serre tree $T(t,\mbc)$. These moves
are actually $\pi_1$-isomorphic continuous maps $\mbc \to \mbc'$ that
map $t$ to $t'$.

\begin{defn}[Type I zip]\label{move:zip}
  Suppose we have a containment of bases $\lambda \subset \mu$ with
  $\dual \lambda \cap \dual\mu= \emptyset$ and
  $\dual \lambda \cap \mu = \emptyset$. The union of $U$ vertical sets
  of $\band\mu$ that intersect $\lambda$ is a rectangle homeomorphic
  to $\band\lambda$.  A \define{type I zip of $\band\lambda$ into
    $\ \band\mu$} consists of the operation of identifying
  $\band\lambda$ to $U$ so that vertical sets are sent
  homeomorphically to vertical sets and $\dual\lambda$ is identified
  to the corresponding subset of $\dual\mu$.
\end{defn}

\begin{defn}[Type II zip, or squish]\label{move:squish}
  Suppose we have the containments of bases $\lambda \subset \mu$ and
  $\dual\lambda \subset \dual\mu.$ Suppose furthermore that there are
  vertical paths $\alpha \subset \band\lambda$ and
  $\beta \subset \band\mu$ such that the concatenation $\alpha*\beta$
  is a nullhomotopic loop.  A \define{type II zip of $\band\lambda$
    into $\band\mu$} is the operation of continuously identifying
  $\alpha$ to $\beta$ and continuously extending this to an
  identification of $\band{\lambda}$ to a union of vertical subsets of
  $\band\mu$. This identification map
  must be injective when restricted to $\band\mu$ and $\band\lambda$
  and must send vertical sets to vertical sets.
\end{defn}

\begin{figure}[htb]
  \centering
  \begin{tikzpicture}[scale=0.3]
  \draw[thick,fill=black!20!white] (0,0)--(2,4)--(10,4)--(8,0)--cycle;
  \draw[thick,fill=black!10!white] (0.5,1)--(1.5,3) -- (9,7)--(8,5)--cycle;
  \draw (11,2) node{$\to$};
  \draw[dashed](2,4)--(10,4);
  \begin{scope}[xshift=12cm]
     \draw[thick,fill=black!20!white] (0,0)--(2,4)--(10,4)--(8,0)--cycle;
  \draw[thick,fill=black!10!white] (0.5,1)--(1.5,3) -- (9.5,3)--(8.5,1)--cycle;
  \end{scope}
  \draw (23,2 )node{$\leftarrow$};
  \begin{scope}[xshift=24cm]
    \draw[ultra thick,fill=black!20!white] (0,0) --(2,4) .. controls +(1,-1) and
    +(-1,-1) .. (10,4)--(8,0) .. controls +(-1,-1) and
    +(1,-1) .. (0,0);
    \draw[thick,fill=black!10!white](0.5,1)--(1.5,3) .. controls +(1,1) and
    +(-1,1) .. (9.5,3)--(8.5,1) .. controls +(-1,1) and
    +(1,1) ..(0.5,1);
    \draw[dashed](2,4) .. controls +(1,-1) and
    +(-1,-1) .. (10,4) ;
  \end{scope}
\end{tikzpicture}
  \caption{Type I and Type II zips. The identifications must send
    tracks to tracks.}
  \label{fig:zips}
\end{figure}
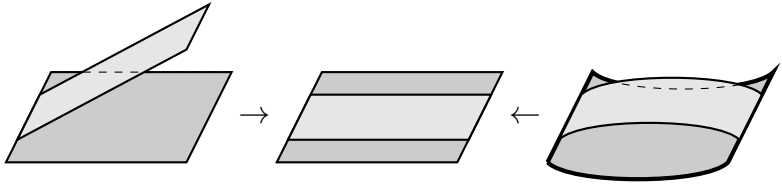

These zipping moves are not used in \cite{Bestvina-Feighn-1995}, but
the type I zip is a step in the transfer (see Definition
\ref{move:transfer} later) or \emph{M4 Slide} of
\cite{Bestvina-Feighn-1995}. The zipping moves will be necessary later
when we will be ``wrapping one band around another'' (see Figure
\ref{fig:merge-zip}.) Also the fact that they are given by explicit
continuous maps is convenient.

\begin{defn}[Collapse a band]\label{move:collapse}
  Let $\band{\mu}$ be a band such that
  $\mu \cap \dual \mu = \emptyset$. Then the \define{collapse of
    $\band{\mu}$ onto $\mu$} is the operation of identifying each
  vertical subset $v$ that intersects $\interior\mu$ to the point
  $x_v = v\cap\mu$.
\end{defn}

\begin{defn}[Annulus]
  An \define{annulus} $\Ann$ is a band $\band\mu$ such that $\mu \bandeq
  \dual\mu$ whose image in the band complex $\mbc$ is homeomorphic to
  an annulus.
\end{defn}

This next move is similar to the type II zip.

\begin{defn}[Crush an annulus]\label{move:crush}
  Let $\Ann \subset \mbc$ be an annulus such that $\gp(\Ann)=\{1\}$,
  then \define{crushing $\Ann$} is the operation of identifying each
  vertical subset of $\Ann$ to its intersection with $\mu$.
\end{defn}

\begin{defn}[Vertically subdivide a band]\label{move:subdivide}
  Let $\Band = [a,b] \times [-1,1]$ be a band in $\mbc$ and let $\{ p\}
  \times [-1,1]$ be a vertical subset. \define{The operation of subdividing
    $\Band$ along $\{ p\} \times [-1,1]$}. Consists of the following
  operations:
  \begin{enumerate}[(1)]
  \item Cut $\Band$ along $\{ p\} \times [-1,1]$ so that we get two bands
    $\Band_1 = [a,p^-] \times [-1,1]$ and $\Band_2 = [p^+,b] \times
    [-1,1]$.
  \item Attach a 2-cell along the loop $(p^-\times [-1,1]) * (p^+ \times
    [-1,1])$, where $*$ denotes concatenation. This 2-cell is called a
    \define{subdivision digon}.
  \end{enumerate}
\end{defn}
\begin{conv}
  We only allow band subdivision of $(t,\mbc)$ if the resulting $t'
  \subset \mbc'$ is efficiently carried. 
\end{conv}
These basic
operations may leave some messiness behind:
\begin{enumerate}[(i)]
\item After performing a zipping move a 2-cell may no longer have an
immersed boundary.
\item \label{it:collapse} After a collapse a 2-cell may have a free face and
perhaps the resulting band complex can be given as $\mbc =
\mbc'*_p\alpha$, i.e. the connected sum at a point $p$ of a band
complex $\mbc'$ and a closed arc $\alpha$.
\item\label{it:sphere} After crushing an annulus, or performing a type II zip, the
  boundary of a 2-cell may map onto a point or an interval resulting in
  a sphere.
\end{enumerate}

We therefore introduce, as basic moves, the following cleaning operations.

\begin{defn}[Delete superfluous cells]\label{move:clean}
  If a 2-cell in $\mbc$ is a sphere, as may occur in
  (\ref{it:sphere}) above, or if it has a free face, remove it. Do the same
  for hanging arcs that occur in (\ref{it:collapse}) above.
\end{defn}

\begin{defn}[Tighten 2-cells]\label{move:tighten}
  If a 2-cell $D$ no longer has an immersed boundary then the attaching
  map $f_D:\partial D \to \mbc$ factors
  as \begin{equation}\label{eqn:tighten-2cell}\partial D \onto
    \left(S^1*_{p_i}\tau_i\right) \stackrel{f'}{\immersion} \mbc,\end{equation}
  where the middle term is a circle with some hanging trees $\tau_i$,
  which arise from the ``pinching'' of the attaching map. This middle
  term is immersed into $\mbc$. We replace $D$ by a 2-cell $D'$ whose
  boundary is identified with $S^1$ in (\ref{eqn:tighten-2cell}) and
  mapped to $\mbc\setminus D$ (abusing notation) via the immersion
  $f'$ in (\ref{eqn:tighten-2cell}). $D'$ is called the \define{tightening} of
  $D$.
\end{defn}

\begin{figure}[htb]
  \centering
  \begin{tikzpicture}[scale=0.3]
  \filldraw[black!10!white,draw=black,thick]
  (-1,-1) -- (4,4) -- (12,4) -- (7,-1) -- cycle;
  \filldraw[black!10!white,draw=black,thick]
  (0,0) -- (2,2) -- (2,8) -- (0,6) -- cycle;
  \draw[dotted] (0.5,0.5) +(1,7) -- +(1,1);
  \draw[ultra thick](-0.1,0.1) -- node[above] {$\mu$} ++(2.1,2.1);
  \draw[ultra thick](-1.1,-1) -- node[right] {$\lambda$} ++(5,5);
  \draw[ultra thick](7.1,-1) -- node[right] {$\overline{\lambda}$}
  ++(5,5);
  \draw[ultra thick] (0,6) --node[left]{$\ol{\mu}$} (2,8);
  
  \begin{scope}[shift={(13,0)}]
    \filldraw[black!10!white,draw=black,thick]
    (-1,-1) -- (4,4) -- (12,4) -- (7,-1) -- cycle;
    \filldraw[black!10!white,draw=black,thick]
    (0,0) -- (2,2) -- (2,8) -- (0,6) -- cycle;
    \draw[ultra thick](0,0.1) -- node[above] {$\mu_-$} ++(2,2);
    \draw[ultra thick] (0,3.1) -- ++ (2,2);
    \draw[ultra thick] (0,6) --node[left]{$\ol{\mu_+}$} (2,8);
    \draw[dotted] (0.5,0.5) +(1,7) -- +(1,1);
    \draw[ultra thick](-1.1,-1) --  ++(5,5);
    \draw[ultra thick](7.1,-1) --  ++(5,5);
    
    \filldraw[black!20!white,draw=black] (0,0) -- (2,2)
    -- (10,2) -- (8,0) --cycle;
  \end{scope}

  \begin{scope}[shift={(26,0)}]
    \filldraw[black!10!white,draw=black,thick]
    (-1,-1) --node[left,black]{$\lambda$} (4,4) -- (12,4) -- (7,-1) -- cycle;
    -- (10,2) -- (8,0) --cycle;
    \filldraw[black!10!white,draw=black,thick]
    (8,0) -- (10,2) -- (2,8) -- (0,6) -- cycle;
    \draw[ultra thick] (0,6) --node[left]{$\ol{\mu}$} (2,8);
    \draw[ultra thick] (8,0) --node[right]{$\mu$} (10,2);
    \draw[dotted] (0.5,0.5) +(1,7) -- +(9,1) -- +(1,1);
  \end{scope}
  
\end{tikzpicture}  
\caption{Transferring $\mu$ from $\lambda$ to
  $\ol{\lambda}$ across $\band\lambda$. The dotted line shows a
  connection (Definition \ref{defn:MBC}(\ref{it:conn-a})). After a transfer this connection may give rise to two
  connections.}\label{fig:transfer}
\end{figure}
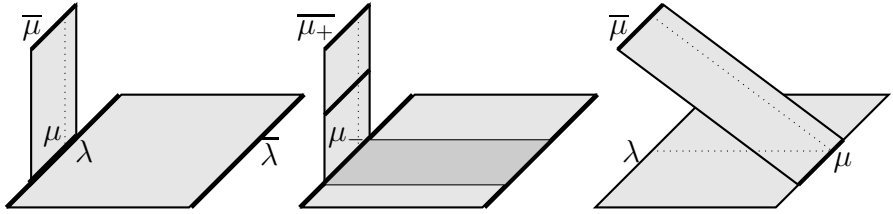

We leave it to the reader to verify that the result of a 2-cell
removal and the tightening move $\mbc \to \mbc'$ can be realized by a
continuous map. The transfer, given below and illustrated in Figure
\ref{fig:transfer} is defined in terms of band subdivisions and
zipping, but we will also treat it as an elementary move.

\begin{defn}[The transfer]\label{move:transfer}
  Let $\mu \subset \lambda$ be bases such that $\mu \neq
  \ol{\lambda}$. The operation of \define{transferring $\mu$ from
    $\lambda$ to $\ol{\lambda}$ across $\band\lambda$} is the
  following sequence of operations.
  \begin{enumerate}[(1)]
  \item Subdivide the band $\band\mu$ horizontally:\[
    \band\mu = \band{\mu_-} \cup \band{\mu_+}
    \] with $\mu \bandeq \mu_-$, $\ol{\mu_-} \bandeq \mu_+$, and
    $\ol{\mu_+} \bandeq \ol{\mu}$.
  \item Zip the band $\band{\mu_-}$ into $\band\lambda$. (By
    hypothesis, this is a type I zip.)
  \item We rename the base $\mu_+$ as $\mu$.
  \end{enumerate}
\end{defn}

\subsection{The preservation property}

We will now give a preservation result for our moves. This result is
stated as a fact at the beginning of \cite[\S
6]{Bestvina-Feighn-1995}. In order to lay out the terminology that is
necessary for our purposes, we will carefully state and prove the
preservation property.

\begin{prop}[The preservation property]\label{prop-preservation}
  Let $m:(\mbc,t) \to (\mbc',t')$ be one of the basic moves given in
  Section \ref{sec:basic-moves}. Then we have an isomorphism of
  fundamental groups\[
  m_{\sharp}:\pi_1(\mbc) \stackrel{\sim}{\to} \fungrp{\mbc'}
  \]
  and a simplicial isomorphism of dual Bass-Serre trees
  $m_T:\dualtree t \mbc \to \dualtree{t'}{\mbc'}$ induced by
  $m$. Furthermore this map is $m_\sharp$-equivariant in the following
  sense, letting $\pi_1(\mbc)$ act naturally on $\dualtree t \mbc$ by
  deck transformations via the quotient map $\quotmap$ (Proposition \ref{prop:tree}((\ref{eq:phi-map})), we
  have\[ m_{\sharp}(g)\cdot m_T(x) = m_T(g \cdot x), \]
  for all $g \in \pi_1(\mbc)$ and all
  $x \in \dualtree t \mbc$.
\end{prop}
\begin{proof}
  We first show prove the proposition for zipping, collapsing and
  crushing moves. 

  We first show that the fundamental groups are isomorphic.  Consider
  first either a type I Zip (Definition \ref{move:zip}) of
  $\band\lambda$ into $\band\mu$ (i.e. with $\lambda \subset \mu$) or
  the collapse of $\band\mu$ onto $\mu$ (Definition
  \ref{move:collapse}). These moves lift to $\fungrp{\mbc}$-equivariant
  moves on on $\tilde\mbc$. Pick a basepoint
  $x\in \mu\subset \mbc$, by the disjointness criteria we see that no
  distinct lifts of $x$ in $\tilde\mbc$ are identified and that
  the resulting complex $\tilde{\mbc'}$ remains simply
  connected. Since the lifts of $x$ in $\tilde{\mbc}$ are in bijective
  correspondence with $\fungrp\mbc$, the isomorphism
  $\fungrp{\mbc} \approx \fungrp{\mbc'}$ follows.

  In the case of a type II zip (Definition \ref{move:squish}) or an
  annulus crush (Definition \ref{move:crush}) the $\pi_1$-triviality
  criteria ensure that we can find lifts
  $\tilde{\band\mu}, \tilde{\band\lambda}$ of $\band\lambda, \band\mu$
  (respectively) such that
  $\tilde{\band\mu} \cup \tilde{\band\lambda}$ is as in the right side
  of Figure \ref{fig:zips}, or that the annulus $\ann\mu$ 
  lifts to $\tilde\mbc$. Arguing as before (taking a basepoint in
  $\mu$) we obtain the isomorphism
  $\fungrp{\mbc} \approx \fungrp{\mbc'}$.

  In all cases the moves map bands to (interiors) of bands, and
  vertical sets to vertical sets. Also points in the complement of the
  union of bands of $\mbc$ are sent to the complement of the union of
  band of $\mbc'$, and the restriction to the complement is
  injective. Lifting to $\tilde\mbc$ we therefore see a bijective
  correspondence between the connected components of $\tilde t$ and
  $\tilde{t'}$. Furthermore if connected components $t_1,t_2 \subset \tilde t$ are dual to
  edges that share a vertex in $\dualtree t \mbc$, then their images
  $t_1',t_2'$ will be dual to edges in $\dualtree{t'}{\mbc'}$ that
  share a vertex. The isomorphism of Bass-Serre trees follows, and
  $m_\sharp$-equivariance of the isomorphism follows from construction.

  The proof for band subdivisions, superfluous 2-cell deletions, and
  tightenings is obvious.
\end{proof}

\subsection{Derived moves}

Having defined basic moves we shall now define the composite, or
derived, moves that constitute the Rips machine. We first introduce
the $\tau$-complexity (originally Makanin's $\xi$-complexity \cite{Makanin-1982}), which is one of the main tools of our
analysis. As we define the derived moves we will show
why they do not increase this $\tau$-complexity.

\begin{defn}\label{defn:sec-complexity}
  Let $\sigma \subset \mbc$ be a maximal section (Definition
  \ref{defn:maximal-section}), let $b(\sigma)$ be the number of
  unmatched bases contained in $\sigma$. We define the
  $\tau$-complexity of a section to be
  \[ \tau(\sigma) = \max \big(b(\sigma)-2,0\big)
  \]
\end{defn}

\begin{defn}[{\cite[Definition 4.3]{Bestvina-Feighn-1995},
\cite[\S 5]{KM-IrredII}} ]\label{defn:tau-complexity}
Let $J \subset \mbc$ be a union of maximal sections, then we define
the \define{$J$-relative $\tau$-complexity} to
be \begin{equation}\label{eqn:tau-defn}\tau(\mbc,J) = \sum_{\sigma
    \not\subset J}\tau(\sigma).\end{equation} If $J=\emptyset$ write
$\tau(\mbc)$ instead of $\tau(\mbc,\emptyset)$.
\end{defn}

\begin{defn}\label{defn:gamma}
  For a point $x \in \mu$ we denote by $\gamma(x)$ the number of unmatched
  bases $\lambda$ such that $\lambda \ni x.$
\end{defn}

\begin{defn}\label{defn:vert-length}
  The \define{vertical length} of the attaching map
  $\partial D \immersion \mbc$ of a 2-cell $D$ is the number of
  connected components of the preimages of the connections ((\ref{it:conn-a}) of
  Definition \ref{defn:measured-band}.) Equivalently, this is the
  number of times the attaching map travels through a band.
\end{defn}

\subsubsection{The Möbius move}
\begin{defn} A dual pair $(\mu,\ol{\mu})$ such that $\mu \bandeq
  \dual\mu$ and $\band\mu$ forms a Möbius band is called a
  \define{Möbius pair}.
\end{defn}

This next move is described in \cite[Lemma
6.4]{Bestvina-Feighn-1995}. Note that since we require tracks to be
two-sided, a track $t \subset \mbc$ can never intersect the core of a
Möbius band. It follows that we can always subdivide a Möbius band
along its core and the resulting band complex will still efficiently
carry $t$.

\begin{defn}[The Möbius move]
  Given a Möbius pair $(\mu,\ol{\mu})$, we subdivide the band
  $\band{\mu} = \mu \times [-1,1]$ along $\{m\}\times [-1,1],$ the core of the
  Möbius band. Call the resulting bands $\band{\mu_0}$ and
  $\band{\mu_1}$. We then transfer the base $\mu_0$ across the band
  $\band{\mu_1}$. The dual pair $(\mu_0, \ol{\mu_0})$ now forms an
  annulus $\ann{\mu_0}$ and the pair $\mu_1, \ol{\mu_1}$ intersect at
  a point. If the annulus $\ann{\mu_0}$ is $\pi_1$-trivial, we crush
  it. We rename $\mu_1,\ol{\mu_1}$ as $\mu, \ol{\mu}$
  respectively.
\end{defn}

Straightforward verification gives the following result:

\begin{lem}\label{lem:post-mobius}
  Let $(\mbc,t)$ be a band complex efficiently carrying a track. Let
  $(\mu,\ol{\mu})$ form a Möbius pair. After the corresponding
  Möbius move $m:(\mbc,t)\to (\mbc,t')$ we have the following:
  \begin{enumerate}[(i)]
  \item The hitting measure of $\mu$ drops by a half,
    i.e. \[|\mu|_{t'} = \frac 1 2 |\mu|_t.\]
  \item The hitting measure of $|\mu|_{t'} $ is the
    same as the measures of the matched base $|\mu_0|_{t'}$.
  \item The $\tau$-complexity did not increase.
  \end{enumerate}
\end{lem}

\subsubsection{The thinning move}
\label{sec:thinning}
The thinning move is applied whenever there is some point $x$ in a
base $\mu$ such that $\gamma(x)=1$ (Definition \ref{defn:gamma}). In
\cite{Bestvina-Feighn-1995} this is the move for Process I, but
modified to keep the number of unmatched bases non-increasing (we
prevent the ``long bands'' of \cite{Bestvina-Feighn-1995} from occurring.)
In \cite{KM-IrredII} this is the move applied in cases 7-10.

\begin{defn}\label{defn:maximal-naked-seg}
  Let $x \in \mu$ be a point such that $\gamma(x)=1$. Let
  $x \subset \mu_0 \subset \mu$ be the maximal segment that contains
  only points $x'$ with $\gamma(x')=1$.  We call $\mu_0$ a
  \define{maximal naked segment}.
\end{defn}

Suppose that a base $\mu$ contains a naked segment we now describe the
\define{thinning move starting at $\mu$}:
\begin{enumerate}[(1)]
\item\label{it:th-sub-mu} \textit{Subdivide $\mu$.} subdivide $\mu$ into segments
  $$\mu = \mu_{-1} \cup \mu_0 \cup \mu_{1}$$ with $\mu_\mo$ or $\mu_1$
  possibly empty, and $\mu_0$ a maximal naked segment of $\mu$. This
  results in a subdivision of $\band\mu$. Denote by $D_{-1}$ and
  $D_{1}$ be the added subdivision digons.
\item\label{it:th-sub-ann} \textit{Subdivide annuli and clean $\mu_0$.}  For every matched
  base pair $\basepair{\lambda}$ such that $\lambda$ intersects
  $\interior{\mu_0}$  vertically subdivide
  \[\band{\lambda} = \band{\lambda_{-1}}\cup\band{\lambda_0} \cup \band{\lambda_1}\]
  so that $\lambda_0 \subset \mu_0$ and $\lambda_{\pm 1}$ doesn't intersect
  $\interior{\mu_0}$.

  Next,  take all the resulting base pair $\basepair{\lambda_0}$
  with $\lambda_0\subset \mu_0$ and transfer $\lambda_0$ and
  $\ol{\lambda_0}$ from $\mu_0$ to $\ol{\mu_0}$ through
  $\band{\mu_0}$. Now $\mu_0$ doesn't intersect any other bases.
\item\label{it:th-coll}  \textit{Collapse the naked segment and the added subdivision
    digons.} Collapse the band $\band{\mu_0}$ onto
  $\dual{\mu_0}$. Delete the subdivision digons $D_{-1},D_1$ that were
  added in Step 1 since they now have free faces.
\item\label{it:th-remove} \textit{Remove long bands.} A \define{long
    band} is a union of two bands \[\band\mu \cup \band\lambda\]
  with $\dual\mu \bandeq \lambda$. If a long band is created, first
  transfer $\dual{\mu}$, and all other bases contained in $\lambda$,
  from $\lambda$ to $\dual\lambda$ through $\band\lambda$ and then
  collapse the band $\band\lambda$ onto $\ol\lambda$.
\item\textit{Clean up.} Crush any remaining $\pi_1$-trivial annuli.
\end{enumerate}

A proof of this next fact for generalized equations can be found in
\cite{KM-IrredII}, instead of adapting it we simply give another
proof.

\begin{lem}\label{lem:post-thinning}
  After applying a thinning move starting at a base $\mu$, the
  $\tau$-complexity did not increase. If no annuli were subdivided,
  the number of 2-cells did not increase, nor did the vertical lengths
  of 2-cell attaching maps.
\end{lem}
\begin{proof}
  From (\ref{it:th-coll}), the number of 2-cells did not increase, if no annuli
  were subdivided. Furthermore, since the only transfer move is
  immediately followed by a collapse of the transfer band, the
  vertical lengths of 2-cell attaching maps could not
  increase.

  It remains to show that the $\tau$-complexity doesn't
  increase. Suppose that we performed (\ref{it:th-sub-ann}), so that $\mu_0$
  doesn't intersect any other bases (we could also do
  (\ref{it:th-sub-mu}), (\ref{it:th-sub-ann}) and then recombine the
  bands created in (\ref{it:th-sub-mu}).) Then  the number
  of maximal sections may have increased, but the number of unmatched
  bases did not; thus the $\tau$-complexity did not increase. 

  Let $\sigma$ be the maximal section containing $\mu$ and $\sigma'$
  the maximal segment containing $\ol\mu$, after the subdivision of
  matched bases.

  Suppose first that $\mu$ is completely naked (i.e. $\mu_0 =
  \mu$).
  Then (\ref{it:th-sub-mu}) doesn't occur and after the collapse of
  $\band\mu$ in (\ref{it:th-coll}), the total number of unmatched
  bases goes down by 2. It therefore clearly follows from Definition
  \ref{defn:tau-complexity}(\ref{eqn:tau-defn}) that the
  $\tau$-complexity decreased.

  Suppose now that $\mu$ gets subdivided into 2 segments
  $\mu_0,\mu_1$. Take $\sigma$ to be co-initial with
  $\mu_0$. After (\ref{it:th-sub-mu}) and
  (\ref{it:th-sub-ann}) $\sigma$ gets subdivided into
  $\sigma_0,\sigma_1$ with $b(\sigma_0)=1$ and
  $b(\sigma_1)=b(\sigma)$. In $\sigma'$ the base $\dual{\mu}$ gets
  replaced by $\dual{\mu_0} \cup \dual{\mu_1}$ this increases
  $b(\sigma')$ by 1
  (Definition \ref{defn:tau-complexity}(\ref{eqn:tau-defn})), but in
  (\ref{it:th-coll}) we collapse $\band{\mu_0}$ onto $\dual{\mu_0}$,
  this deletes $\dual{\mu_0}$ from $\sigma'$ so $b(\sigma')$ goes back
  down.

  Suppose finally that $\mu$ gets subdivided into
  $\mu_{-1},\mu_0, \mu_1$. After (\ref{it:th-sub-mu}) the maximal
  section $\sigma$ gets split into
  $\sigma_{-1},\sigma_0,\sigma_1$. With $b(\sigma_0) = 1$ and
  $b(\sigma_{\pm 1}) \geq 2$. If we look at the contribution of what
  is left of $\sigma$ we have a decrease in the contribution of
  $\tau$-complexity of at
  least\[ \tau(\sigma) -
  \big(\tau(\sigma_{-1})+\tau(\sigma_0)+\tau(\sigma_{1})\big) \geq 1.
  \] On the other hand if we look at $\sigma'$ we see that $\ol{\mu}$
  gets subdivided into three bases and $\ol{\mu_0}$ gets deleted after
  the collapse in( \ref{it:th-coll}) we therefore have an increase in the
  contribution to the $\tau$ complexity of resulting sections that
  constitute $\sigma'$, which may have been subdivided, is at most
  $\tau(\sigma')+1$, thus the total $\tau$-complexity did not
  increase.

  Note further that in all the cases above, if we were working with a
  $J$-relative complexity, with $\dual\mu \subset J$ then the
  $J$-relative $\tau$ complexity also did not increase.
  
\end{proof}

\begin{figure}[htb]
\centering
\begin{tikzpicture}[scale={0.6}]
   \filldraw[black!10!white,draw=black]
  (-1,0) -- (-1,-3) -- (3,-3) --(3,0) --cycle
  (-1,0) -- (-2,2) -- (-1,2) -- (0,0) -- cycle
  (2,0) -- (3.5,0) -- (5.5,2) -- (4,2) -- cycle 
  ;
  \draw (1,-1) node {$B(\mu)$};
  \filldraw[black!10!white,draw=black] (0,0) .. controls (0.25,0) 
  and (0.25,1) .. (0,1) -- (0,1) .. controls (-0.25,1) and (-0.25,0)
  .. (0,0);
  
  \filldraw[black!20!white,draw=black] 
  (0,0) .. controls (0.25,0) and (0.25,1) .. (0,1) -- 
  (2,1) .. controls (2.25,1) and (2.25,0) .. (2,0) -- cycle;
  
  \draw[dotted] 
  (2,1) .. controls (1.75,1) and (1.75,0) .. (2,0)
  (-1,2) -- (0,0);
  
  \begin{scope}[shift = {(8,0)}]
    \filldraw[black!10!white,draw=black]
  (-1,0) -- (-1,-3) -- (0,-3) -- (0,0) --cycle
  (2,0) -- (2,-3) --  (3,-3) --(3,0) --cycle
  (-1,0) -- (-2,2) -- (-1,2) -- (0,0) -- cycle
  (2,0) -- (3.5,0) -- (5.5,2) -- (4,2) -- cycle 
  ;
  \begin{scope}[shift={(0,-3)}]
    \filldraw[black!10!white,draw=black] (0,0) .. controls (0.25,0) 
    and (0.25,1) .. (0,1) -- (0,1) .. controls (-0.25,1) and (-0.25,0)
    .. (0,0);
    
    \filldraw[black!20!white,draw=black] 
    (0,0) .. controls (0.25,0) and (0.25,1) .. (0,1) -- 
    (2,1) .. controls (2.25,1) and (2.25,0) .. (2,0) -- cycle;
    
    \draw[dotted] 
    (2,1) .. controls (1.75,1) and (1.75,0) .. (2,0)
    (0,0) -- (0,1);
    
    \draw[thick, dotted] (0,0) -- (0,3);
  \end{scope}
  \draw (-2.5,-1) node {$B(\mu_\mo)$}
  (4,-1) node {$B(\mu_1)$};  
\end{scope}
\end{tikzpicture}
\caption{A thinning move (the 2-cells are not shown)}\label{fig:thinning}
\end{figure}
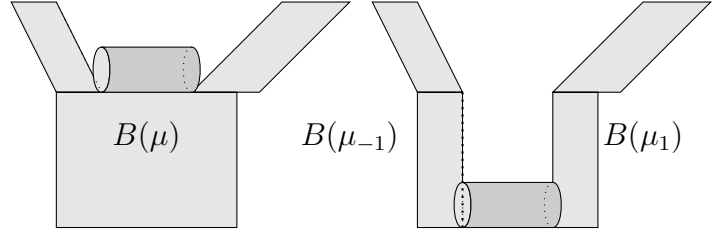

\subsubsection{The entire transformation}
\label{sec:entire-trans}
This   is   the   entire  transformation   given   before   \cite[Case
12]{KM-IrredII}, it also  constitutes the move used in  Process II in
\cite{Bestvina-Feighn-1995}.  One  of   the   disadvantages  of   the
topological setting is that dealing  with closures of maximal sections
is awkward,  due to  the fact that  endpoints of bases  can lie  in the
closures of distinct  maximal sections.

\begin{defn}\label{defn-ordering}
  An \define{ordering $<$ on a band complex $\mbc$} is an ordering $<$
  on the union $U$ of maximal sections (Definition
  \ref{defn:maximal-section}) that is compatible with some embedding
  $U \into \mathbb{R}$. An unmatched base $\mu$ whose interior is
  $<$-coinitial is called a \define{leading base.} An endpoint of
  $\mu$ is called initial (terminal) if it is the limit of a
  $<$-decreasing ($<$-increasing) sequence of points in
  $\interior\mu$.
\end{defn}

\begin{conv}\label{conv:ordering}
  Band complexes will always be assumed to be equipped with an
  ordering.
\end{conv}

We now describe the \define{entire transformation with carrier
  $\mu$}. Let $\mu$ be a maximal leading base.
\begin{enumerate}[(1)]
\item\label{it:et-sub}\textit{Subdivide matched bases.} If the base
  $\mu$ intersects any annuli $\ann{\lambda_i}$, we vertically subdivide them
  so that the resulting annuli are either contained in $\mu$ or do not
  intersect the interior of $\mu$.
\item\label{it:et-move}\textit{Move bases to the right.} Transfer
  every other leading base $\lambda \subset \mu$ (except $\mu$) onto
  $\ol\mu$ through $\band{\mu}$.
\item\label{it:et-collapse} \textit{Collapse the naked initial
    segment.}  Let $\mu_0$ be closure of the $<$-coinitial maximal
  naked subsegment of $\mu$. Subdivide $\band\mu$ into
  $\band{\mu_0} \cup \band{\mu_1}$. Denote by $D_1$ the added
  subdivision digon. Collapse $\band{\mu_0}$ onto $\ol{\mu_0}$ and
  delete $D_1$ because it has a free face.
\item\label{it:et-clean}\textit{Rename and clean up.} We rename $\mu_1$
as $\mu$ and and crush any $\pi_1$-trivial annuli.
\end{enumerate}

\begin{defn}\label{defn:carrier}
  The leading base $\mu$ given in the definition of the entire
  transformation is called the \define{carrier base.}
\end{defn}

An illustration of the result of an entire transformation is given in
Figure \ref{fig:different-outcomes}. This next result follows from a
counting argument.

\begin{lem}[{\cite[Proposition 7.5]{Bestvina-Feighn-1995}}]
  \label{lem:entire-trans}
  After applying an entire transformation the $\tau$-complexity did
  not increase.
\end{lem}

Unlike in the thinning case, the vertical lengths of 2-cell attaching
maps may increase.

\subsection{The Rips machine}
\label{sec:Rips-machine}
The \define{Rips machine}, which was first described in
\cite{Bestvina-Feighn-1995}, is a geometric adaptation of Makanin's
algorithm which takes an measured band complex $\mbc$ and produces a
sequence called the \define{Rips sequence}:
\[\mbc=\mbc_0 \rightarrow \mbc_1 \rightarrow \ldots\]
of measured band complexes, constructed inductively. We assume that
$\mbc$ is equipped with an order $<$, as given in Definition \ref{defn-ordering}.

\begin{lem}\label{lem:order-preserve}
  Let $\mbc$ be equipped with an order $<$ and let
  $m:(\mbc,t) \to (\mbc',t')$ be a Möbius move, a thinning move or an
  entire transformation. Then there is a natural induced order $<'$ on $\mbc'$.
\end{lem}
\begin{proof}
  A Möbius move may split a maximal section $\sigma$ (i.e. delete a
  point) into two maximal sections $\sigma_0,\sigma_1$ $<$ therefore
  restricts to $<'$ on the new union of maximal sections.

  Denote by $U,U'$ the union of maximal sections in $\mbc,\mbc'$,
  respectively.  If a thinning move is applied, then some maximal base
  is collapsed onto its dual so that $\mbc \to \mbc'$ is actually a
  retraction. It follows that after splitting some maximal sections
  the restriction $U \to U'$ is also a retraction, so there is a
  natural restriction $\mbc'$. For entire transformations, the initial
  subdivisions and transfers will at most split $U$ into more open
  intervals. This is immediately followed by a collapse, so there is a
  well defined $<'$ as before.
\end{proof}

We now describe the Rips sequence for $\mbc$ induced by a track $t$
efficiently carried by $\mbc$.

\begin{enumerate}[(1)]
\item\label{it:cruch} Crush any $\pi_1$-trivial annuli and delete any superfluous
  2-cells.
\item If there is point $\mbc_i$ with $\gamma(x)=1$, we apply a
  thinning move collapsing the $<$-minimal maximal naked segment
  (Definition \ref{defn:maximal-naked-seg}), to obtain
  $\mbc_{i+1}$.
\item Otherwise, if there are any unmatched bases
  \begin{enumerate}
  \item if possible, apply a Möbius move on a $<$-minimal Möbius band,
    or
  \item apply an entire transformation, then
  \end{enumerate}
 tighten all 2-cells attaching maps.
\item Once all the bases are matched stop.
\end{enumerate}

We note that our choice of ordering $<$ on $\mbc$ is by no means
canonical. However once it is made, the Rips sequence becomes
deterministic.

\begin{prop}\label{prop:rips-machine-terminate}
  Let $t \subset \complex C$ be a track in a 2-complex and let $\mbc$
  be the corresponding measured band complex given in Proposition
  \ref{prop:construct-bc}. Then after finitely many steps the Rips
  machine terminates on a band complex $(\mbc_T,t_T)$, with all bases
  matched and coinciding. The interior of the union of the bases is a
  regular neighbourhood of $t_T$ homeomorphic to $t\times (-1,1)$, $t$
  is a wedge sum of circles, and there are no connections in the
  interior of the bands.
\end{prop}
\begin{proof}
  Since all bases start off with finite integer valued length, the
  Rips machine eventually stops since every step decreases the length
  of some base by a positive integer.

  All bases are matched, and they must all coincide, since $t_T$ is
  connected and it's efficiently carried by $\mbc_T$. Now, because all
  the bases are matched, if $\mu$ is any base $|\mu|_{t_T}$, the
  number of connected components of $t_T$, must equal 1. Efficient
  carrying  also excludes the possibility of connection in
  the interior of a band. 
\end{proof}

Proposition \ref{prop-preservation} implies that the
final dual Bass-Serre tree $\dualtree{t_T}{\mbc_T}$ is
$\fungrp{\mbc}$-equivariantly isomorphic to the original
$\dualtree{t}{\complex C}$. The following definition is important for
the next section

\begin{figure}[htb]
  \centering
  \begin{tikzpicture}[scale=0.6, every node/.style={scale=0.6}]
    \node [draw, cloud, cloud puffs=6, minimum width = 5cm, minimum height=5cm, fill=gray!60] at (-1,0) {};
    \node [draw, cloud, cloud puffs=8, minimum width = 5cm, minimum height=5cm, fill=gray!60] at (5,0) {};
    \draw[fill=white] (0,0) .. controls +(-0.5,0) and +(-0.5,0) .. (0,2) -- 
    (4,2) ..controls +(0.5,0) and +(0.5,0) ..(4,0) -- (0,0);
    \draw[dashed,color=gray] (0,0) .. controls +(0.5,0) and +(0.5,0) .. (0,2) 
    (4,2) ..controls +(-0.5,0) and +(-0.5,0) ..(4,0);
    \draw[ultra thick,color=black] (2,0) .. controls +(-0.5,0) and
    +(-0.5,0) .. (2,2);
    \draw[ultra thick,dashed,color=black] (2,2) .. controls +(0.5,0) and +(0.5,0) .. (2,0);
    \begin{scope}[yshift=-2cm]
      \draw[fill=white] (0,0) .. controls +(-0.5,0) and +(-0.5,0) .. (0,2) -- 
      (4,2) ..controls +(0.5,0) and +(0.5,0) ..(4,0) -- (0,0);
      \draw[dashed,color=gray] (0,0) .. controls +(0.5,0) and +(0.5,0) .. (0,2) 
      (4,2) ..controls +(-0.5,0) and +(-0.5,0) ..(4,0);
      \draw[ultra thick,color=black] (2,0) .. controls +(-0.5,0) and
    +(-0.5,0) .. (2,2);
    \draw[ultra thick,dashed,color=black] (2,2) .. controls +(0.5,0) and +(0.5,0) .. (2,0);
    \end{scope}
  \end{tikzpicture}
  \caption{A band complex in terminal form, i.e. all bases are
    matched. As a graph of spaces, the edge space is clearly
    visible. There is a unique track, drawn in black, that efficiently
    carried by this band complex. (Recall that tracks are connected by
    definition.)}
  \label{fig:terminal-form}
\end{figure}
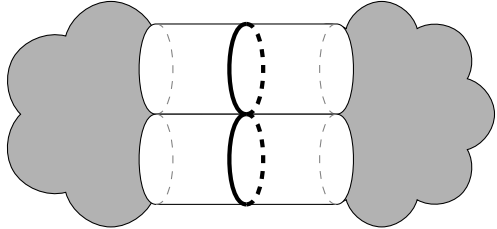

\begin{defn}\label{defn:terminal-form}
  A band complex $\mbc$ is in \define{terminal form} if it is as
  described in conclusion of Proposition
  \ref{prop:rips-machine-terminate}, see Figure \ref{fig:terminal-form}
\end{defn}
 
\section{The elimination process}
\label{sec:EP}
We will now turn our attention to the set of all tracks that are
efficiently carried by a band complex.

Let $\mbc$ be a band complex equipped with an ordering $<$ on the
union of maximal sections (Definition \ref{defn-ordering}). Then,
given a track $t \subset \mbc$, the Rips machine (Section
\ref{sec:Rips-machine}) will perform a specific derived transformation
$(\mbc,t) \to(\mbc',t')$. The type of transformation, either an
annulus crush, a thinning move, a Möbius move, or an entire
transformation, is determined by $\mbc$ and the ordering $<$, but not
the track $t$ it carries.

\begin{figure}[htb]
  \centering
  \begin{tikzpicture}[yscale=0.7, every node/.style={scale=0.8},inner sep=0.5mm,scale=0.8]
  \draw[thick,fill=gray!60] 
  (0,0)--(1,2)--(3,2)--(2,0)--cycle 
  (3,2)--(5,2)--(6,0)--(4,0)--cycle
  (1,2)--(1,5)--(5,5)--(5,2)--cycle
  (1,5)--(0,7)--(2,7)--(3,5)--cycle
  (3,5)--(4,7)--(6,7)--(5,5)--cycle;
  \draw[ultra thick,dashed] (1,7)--(2,5)--(2,2)--(1,0);
  \draw[ultra thick,dashed] (4.5,7)--(3.5,5)--(3.5,2)--(4.5,0);
  \draw[ultra thick,dashed] (5.5,7)--(4.5,5)--(4.5,2)--(5.5,0);
  \node[fill=gray!60] (delta) at (3,4.5) {$\mu$};
  \node[fill=gray!60] (dualdelta) at (3,2.5) {$\dual\mu$};
  \node[fill=gray!60] (dualdelta) at (2,5.25) {$\delta$};
  \node[fill=gray!60] (dualdelta) at (4,5.25) {$\eta$};
  \node[fill=gray!60] (dualdelta) at (2,1.75) {$\lambda$};
  \node[fill=gray!60] (dualdelta) at (4,1.75) {$\rho$};

  \begin{scope}[yshift=-2cm,scale=0.8] 
    \base{0}{0}{4}{\mu}
    \base{1}{0}{2}{\delta}
    \base{1}{2}{2}{\eta}

    \opdualbase{1}{5}{4}{\mu}
    \base{0}{5}{2}{\lambda}
    \base{0}{7}{2}{\rho}
  \end{scope}

  \begin{scope}[xshift=8cm,yshift=2cm]
    \draw[thick,fill=gray!60] 
    (0,0)--(1,2)--(3,2)--(2,0)--cycle 
    (3,2)--(5,2)--(6,0)--(4,0)--cycle
    (3,2)--(3,5)--(5,5)--(5,2)--cycle
    (1,2)--(0,4)--(2,4)--(3,2)--cycle
    (3,5)--(4,7)--(6,7)--(5,5)--cycle;
    \draw[ultra thick,dashed] (4.5,7)--(3.5,5)--(3.5,2)--(4.5,0);
    \draw[ultra thick,dashed] (5.5,7)--(4.5,5)--(4.5,2)--(5.5,0);
    \draw[ultra thick,dashed] (1,4)--(2,2)--(1,0);

    \begin{scope}[yshift=-2cm,scale=0.8] 
    \base{0}{2}{2}{\mu}
    
    \base{1}{2}{2}{\eta}

    \opdualbase{1}{7}{2}{\mu}
    \base{1}{5}{2}{\delta}
    \base{0}{5}{2}{\lambda}
    \base{0}{7}{2}{\rho}
  \end{scope}
  \end{scope}

  \begin{scope}[yshift=-10 cm]
    \draw[thick,fill=gray!60] 
    (0,0)--(1,2)--(3,2)--(2,0)--cycle 
    (3,2)--(5,2)--(6,0)--(4,0)--cycle
    (1,2)--(1,5)--(5,5)--(5,2)--cycle
    (1,5)--(0,7)--(2,7)--(3,5)--cycle
    (3,5)--(4,7)--(6,7)--(5,5)--cycle;
    \draw[ultra thick,dashed] (1,7)--(2,5)--(1.5,2)--(0.5,0);
    \draw[ultra thick,dashed] (4.5,7)--(3.5,5)--(2.5,2)--(1.5,0);
    \draw[ultra thick,dashed] (5.5,7)--(4.5,5)--(4.5,2)--(5.5,0);
    \begin{scope}[xshift=8cm,yshift=-2cm]
      \draw[thick,fill=gray!60] 
      (0,0)--(1,2)--(3,2)--(2,0)--cycle 
      (3,2)--(5,2)--(6,0)--(4,0)--cycle
      (2,2)--(3,5)--(5,5)--(5,2)--cycle
      (1,2)--(0,4)--(2,4)--(2,2)--cycle
      (3,5)--(4,7)--(6,7)--(5,5)--cycle;
      \draw[ultra thick,dashed] (5.5,7)--(4.5,5)--(4.5,2)--(5.5,0);
      \draw[ultra thick,dashed] (4.5,7)--(3.5,5)--(2.5,2)--(1.5,0);
      \draw[ultra thick,dashed] (1,4)--(1.5,2)--(0.5,0);
      
      \begin{scope}[yshift=+8cm,scale=0.8] 
        \base{0}{2}{2}{\mu}
        \base{1}{2}{2}{\eta}
        
        \opdualbase{1}{6}{3}{\mu}
        \base{1}{5}{1}{\delta}
        \base{0}{5}{2}{\lambda}
        \base{0}{7}{2}{\rho}
  \end{scope}

    \end{scope}
  \end{scope}

  \draw[->, ultra thick] (8,-1) -- (9,0);
  \draw[->, ultra thick] (8,-2) -- (9,-3);
 
\end{tikzpicture}
\caption{Two different tracks carried by the same band complex $\mbc$
  give distinct combinatorial outcomes after applying an entire
  transformation (see Section \ref{sec:entire-trans}.) Here the
  carrier $\mu$ moves base $\delta$ onto its dual. The base diagrams
  illustrate the ordering $<$ (see Definition \ref{defn-ordering}) on
  the union of maximal sections.}
\label{fig:different-outcomes}
\end{figure}
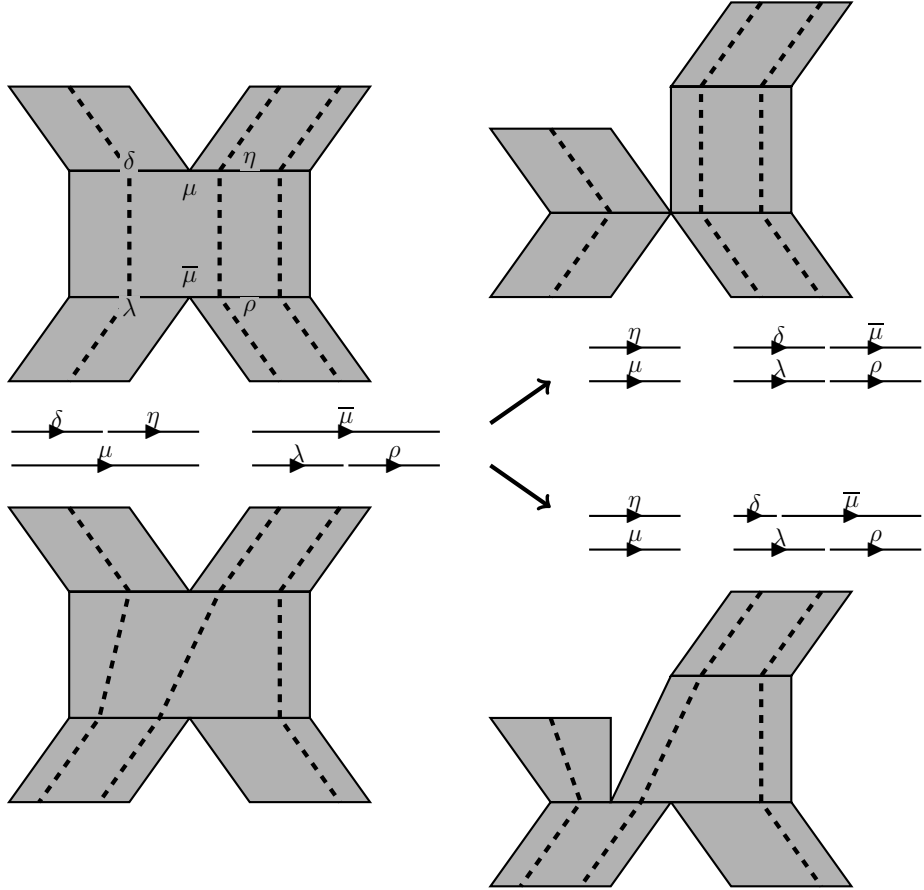

Although the type of transformation doesn't depend on the track $t$
carried by $\mbc$, the combinatorial equivalence class of the
resulting $\mbc'$ does depend on the track $t$. For example consider
Figure \ref{fig:different-outcomes} which shows two different
combinatorial outcomes coming from two different tracks carried by the
same band complex.
Seeing as we want to study the set of all tracks carried by $\mbc$ we
must consider all these combinatorial outcomes simultaneously.

\subsubsection{The elimination tree $\ETone{C}$}\label{sec:ETone}

Given a band complex $\mbc$ with a track $t \subset\mbc$, the Rips
machine gives a sequence $(\mbc,t) \to \ldots \to (\mbc_T,t_T)$ with
$(\mbc_T,t_T)$ in terminal form. If we want to consider all tracks
carried by $\mbc$, then we must have a branching sequence or, in other
words, a rooted directed tree.

\begin{defn}\label{defn:tracks}
  Let $\mbc$ be a band complex. Denote by $\effTracks{\mbc}$ the set of
  tracks efficiently carried by $\mbc$.
\end{defn}

\begin{defn}[Combinatorially equivalent derived transformations]\label{defn:combinatorially-equiv}
  Let $t_1$ and $t_2$ be two tracks efficiently carried by a band
  complex $\mbc$. The derived transformations
  $(\mbc,t_1) \to (\mbc_1',t_1')$ and $(\mbc,t_2) \to (\mbc_2',t_2')$
  are \define{combinatorially equivalent} if there is a commuting
  homeomorphism $e$
  \[
    \begin{tikzpicture}
      \node (mbc1) at (-1,-1) {$\mbc_1'$};
      \node (mbc) at (0,0) {$\mbc$};
      \draw[->] (mbc)--(mbc1);
      \node (mbc2) at (1,-1) {$\mbc_2'$};
      \draw[->] (mbc)--(mbc2);
       \draw[->] (mbc1)--node[above]{$e$}(mbc2);
   \end{tikzpicture}
 \] that is a combinatorial equivalence of band complexes in the sense of
 Definition \ref{defn:BC-equiv}.
\end{defn}

\begin{conv}\label{conv:equiv-is-equal}
  In section \ref{sec:Rips-machine} it was convenient to consider band
  complexes as being measured in order to precisely describe
  continuous quotient maps. For the rest of the paper, unless stated
  otherwise, band complexes $\mbc$ without tracks will be considered
  equal if they are combinatorially equivalent in the sense of
  Definition \ref{defn:BC-equiv}.
\end{conv}

The entire transformations shown in Figure \ref{fig:different-outcomes} are
not combinatorially equivalent. We note that, in the notation of the
above definition, the tracks $t_1',t_2'$ need no longer give
combinatorially equivalent derived transformations of
$\mbc_1' = \mbc_2'$.

\begin{defn}[Derived transformations of a band
  complexe]\label{defn:Derived-BC}
  Let $\mbc$ be a band complex equipped with an ordering $<$ of the
  union of its maximal sections (Definition \ref{defn-ordering}). For
  every $t \in \effTracks{\mbc}$ there is a corresponding derived
  transformation $(\mbc,t)\to(\mbc',t')$ with resulting band complex
  $\mbc'$. The type of this transformation, either an annulus crush, a
  Möbius move, a thinning move, or an entire transformation, depends on
  the underlying band complex $\mbc$ (and the ordering). Let
  $\mbc_1,\dots,\mbc_n$ denote the finite set of combinatorial
  equivalence classes of the resulting band complexes. A \define{derived
    transformation of a band complex} $\mbc$ is the operation that
  produces the finite collection of continuous
  maps\begin{equation} \label{eqn:EP-descendants}
    \begin{tikzpicture}
      \node (1) at (0,0) {$\mbc$};
      \node (2) at (-1,-0.75) {$\mbc_1$};
      \node (3) at (1,-0.75) {$\mbc_n$};
      \node (4) at (0,-0.75) {$\cdots$};
      \draw[->] (1) -- (2);
      \draw[->] (1) -- (3);
    \end{tikzpicture}
    .
  \end{equation}
  
\end{defn}

We remark that we can algorithmically construct the set
$\mbc_1,\ldots,\mbc_n$.

\begin{lem}\label{lem:induced-order}
  The if $<$ is an ordering on $\mbc$ then all its children obtained
  by a derived transformation have well-defined induced orderings.
\end{lem}
\begin{proof}
  The induced ordering given in Lemma \ref{lem:order-preserve} does
  not depend on the track, only on the continuous map between the
  underlying band complexes.
\end{proof}

\begin{conv}\label{conv:order}
  For the rest of the paper, unless stated otherwise, we will assume
  that a band complex $\mbc$ comes equipped with such an ordering $<$
  of the union of its maximal sections.
\end{conv}

The \define{elimination process} for a 2-complex $\complex C$ is the
construction of the elimination tree $\ETone{\complex C}$, a directed
rooted tree defined inductively as follows:

\begin{enumerate}[(1)]
\setcounter{enumi}{-1}
\item The root of $\ETone{C}$ is the polygonal complex $\complex C$.
\item\label{it:et-first-children} The set of children of $C$ is the finite collection of band
  complexes $\mbc_1,\ldots,\mbc_{n_\complex{C}}$ provided by
  Proposition \ref{prop:finite-ubcs} that can efficiently carry all
  tracks of $\complex C$. For each
  $\mbc_i, 1\leq i \leq n_\complex{C} $, we equip the union of maximal
  sections (Definition \ref{defn:maximal-section}) with an ordering
  $<$ as in Definition \ref{defn-ordering}.
\item\label{it:terminal-leaf} If a band complex $\mbc_v$ in $\ETone{\complex C}$ is in terminal form
  (Definition \ref{defn:terminal-form}), then it is called a
  \define{terminal leaf.}
\item\label{it:inadmissible} If a band complex $\mbc_v$ in $\ETone{\complex C}$ can not be brought to
  terminal form via derived moves because the union of the bands is not
  connected or all bases are matched, but there are 2-cell attaching
  maps that intersect the interior of the bands, then it is a called
  an \define{inadmissible leaf.}
\item\label{it:ET-otherwise} Otherwise we continue to grow
  $\ETone{\complex C}$ at a band complex $\mbc_v$ by adding its
    descendants with a derived transformation (Definition
    \ref{defn:Derived-BC}.)  Equip the union of bases of each
    descendant $\mbc_{v'}$ of $\mbc_v$ with the induced order $<$
    given by Lemma \ref{lem:induced-order}.
\end{enumerate}

As will be explained in the next section $\ETone{\complex C}$ gives a
way to encode the set of tracks that can lie in the polygonal complex
$\complex C$. It follows that in general it is infinite.

\subsection{The sets of tracks in a band complex 
organized by open
  neighbourhoods.}\label{sec:set-of-tracks}

\begin{lem}\label{lem:going-backwards}
  Let $\mbc \to \mbc'$ be one of the continuous maps of the derived
  transformation on $\mbc$ (Definition \ref{defn:Derived-BC}
  (\ref{eqn:EP-descendants})). Suppose that $\mbc'$ is equipped with a
  measure (Definition \ref{defn:measured}.) Then there is a
  well-defined pullback measure on $\mbc$. Furthermore, as long as the
  union of the interiors of the bands in $\mbc'$ is connected, the
  maximal measure of each base $\mu$ of $\mbc$ is no more than the sum
  of the measures of the bases in $\mbc'$ that are in the image of
  $\mu$ via the map $\mbc \to \mbc'$.
\end{lem}
\begin{proof}
  It is now enough to consider the basic moves in Section
  \ref{sec:basic-moves}, i.e. the zips, collapses, annulus crushes,
  vertical and horizontal subdivisions. In all cases given such a
  transformation $\hat\mbc \to \mbc'$ as a continuous map, there is a
  unique pull back measure we can put on $\hat\mbc$.  The upper bound
  on the measure of the bases is of $\mbc$ is obvious from the
  definitions of the derived moves.
\end{proof}

\begin{cor}[Going backwards]\label{cor:going-backwards}
  Let $\mbc \to \mbc'$ be one of the continuous maps of the derived
  transformation on $\mbc$. Suppose that $\mbc'$ efficiently carries a
  track $t$ (Definition \ref{def:efficiently-carry}.) Then there is a
  unique track $\hat t$ that is efficiently carried by $\mbc$ such
  that $(\mbc,\hat t) \to (\mbc',t)$ is a derived transformation in
  the Rips machine.
\end{cor}
\begin{proof}
  Since $\mbc'$ efficiently carries $t$, we can put a measure on
  $\mbc'$ (Definition \ref{defn:measured}) that corresponds to
  $t$. Lemma \ref{lem:going-backwards} gives a pullback measure on
  $\mbc$, which induces a track $\hat t \subset \mbc$ which induces
  the derived transformation. Furthermore it is routine to check for
  the basic moves in Section \ref{sec:basic-moves} that $\hat t$ is
  indeed efficiently carried by $\mbc$, provided $t$ is efficiently
  carried by $\mbc'$.
\end{proof}

Any band complex in terminal form (Definition
\ref{defn:terminal-form}) efficiently carries a unique track. We will
now show how repeatedly going backwards enables us to use
$\ETone{\complex C}$ to organize the collection of tracks that are
efficiently carried by $\complex C$.

\begin{prop}\label{prop:basis-for-tracks}
  There is a bijective correspondence between the set of tracks in
  $\complex C$ and $\{\mbc_l\}$ the set of terminal leaves (Section
  \ref{sec:ETone} (\ref{it:terminal-leaf})) of $\ETone{\complex C}$.
\end{prop}
\begin{proof}
  Let $\mbc$ be a band complex at the top level of
  $\ETone{\complex C}$. As a topological space it is homeomorphic to
  $\complex C$. Any track efficiently carried by $\mbc$ is obviously a
  track in $\complex C$.

  If $t \subset \mbc$ is a track, by Proposition
  \ref{prop:rips-machine-terminate} the Rips sequence for $(\mbc, t)$
  will give a path in $\ETone{\complex C}$ from $\mbc$ to some band
  complex in terminal form $\mbc_l$. This map from tracks to leaves is
  injective since, having fixed $<$, the outcome of a particular
  derived move $(\mbc,t) \to (\mbc',t')$ move depends only on the track
  $t$.

  On the other hand let $\mbc_l$ be a terminal leaf of
  $\ETone{\complex C}$. We may metrize all its bases to have length 1,
  so that each band is explicitly parameterized as
  $[0,1]\times [-1,1]$. $\mbc_l$ efficiently carries the track $t_l$
  which intersects each band as $\{\frac{1}{2}\} \times [-1,1]$. Now
  starting at $(\mbc_l,t_l)$ and repeatedly going backwards
  (Corollary \ref{cor:going-backwards}) in $\ETone{\complex C}$, we
  obtain $(\mbc,t)$, where $t \subset \mbc$ is an efficiently
  carried track.
\end{proof}

If $\mbc$ is some band complex not in terminal form(as in Section
\ref{sec:EP} (\ref{it:terminal-leaf})) then any efficiently carried
track $t\subset\mbc$ gives rise to one of the children in the
elimination tree
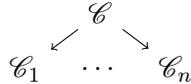
\[
\begin{tikzpicture}
  \node (1) at (0,0) {$\mbc$};
  \node (2) at (-1,-0.75) {$\mbc_1$};
  \node (3) at (1,-0.75) {$\mbc_n$};
  \node (4) at (0,-0.75) {$\cdots$};
  \draw[->] (1) -- (2);
  \draw[->] (1) -- (3);
\end{tikzpicture}
\]

It follows that there are injective maps
\begin{equation}\label{eqn:iotamap}
  \iota_j:\effTracks{\mbc_j}\into \effTracks{\mbc}
\end{equation}
whose images give a cover\[
\effTracks{\mbc} = \bigcup_j \iota_j\left(\effTracks{\mbc_j}\right).\]

\begin{defn}\label{defn:complex-NHD}
  Let $\mbc_v$ be a band complex in an elimination tree
  $\ETone{\complex C}$. We denote the \define{track
    neighbourhood}
  \[\trnhd{\ETone{\complex C}}{\mbc_v} \subset \effTracks{\mbc}\]
  to be the set of tracks carried by $\mbc$ obtained by composing the
  maps (\ref{eqn:iotamap}) going from $\mbc_v$ all the way back to
  $\mbc$.
\end{defn}

Equivalently if $\mbc_v$ is a band complex in $\ETone{\complex C}$, then
there is a natural inclusion
$\ETone{\mbc_v} \subset \ETone{\complex C}$. Proposition
\ref{prop:basis-for-tracks} immediately gives the inclusion
$\effTracks{\mbc_v} \into \effTracks{\mbc}$ obtained by iterating (\ref{eqn:iotamap}).

\subsection{Analogies with surface train tracks}
If the reader has some familiarity with surface train tracks, the
following analogies may be helpful.

If $\Sigma$ is a surface then a train track $\tau\subset \Sigma$ is
analogous to a band complex structure $\mbc$ on a 2-complex
$\complex{C}$. If we assign positive integer weights to the branches
of a train track $\tau$ satisfying the switch equations then we get a
multicurve in $\Sigma$, which is analogous to a pattern in
$\complex{C}$. For us a track in a 2-complex is analogous to a simple
closed curve.

The assignment of weights to branches of a train track gives rise to a
splitting sequence \[\tau = \tau_0 \to \cdots \] which will eventually
split $\tau$ into a multicurve if the weights are positive integers. This is
analogous to a Rips process.

On the other hand if we put a measure on a train track then we can
consider all possible train tracks \[
\begin{tikzpicture}
  \node (1) at (0,0) {$\tau$};
  \node (2) at (-1,-0.75) {$\tau_1$};
  \node (3) at (1,-0.75) {$\tau_n$};
  \node (4) at (0,-0.75) {$\cdots$};
  \draw[->] (1) -- (2);
  \draw[->] (1) -- (3);
\end{tikzpicture}
\] that can be obtained from $\tau$ via a splitting move. Iterating,
this gives an analogue to the elimination tree. In fact we will get an
actual tree if we impose some kind of order $<$ which specifies at
which switch to split at each step.  If eventually the train track has
split itself into a simple closed curve $\tau_T$ then we have a train
track in terminal form. Assigning weight 1 to the branch and working
backwards (i.e. using folding sequences) gives us a ``complicated''
simple closed curve in $\Sigma$.

If we were to consider the set of projectivized measured laminations,
then irrational laminations would give infinite splitting
sequences. Furthermore the neighbourhoods of Definition
\ref{defn:complex-NHD} are somewhat analogous to open neighbourhoods
in the Hausdorff topology on laminations. Indeed, two laminations are
``close'' if the corresponding splitting sequences coincide for a long
time. This all carries through to measured laminations on cell
complexes, but this technology is not needed, and the ordering $<$
will cause us to stay stuck in a single minimal component.

\subsection{Inadmissibility from $\kappa$-acylindricity}
\label{sec:kappa-inadmissible}
Up to now the $\kappa$-acylindricity of the dual Bass-Serre tree $\dualtree t
\mbc$, has not been used at all.

\begin{defn}\label{defn:kappaTracks}
  A track $t$ efficiently carried by a band complex $\mbc$ is called a
  \define{$\kappa$-track} if the dual tree $\dualtree t \mbc$ is
  $\kappa$-acylindrical. We denote by $\kappaTracks{\mbc}$ the set of
  $\kappa$-tracks efficiently carried by $\mbc$.
\end{defn}

We give two extra criteria to exclude vertices $\mbc_v$ of
$\ETone{\complex C}$ because their track neighbourhood
$\trnhd{\mbc}{\mbc_v}$ cannot contain any $\kappa$-tracks.

\begin{lem}\label{lem:short}
  If $\fungrp{\mbc}$ has no elements of order 2 and
  $t \in \kappaTracks\mbc$, then if $\mu \bandeq \dual\mu$, then
  either $\basepair\mu$ forms an annulus that can be crushed (which
  decreases the number of bands,) or $|\mu|_t \leq \kappa$ where
  $|\mu|_t$ is the hitting measure.
\end{lem}
\begin{proof}
  By the $\kappa$-acylindricity assumption any element that fixes an
  arc of length $\kappa+1$ in the dual tree $\dualtree t \mbc$ must be
  trivial. By assumption $\band\mu$ either forms an annulus
  $\ann{\mu}$ or a Möbius strip $\mob\mu$. In both cases
  $\gp\left({\ann\mu}\right)$ or $\gp\left({\mob\mu}\right)$ is
  generated by an element $g\neq 1$, since clause (\ref{it:cruch}) in
  the Rips machine (Section \ref{sec:Rips-machine}) crushes
  $\pi_1$-trivial annuli and $\gp\left({\mob\mu}\right)$ must act
  nontrivially on $\dualtree t \mbc$.

  We pass to the universal cover and consider the equivariant map
  $\quotmap$ of Proposition\ref{prop:tree}(\ref{eq:phi-map}), and we see that in the annulus case
  $g$ fixes an arc of $\dualtree t \mbc$ of length $|\mu|_t$. In the
  Möbius strip case $g^2\neq 1$ fixes an arc of length $|\mu|_t$ in
  $\dualtree t \mbc$. In the Möbius band case this forces $|\mu|_t$ to
  be at most $\kappa$. In the annulus case, if $|\mu|_t > \kappa$ then
  we can crush it since $\gp\left({\ann\mu}\right) = \{1\}$.
\end{proof}

As an immediate Corollary we have:

\begin{prop}\label{prop:kappa-inadmissible}
  Let $\mbc_v$ be a band complex in $\ETone{\complex C}$. If along
  some path\[
  p:\mbc_u \to \cdots \to \mbc_v
  \]
  in $\ETone{\complex C}$ either some annulus gets subdivided more
  than $\kappa+1$ times or some base $\mu$ that formed a Möbius pair
  with its dual gets shortened or vertically subdivided more than
  $\kappa+1$ times, then $\trnhd{\ETone{\complex C}}{\mbc_v}$ doesn't
  contain any $\kappa$-acylindrical tracks.
\end{prop}

\begin{defn}\label{defn:kappa-inadmissible}
  We call a path $p$ in $\ETone{\complex C}$ such as the one given in
  Proposition \ref{prop:kappa-inadmissible}
  \define{$\kappa$-inadmissible}.
\end{defn} 

\subsection{Automorphic minimality and repetitions}

\begin{defn}[Size and minimality]\label{defn:track-size}
  The \define{size} of a track $t\subset \mbc$ in a band complex is
  the finite sum\[ \size{t}=\sum_{\mu} |\mu \cap t| \]
  where $\mu$ ranges over the bases of $\mbc$. A track $t$ is called
  \define{automorphically minimal} if among all other tracks $t'$ such
  that $t \sim_{\aut(\fungrp{\mbc})} t'$, $\size{t} \leq \size{t'}$.
\end{defn}

This next lemma is easy to prove from the definitions of the basic moves.

\begin{lem}\label{lem:size}
  Let $m:\mbc \to \mbc'$ be a basic transformation. If $t_1'$ is a
  track efficiently carried by $\mbc'$ then there is a corresponding
  track $t_1$ efficiently carried by $\mbc$ such that
  $m:(\mbc,t_1) \to (\mbc',t_1')$. Furthermore if $t_2'$ is
  efficiently carried by $\mbc'$ and $\size{t_1'}<\size{t_2'}$ then
  $\size{t_1} < \size{t_2}$ where $t_2$ is the track efficiently
  carried by $\mbc$ corresponding to $t_2$.
\end{lem}

This notion of automorphic minimality may seem convoluted, but  the
proof of the following proposition may clear things up for the reader.

\begin{prop}[Repetitions and minimality]\label{prop:repetition}
  Let $\mbc_v$ be a band complex in $\ETone{\complex C}$. If along
  some path \[p:\mbc_u \to \cdots \to\mbc_v\] in $\ETone{\complex C}$
  there are two combinatorially equivalent band complexes (Definition
  \ref{defn:BC-equiv}) $\mbc_u \stackrel{\approx}{\to} \mbc_v$, i.e. a
  \define{repetition}, then $\trnhd{\ETone{\complex C}}{\mbc_v}$
  cannot contain any minimal tracks.
\end{prop}
\begin{proof}
  Suppose towards a contradiction that there was a track
  $t \subset \mbc$ in $\trnhd{\ETone{\complex C}}{\mbc_v}$ that was
  minimal. Let $q_v:\mbc_v \to \cdots \to \mbc_l$ be the path to the
  terminal leaf in $\ETone{\complex C}$ corresponding to $t$ (recall
  Proposition \ref{prop:basis-for-tracks}.)

  Let $p_u$ be the concatenation of paths $p$ and $q_v$, i.e.
  $p_u:\mbc_u \to \cdots \to \mbc_v \to \cdots \to \mbc_l$. Since
  $\mbc_u \approx \mbc_v$, we can attach the path $q_v$ to $\mbc_u$ to
  get a corresponding path $q_u:\mbc_u \to \cdots \to \mbc_{l'}$,
  where $\mbc_{l'} \approx \mbc_{l}$. Let $t' \subset \mbc$ be the
  track corresponding to $\mbc_{l'}$, and denote by $t'_i$
  (respectively $t_i$) the image of $t'$ (respectively $t$) in
  $\mbc_i$, should there be such an image.  Derived moves always
  decrease the lengths of bases; thus, in $\mbc_u$,
  $\size{t'_u} < \size{t_u}$. Working backwards in
  $\ETone{\complex C}$, i.e. repeatedly applying Corollary
  \ref{cor:going-backwards}, all the way back to a direct descendant
  $\mbc$ of $\complex C$ yields $\size{t'}< \size{t}$.

  On one hand, by the definitions of $q_v$ and $q_u$, we have
  $(\mbc_v,t_v) \stackrel{\approx}{\to} (\mbc_u,t'_u)$. On the other
  hand we have a sequence of derived
  moves\[ (\mbc_u,t_u) \to \cdots \to (\mbc_v,t_v)
  \] this gives a composition of continuous maps\[
  \varphi: \mbc_u \to \cdots \to \mbc_v \stackrel{\approx}{\to} \mbc_u 
  \]
  which by Proposition \ref{prop-preservation} induces an isomorphism
  on $\pi_1$, hence $\varphi_\sharp \in
  \aut{(\pi_1(\mbc_u))}$.
  Furthermore, by construction $\varphi(t_u)=t'_u$, so again Proposition
  \ref{prop-preservation} gives us that $\dualtree{t_u}{\mbc_u}$ is
  $\varphi_{\sharp}$-equivariantly isomorphic to $\dualtree{t'_u}{\mbc_u}$
  where the action on the second tree is given by
  $(g,x) \mapsto \varphi_{\sharp}(g)\cdot x$. It follows that $t'$ and
  $t$ are automorphically equivalent, contradicting the minimality of
  $t$.
\end{proof}

\begin{defn}\label{defn:rep-inadmissible}
  A path $p$ in $\ETone{\complex C}$ that satisfies the hypotheses of
  Proposition \ref{prop:repetition} is called
  \define{repetition-inadmissible.}
\end{defn}

\subsection{Restricted elimination processes}\label{sec:restricted-ep}
In order to construct $\ETone{\complex C}$ we will sometimes have to
construct auxiliary elimination trees that are rooted at band
complexes $\mbc$ in $\ETone{\complex C}.$

Let $J \subset \mbc$ be a union of maximal sections. We redefine the
order $<$ so that the maximal sections in $J$ are terminal. The
restricted elimination tree is used to study how the bases of $\mbc$
can be moved into $J$. If $\mbc \to \mbc'$ is a derived
transformation, which is a continuous map, then $J$ has a well-defined
image in $\mbc'$ which we also denote by $J$. $\ETone{\mbc,J}$ is
constructed as follows:
\begin{enumerate}[(1)]
\setcounter{enumi}{-1}
\item The root of $\ETone{\mbc,J}$ is $\mbc$.
\item If every base of a band complex $\mbc_l$ in $\ETone{\mbc,J}$ is
  contained in $J$ then $\mbc_l$ is called a \define{$J$-terminal leaf}.
\item \define{Inadmissible leaves} are defined the same way as for
  $\ETone{\complex C}$.
\item Otherwise we apply a corresponding derived transformation,
  either a Möbius move on $\basepair\mu$, a thinning move starting
  at $\mu$, or an entire transformation with leading base $\mu$ to create the children of $\mbc_v$. Equip
  the union of bases of every child $\mbc_v'$ of $\mbc_v$ with the
  induced order $<$.
\end{enumerate}

Restricted elimination processes will be required for some
subprocesses of our main algorithm. It is obvious that the various
inadmissibility criteria for a standard elimination tree
$\ETone{\complex C}$ also hold for restricted elimination trees.

\begin{conv}
  Many statements about restricted elimination trees will also follow
  for the standard elimination tree by replacing $\ETone{\complex C}$
  by the elimination trees
  $\ETone{\mbc_1,\emptyset},\ldots\ETone{\mbc_{n_{\complex
        C}},\emptyset}$ where $\mbc_1,\ldots, \mbc_{n_{\complex C}}$
  are the children of $C$ (see step (\ref{it:et-first-children}) in
  Section \ref{sec:ETone}.) These elimination trees are contained in
  $\ETone{\complex C}$. We will therefore assume that results about
  restricted elimination processes will apply to the standard
  elimination process, the latter being a special case.
\end{conv}

\subsection{The infinite branches of $\ETone{\complex C}$.}

One of Makanin's key observations is that every infinite branch of
$\ETone{\complex C}$ stabilizes into one of three cases. The main
ingredient is the following lemma.

\begin{lem}\label{lem:finite-switches}
  If after applying an entire transformation which didn't involve a
  Möbius move we are in the situation where we must make a
  thinning move, then ($J$-relative) $\tau$-complexity decreased.
\end{lem}

In the case of band complexes, this fact is explained between
Proposition 7.5 and Proposition 7.6 of
\cite{Bestvina-Feighn-1995}. The proof consists of a straightforward complexity
counting argument. This next result is ubiquitous whenever the
elimination process/Rips machine is involved.

\begin{thm}[Fundamental classification]\label{thm:branch-types}
  Every infinite branch $\branch{}$ in $\ETone{\complex C}$
  ($\ETone{\mbc,J}$) has a tail $\branch{v} = \mbc_v \to \cdots$ of
  one of the following form:
  \begin{enumerate}[(i)]
  \item\label{it:types-thin} \define{Thinning:} Every derived
    transformation along $\branch v$ is a thinning move.
  \item\label{it:types-quad} \define{Quadratic:} Every derived
    transformation along $\branch v$ is an entire transformation. For
    all but finitely many points in the ($J$-complement of the) union
    of bases of the band complexes $\mbc_w$ along $\branch v$ we have
    $\gamma(x) = 2$ (Definition \ref{defn:gamma})
  \item\label{it:types-super} \define{Superquadratic:} Every derived transformation along
    $\branch v$ is an entire transformation. There is a whole open
    interval of points in the ($J$-complement of the) union of bases of
    the band complexes $\mbc_w$ along $\branch v$ such that $\gamma(x)
    \geq 3.$
  \end{enumerate}
  Furthermore if we require $\trnhd{\ETone{\complex C}}{\mbc_v}$ to contain
  a $\kappa$-track, then we may assume that no Möbius moves or annulus
  subdivisions occur along $\branch v$.
\end{thm}

In (\ref{it:types-quad}) above, we would really like to say that every
point in the union of bases is contained in exactly two unmatched
bases. Points on the boundary of bases, however, may be contained in up
to four distinct bases, but there are only finitely many of them.

\begin{proof}
  If the ($J$-restricted) $\tau$-complexity is 0, because we are not
  allowing long bands (see Step 4 of the thinning move, Section
  \ref{sec:thinning}), all bases are matched (all bases moved onto
  $J$) so we are at a leaf. It therefore follows by Lemma
  \ref{lem:finite-switches} that the infinite branch $\branch{}$
  eventually always consists of thinning moves, or eventually always
  consists of entire transformations.

  Suppose now that $\branch v$ is not of thinning type. If $\mbc_v$ is
  of quadratic type, then after applying an entire transformation
  $\mbc_v'$ is still quadratic. The trichotomy now follows. The fact
  that Möbius moves and annulus subdivisions stop occurring follows
  from Lemma \ref{lem:short}.
\end{proof}

König's Lemma states that every infinite rooted tree with vertices
of finite valency must have an infinite branch. This classification of
infinite branches is the foundation of the construction of a finite
subtree of $\ETone{\complex C}$ containing all the leaves
corresponding to minimal $\kappa$-acylindrical tracks.

\subsection{The admissible subtree
  $\ADone{\mbc,J}$}\label{sec:admissible-subtree}
\begin{defn}\label{defn:admissible-subtree}
  The \define{admissible elimination tree}
  $\ADone{\mbc,J} \subset \ETone{\mbc,J}$ is the subtree obtained by
  forbidding $\kappa$-inadmissible and repetition-inadmissible
  subpaths (Definitions \ref{defn:kappa-inadmissible} and
  \ref{defn:rep-inadmissible}). We similarly define the admissible
  elimination tree $\ADone{\complex C} \subset \ETone{\complex C}$.
\end{defn}

This next proposition enables us to restrict the search for tracks in
the algorithm for Theorem \ref{thm:main} to admissible elimination
trees.

\begin{prop}\label{prop:admissible-contains}
  Let $t$ be an automorphically minimal $\kappa$-track efficiently
  carried by $\mbc$. Any path $\mbc \to \cdots$ in $\ETone{\mbc,J}$
  induced by $t\subset \mbc$ must be contained in $\ADone{\mbc,J}$.
\end{prop}
\begin{proof}
  Otherwise Proposition \ref{defn:kappa-inadmissible} or
  \ref{defn:rep-inadmissible} lead to a contradiction of the
  hypotheses.
\end{proof}

\begin{prop}
  For every $n$, the subtree of radius $n$ of $\ADone{\mbc,J}$ can be
  effectively constructed.
\end{prop}
\begin{proof}
  For any band complex the collection of children (see Definition
  \ref{defn:Derived-BC}(\ref{eqn:EP-descendants})) can be constructed
  effectively and the various inadmissibility conditions can be
  verified effectively.
\end{proof}

\subsection{Reduction to the superquadratic
  case}\label{sec:repetition}
We show that infinite thinning or quadratic branch in $\ETone{\mbc,J}$
contains a repetition. From this it will follow that the admissible
elimination tree $\ADone{\mbc,J}$ does not have any infinite
thinning or quadratic branches.

\begin{lem}[c.f. {\cite[Lemma 15]{KM-IrredII}}]\label{lem:thinning-rep}
  Any sufficiently long thinning path $\mbc_v \to \cdots$ contains a
  repetition, i.e. a subpath $\mbc_u \to \cdots \to \mbc_w$ with
  $\mbc_u \approx \mbc_w$.
\end{lem}
\begin{proof}
  By Theorem \ref{thm:branch-types} we may assume that no more
  annulus subdivision occur. By Lemma \ref{lem:post-thinning}, the
  number of 2-cells and the vertical lengths of the 2-cell attaching
  maps are non-increasing, furthermore since the $\tau$-complexity is
  bounded and there are no maximal sections $\sigma$ with $\tau(\sigma)=0$
  there is a bound on the number of bases.

  A band complex is obtained by gluing bands to a graph $\Gamma$ and
  then attaching other 2-cells. Since derived transformations do not
  change anything in the exterior of union of the bands, then number
  of bands and 2-cells remains bounded, and the combinatorial lengths
  of the attaching immersions of the 2-cells is bounded, any
  sufficiently long thinning path will have a repetition.
\end{proof}

This next lemma will also be used later in Section \ref{sec:frep}.

\begin{lem}[c.f. {\cite[Case 14]{KM-IrredII}}]\label{lem:quadratic-rep}
  Any sufficiently long quadratic path $\mbc_v \to \cdots$ contains a
  repetition.
\end{lem}
\begin{proof}
  By Lemma \ref{lem:entire-trans}, the $\tau$ complexity doesn't
  increase, and we may assume that no Möbius moves or annulus
  subdivision occur. It remains to show that the vertical lengths of 
  2-cell attaching maps remain bounded; the result will then follow as
  in the previous proof.

  Consider Figure \ref{fig:1} with leading bases $\mu,\lambda$, where
  $\mu$ is the carrier.
  \begin{figure}[htb]
    \centering
    \begin{tikzpicture}[scale=0.66]
      \begin{scope}
        \clip[decorate,decoration=zigzag] (-1,-1) rectangle (3,3);
        \draw[fill=black!10!white] (0,0) -- (0,2) -- (3.5,2) -- (3.5,0) --cycle;
      \end{scope}
      \begin{scope}
        \clip[decorate,decoration=zigzag] (4.5,5) -- (3,2) -- (0,1) -- (1,5) --cycle;
        \draw[fill=black!10!white] (1.5,2) -- (2.5,4) -- (5,4) -- (4,2) --cycle;
      \end{scope}
      \draw[fill=black!10!white] (-1,4) -- (0,2) -- (1.5,2) -- (0.5,4)
      --cycle;
      \draw (0.75,2) node[above] {$\lambda$} 
      (0.75,2) node[below] {$\mu$}
      (2.5,2) node[above] {$\eta$}
      (0.75,0) node[above] {$\ol{\mu}$};
      \draw[ultra thick,join=round] (0.5,4) -- (1.5,2) -- (2.5,4)
      (-0.75,4) -- (0.25,2) -- (0.25,0);
      \draw[->] (4,2) -- (5.5,2);
      \draw[shift={(-4,-1)}] (0,-1) --node[above]{$\mu$} (3.5,-1) node[right]{\ldots}
      (4.5,-1) --node[above]{$\ol{\mu}$} (7.5,-1) node[right]{\ldots}
      (6,-1.5) node{$\ldots$}
      (0,-1.5) -- node[below]{$\lambda$} (1.45,-1.5)
      (1.55,-1.5) -- node[below]{$\eta$} (3.5,-1.5) node[right]{\ldots};
      
      \begin{scope}[shift={(7,0)}]
        \begin{scope}
          \clip[decorate,decoration=zigzag] (-1,-1) rectangle (3,3);
          \draw[fill=black!10!white] (1.5,0) -- (1.5,2) -- (3.5,2) -- (3.5,0) --cycle;
        \end{scope}
        \begin{scope}
          \clip[decorate,decoration=zigzag] (4.5,5) -- (3,2) -- (0,1) -- (1,5) --cycle;
          \draw[fill=black!10!white] (1.5,2) -- (2.5,4) -- (5,4) -- (4,2) --cycle;
        \end{scope}
        \begin{scope}[shift={(0,-2)}]
          \draw[fill=black!10!white] (-1,4) -- (0,2) -- (1.5,2) -- (0.5,4) --cycle;
        \end{scope}
        \draw (0.75,0) node[above] {$\lambda$} 
        (2.5,2) node[below] {$\mu$}
        (2.5,2) node[above] {$\eta$}
        (2.5,0) node[above] {$\ol{\mu}$};
        \draw[ultra thick,join=round] (0.5,2) -- (1.5,0) --(1.5,2) -- (2.5,4)
        (-0.75,2) -- (0.25,0);
        \draw[shift={(-2,-1)}] (1.5,-1) --node[above]{$\mu$} (3.5,-1) node[right]{\ldots}
        (6.05,-1) --node[above]{$\ol{\mu}$} (7.5,-1)
        node[right]{\ldots} 
        (6,-1.5) node{$\ldots$}
        (4.5,-1) -- node[above]{$\lambda$} (5.95,-1)
        (1.5,-1.5) -- node[below]{$\eta$} (3.5,-1.5) node[right]{\ldots};
        
      \end{scope}
    \end{tikzpicture}
    \caption{Above, the attaching maps of 2-cells under an entire
      transformation in the quadratic case. Below, the resulting
      ordered base configurations.}\label{fig:1}
  \end{figure}
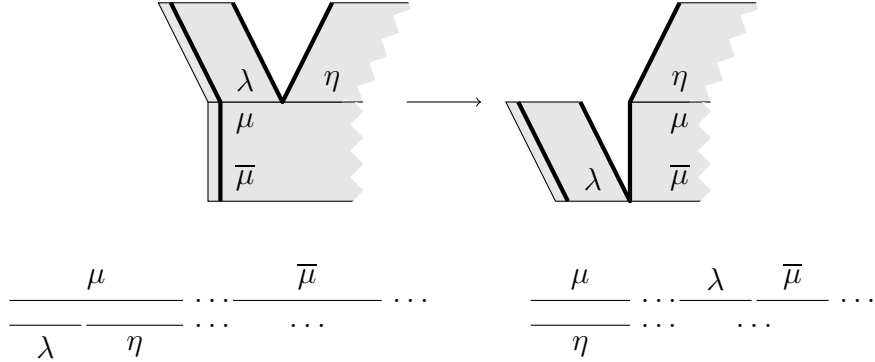
  The only way the attaching map of a 2-cell do decrease in vertical
  length is if it it has a subpath as in the shaded path on the left
  of Figure \ref{fig:1} that travels through
  $\band{\lambda} \cup \band{\mu}$. The only way for a segment in the
  boundary of a 2-cell to get ``stretched'' is if travels through the
  vertical sides of $\band{\lambda}$ and $\band{\eta}$. In this case
  we have a segment $\sigma$ of length 2 that gets stretched to a
  segment $\sigma'$ of length 3. After the entire transformation,
  however, $\eta$ and $\mu$ are now leading bases, this means that
  after the next entire transformation there is a subsegment of length
  2 $\sigma'$ that gets shortened again back to length 1.

  It therefore follows that the vertical lengths of the boundaries of
  2-cells remain bounded throughout the quadratic path and the result follows.
\end{proof}

These two Lemmas immediately imply the following corollary, which
pretty much sets the tone for the rest of the paper.

\begin{cor}\label{cor:must-be-superquadratic}
  The admissible elimination subtree
  $\ADone{\mbc,J} \subset \ETone{\mbc,J}$ does not have any infinite paths
  of thinning or quadratic type. Equivalently, all infinite paths in
  $\ADone{\mbc,J}$ have superquadratic tails.
\end{cor}

\section{Overlapping pairs and periodic mergers}\label{sec:superquadratic}

We start our attack of the superquadratic case by examining
overlapping pairs and by introducing a new move: the periodic merger
(precisely defined in Section \ref{sec:periodic-mergers}). Throughout
this section $\quotmap$ will denote the map
$\quotmap: \tilde{\mbc} \to \dualtree{t}{\mbc}$ given in Proposition
\ref{prop:tree}(\ref{eq:phi-map}).

\begin{defn}\label{defn:overlapping}
  A dual pair $\basepair\mu$ is an \define{overlapping pair} if
  $\interior{\mu \cap \ol{\mu}} \neq \emptyset$ and the pair is
  \define{orientation preserving}, i.e.  if the image of
  $\band\mu \subset \mbc$ does not contain an embedded Möbius band.
\end{defn}

\begin{conv}
  When $\basepair\mu$ is an overlapping pair we will assume that
  $\mu < \ol{\mu}$, where $<$ is the ordering on $\mbc$ (Definition
  \ref{defn-ordering}.)
\end{conv}

\begin{defn}\label{defn:tubular-class}
  Let $\band{\mu}$ be a band such that $\basepair\mu$ is an
  overlapping pair. Let $p \in \mu$ the $<$-initial point (Definition
  \ref{defn-ordering}.)  The \define{tubular loop} $\tube{\mu}$ is the loop
  $\alpha*\beta$ where $\alpha$ is the path in the side of $\band\mu$
  starting at $p$ and going from $\mu$ to $\dual\mu$ and $\beta$ is
  the path in $\mu$ connecting the endpoint of $\alpha$ to $p$. See
  Figure \ref{fig:tubular}.
\end{defn}
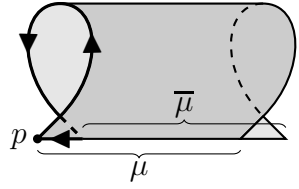
\begin{figure}[htb]
\centering
\begin{tikzpicture}[scale=0.3]
  \fill[black!20!white] (2,0) -- (1,1) .. controls (3,3) and (2.8,6)
  .. (1,6) -- (10,6) .. controls (11.8,6) and (12,3) .. (10,1) --
  (9,0) -- cycle;
  \draw[thick] (1,1) .. controls (3,3) and (2.8,6)
  .. (1,6) -- (10,6) .. controls (11.8,6) and (12,3) .. (10,1) --
  (9,0) -- (2,0);
  \draw[dashed,thick] (10,6) .. controls (8.2,6) and (8,3) .. (10,1);
  
  \filldraw[black!10!white,draw=black,thick] (9,0) -- (10,1) -- (11,0)
  --cycle;
  \filldraw[black!10!white,draw=black,thick](1,1) .. controls (3,3) and
  (2.8,6) .. (1,6) (1,6) .. controls (-0.8,6) and (-1,3) .. (1,1);
  \fill[black!20!white] (0,0) -- (1,1) -- (2,0);
  \draw[dashed,thick] (1,1) -- (2,0);
  \draw[thick] (1,1) -- (0,0) -- (2,0);
  
  \draw[very thick] (0,0) -- (1,1) (1,1) (1,1) .. controls (3,3) and
  (2.8,6) .. node[sloped]{$\blacktriangleleft$} (1,6) (1,6)
  .. controls (-0.8,6) and (-1,3) ..  node[sloped]
  {$\blacktriangleleft$}(1,1) (2,0) --
  node[sloped]{$\blacktriangleleft$} (0,0);
  \draw[dashed, very thick] (1,1) -- (2,0);
  \fill[black] (0,0) circle  (0.2);
  \node at (0,0) [left]{$p$};
  
  \draw[decorate,decoration={brace,amplitude=0.1cm}] (9,-0.5) --
  node[below] {$\mu$} (0,-0.5);
  \draw[decorate,decoration={brace,amplitude=0.1cm}] (2,0.5) --
  node[above] {$\overline{\mu}$} (11,0.5);
\end{tikzpicture}

\caption{An overlapping pair. The dark loop indicated on the left is
  the tubular loop $\tube{\mu}$.}\label{fig:tubular}
\end{figure}

We note that our definition of a tubular loop is an oriented based
loop. This gives rise to an element $g\in\fungrp\mbc$ that we will
call a \define{$\mu$-tubular element}. For the rest of the paper we
will avoid mentioning the basepoint.

\begin{figure}[htb]
\centering
\begin{tikzpicture}
  \begin{scope}[scale=0.75]
  \clip[decorate,decoration=zigzag] (-2,-2) -- (-3,5) -- (5,5) --
  (5,-2) --cycle ;

  \filldraw[black!20!white,draw=black]
  (0,0) -- (0,3) -- (3,3) -- (3,0) -- cycle
  (-1,3) -- (-1,6) -- (2,6) -- (2,3) -- cycle
  (1,0) -- (4,0) -- (4,-3) -- (1,-3) -- cycle;
  \fill[black] (0,0) circle (0.1)
  (-1,3) circle (0.1);

  \draw[thin,dashed] 
  (0.25,0) -- (0.25,3)
  (0.75,0) -- (0.75,3)
  (1.25,0) -- (1.25,3)
  (1.75,0) -- (1.75,3)
  (2.25,0) -- (2.25,3)
  (2.75,0) -- (2.75,3);
  \begin{scope}[shift={(-1,3)}]    
    \draw[thin,dashed] 
    (0.25,0) -- (0.25,3)
    (0.75,0) -- (0.75,3)
    (1.25,0) -- (1.25,3)
    (1.75,0) -- (1.75,3)
    (2.25,0) -- (2.25,3)
    (2.75,0) -- (2.75,3);
  \end{scope}  
  \begin{scope}[shift={(1,-3)}]    
    \draw[thin,dashed] 
    (0.25,0) -- (0.25,3)
    (0.75,0) -- (0.75,3)
    (1.25,0) -- (1.25,3)
    (1.75,0) -- (1.75,3)
    (2.25,0) -- (2.25,3)
    (2.75,0) -- (2.75,3);
  \end{scope}

  \draw[very thick](0,0) --
  node[above,fill=black!20!white]{$\tilde{\mu}$} +(3,0);
  \draw[very thick] (-1,3) + (0,0) --
  node[above,fill=black!20!white]{$\tube{\mu}\cdot\tilde{\mu}$} +(3,0) (0,3) +
  (0,0) -- node[below,fill=black!20!white]{$\tilde{\overline{\mu}}$}
  +(3,0);

  \draw (0,0) node[left]{$\tilde{p}$} (-1,3)
  node[left]{$\tube{\mu}\cdot\tilde{p}$};
  \end{scope}

  \begin{scope}[shift={(5,2)}]
    \draw (-1.5,0) -- (4.5,0);
    \draw[decorate,decoration={brace,amplitude=0.1cm}]
    (2,-0.4) -- node[below]{$\tube{\mu} \cdot \quotmap(\tilde\mu)$}(-1,-0.4);
    \draw[decorate,decoration={brace,amplitude=0.1cm}]
    (0,0.4) -- node[above]{$\quotmap(\tilde{\mu})$}(3,0.4);
    \draw[fill=black] (-1,0) circle (0.05)
    (-0.5,0) circle (0.05)
    (0,0) circle (0.05)
    (0.5,0) circle (0.05)
    (1,0) circle (0.05)
    (1.5,0) circle (0.05)
    (2,0) circle (0.05) 
    (2.5,0) circle (0.05)
    (3,0) circle (0.05)
    (3.5,0) circle (0.05)
    (4,0) circle (0.05)
    (-1.5,0) node[left]{\ldots} (4.5,0) node[right]{\ldots};
  \end{scope}
\end{tikzpicture}
\caption{On the left, the action of a tubular element on the universal
cover by deck transformations; on the right the corresponding action
on the tree $T(P,\complex{C})$. Tracks are shown as dashed lines}\label{fig:tubular-action}
\end{figure}
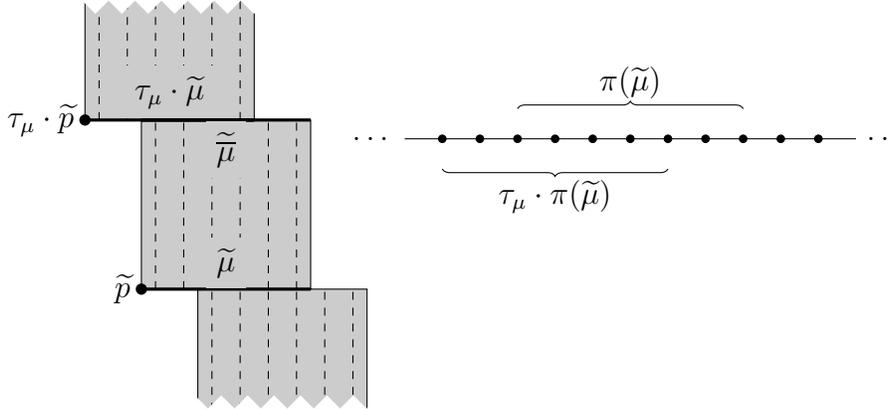

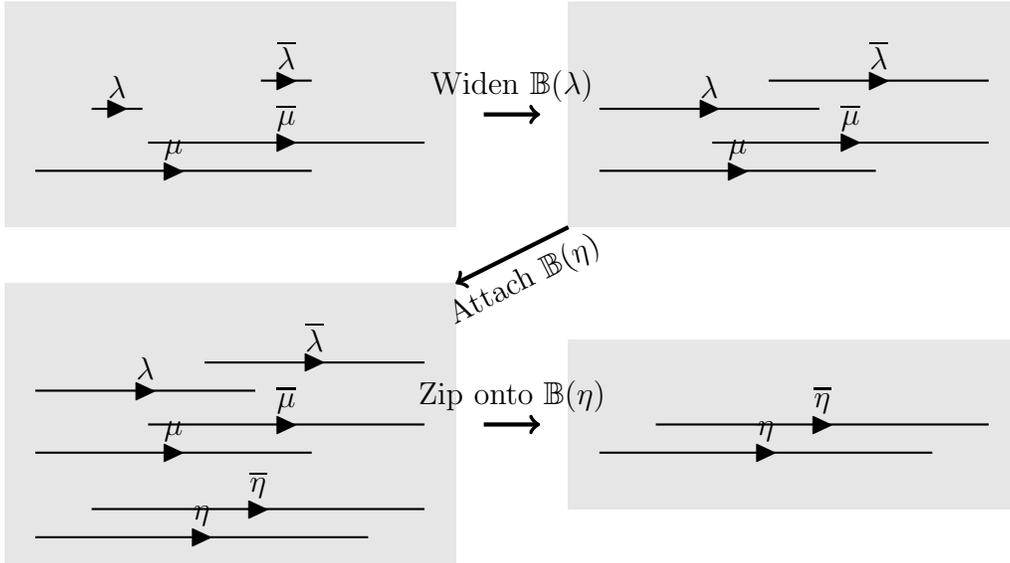
\begin{figure}[htb]
  \centering
  \begin{tikzpicture}[scale=0.75]
    \fill[black!10!white] (-0.5,-1) rectangle (7.5,3);
    \opolpair{0}{0}{5}{2}{\mu}
    \opolpair{1.1}{1}{1}{3}{\lambda}
    \begin{scope}[xshift=10cm]
      \fill[black!10!white] (-0.5,-1) rectangle (7.5,3);
      \opolpair{0}{0}{5}{2}{\mu}
      \opolpair{1.1}{0}{4}{3}{\lambda}
    \end{scope}
    \begin{scope}[yshift=-5cm]
      \fill[black!10!white] (-0.5,-2) rectangle (7.5,3);
      \opolpair{0}{0}{5}{2}{\mu}
      \opolpair{1.1}{0}{4}{3}{\lambda}
      \opolpair{-1.5}{0}{6}{1}{\eta}
    \end{scope}
    \begin{scope}[xshift=10cm,yshift=-5cm]
      \fill[black!10!white] (-0.5,-1) rectangle (7.5,2);
      \opolpair{0}{0}{6}{1}{\eta}
    \end{scope}
    \draw[ultra thick,->] (8,1) -- node[above]{Widen $\band\lambda$}
    (9,1);
    \draw[ultra thick,->] (8,-4.5) -- node[above]{Zip onto $\band\eta$} (9,-4.5);
    \draw[ultra thick,->] (9.5,-1) -- node[below,sloped]{Attach 
      $\band\eta$} (7.5,-2);
  \end{tikzpicture}
  \caption{A periodic merger (defined in Section
    \ref{sec:periodic-mergers}.) In this example
    $\tr{\tube\mu}=2, \tr{\tube\lambda}=3$ and
    $\tr{\tube\eta}=\gcd(2,3)=1.$}
  \label{fig:PM}
\end{figure}

Suppose that $\mbc$ carries a track $t$. Figure
\ref{fig:tubular-action} illustrates the action of the $\mu$-tubular
element $\tube\mu$ on $\dualtree t \mbc$ given by the quotient map
$\quotmap$ of Proposition \ref{prop:tree}. 
$p \in \mbc$ is as in Definition \ref{defn:tubular-class} and
$\tilde p \in \tilde \mu$ is a lift in the universal cover
$\tilde\mbc$ of $p \in \mu \subset \mbc$.
It is evident, for any track $t$ efficiently carried by $\mbc$, that
$\tube\mu$ must fix some axis of $\dualtree t \mbc$. Closer
examination immediately yields:

\begin{lem}\label{lem:tubular-description}
  If $\mbc$ efficiently carries a track $t$ and $\basepair\mu$
  is an overlapping pair then the $\mu$-tubular element $\tube\mu$ acts on
  $\dualtree t \mbc$ hyperbolically with
  translation length denoted
  \[ \tr\mu = |\mu \setminus \dual\mu|_t. \]
\end{lem}
\begin{figure}[htb]
  \centering
  \begin{tikzpicture}
      \draw[ultra thick, decorate,decoration={brace,amplitude=0.1cm}] (5,0) --
  node[below] {$\inter\mu$} (0,0);
    \opolpair {0.5}{0}{4}{1}{\mu}
    \draw[decorate,decoration={brace,amplitude=0.1cm}] (0,1.2) --
    node[above]{$\tr\mu$} (1,1.2);
  \end{tikzpicture}
  \caption{The translation length $\tr{\mu}$ and $\inter\mu$, the
    section corresponding to $\mu$ given in Definition
    \ref{defn:overlapping-interval}.}
  \label{fig:tr}
\end{figure}
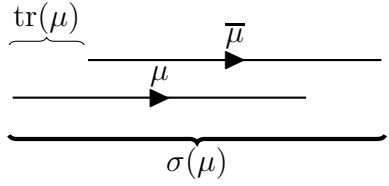
One of the principal features of a group acting acylindrically on a
tree is that infinite line stabilizers are cyclic. It follows that if
there are two overlapping pairs that themselves overlap sufficiently,
the corresponding tubular elements must fix a common axis, and
therefore must lie in a common cyclic subgroup. In Section
\ref{sec:periodic-mergers} we will describe the periodic merger, a
move from \cite{TAH-2006}, which will replace these two overlapping
pairs by a single overlapping pair. This is illustrated in Figure
\ref{fig:PM}.
There are two subtleties to this operation, which do not occur when
deciding if an equation has a solution, or if we simply want to count
orbits as in \cite{TAH-2006}.
\begin{itemize}
\item We must merge these bands into one while preserving the
  fundamental group of the band complex as well as the dual Bass-Serre
  tree.
\item We must be able to \emph{algorithmically} produce a finite list containing all
  combinatorial outcomes of a periodic mergers for $(\mbc,t)$, where
  $t$ ranges over $\kappaTracks{\mbc}$. (Definition
  \ref{defn:kappaTracks}.)
\end{itemize}
There is one outstanding difficulty: given two commuting elements
$g,h \in \fungrp{\mbc}$, decide if they lie in a common cyclic
subgroup. Only being able to solve the word problem in $\fungrp\mbc$
is insufficient to solve this problem in general. To overcome this
impasse we will use the author's generalized Bulitko Lemma
\cite{touikan-bulitko} in a way that is completely different from its
usual purpose.

\subsection{Interactions with tubular
  elements: entanglement}\label{sec:tubularinteraction}

\begin{defn}\label{defn:overlapping-interval}
  If $(\mu,\ol{\mu})$ is an overlapping pair then we denote the
  \define{section corresponding to $(\mu,\ol{\mu})$} as \[ \inter{\mu}
  = \mu \cup \dual{\mu}. \]
\end{defn}

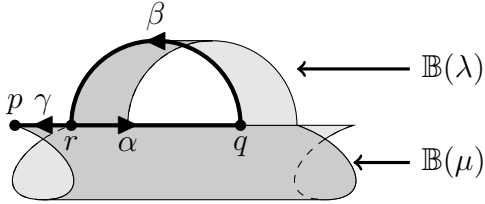
\begin{figure}[htb]
  \begin{center}
    \begin{tikzpicture}[scale=0.75]
      \begin{scope}[yscale=-0.66]
        \olband{0}{0}{5}{1} 
        \draw[very thick, ->] (7,1)
        node[right]{$\band\mu$} -- (6,1);
      \end{scope}
      \opband{1}{0}{1}{3}
      \draw[ultra thick](0,0) --node{$\blacktriangleright$} node[below]{$\alpha$} (4,0);
      \opside{1}{0}{3}
      \draw (2.5,1.5) node{$\blacktriangleleft$} node[above]{$\beta$};
      \draw (0.5,0) node{$\blacktriangleleft$} node[above]{$\gamma$};
      \fill[black] (0,0)node[above]{$p$} circle (0.1);
      \fill[black] (1,0)node[below]{$r$} circle (0.1);
      \fill[black] (4,0)node[below]{$q$} circle (0.1);
      \draw[very thick, ->]  (7,1) node[right]{$\band\lambda$} -- (5,1);
    \end{tikzpicture}
  \end{center}
  \caption{The band $\band\lambda$ has both bases lying in the section
    $\inter\mu$. The $\mu$-relative loop $\relloop\mu\lambda$ is the
    loop $\alpha*\beta*\lambda$.}
  \label{fig:relloop}
\end{figure}

Suppose now that here is another band $\band\lambda$ whose 
unmatched bases both lie in $\inter\mu$, as in Figure
\ref{fig:relloop}.
\begin{conv}
  We write $\basepair\lambda \subset \inter\mu$ to signify
  $\lambda \cup \dual\lambda \subset \inter\mu$. We will always assume
  that that base pairs $\basepair\lambda$, $\basepair\mu$ are
  unmatched.
\end{conv}
To study how these bands interact we have the following

\begin{defn}\label{defn:relative-class}
  Let $\basepair\mu$ be overlapping and let $\basepair\lambda  \subset \inter{\mu}$ with $\lambda < \ol{\lambda}$. Let $p$ be as in
  Definition \ref{defn:tubular-class}. The \define{$\mu$-relative
    loop} $\relloop{\mu}{\lambda}$ is the concatenation
  $\alpha*\beta*\gamma$ shown in Figure \ref{fig:relloop} where
  $\alpha$ is the path from $p$ to $q$ in $\inter\mu$, $\beta$ is the
  path from $q$ to $r$ travelling along a side of $\band\lambda$, and
  $\gamma$ is the path in $\inter\mu$ from $r$ to $p$.
\end{defn}

\begin{figure}[htb]
  \centering
  \begin{tikzpicture}[scale=0.75]
    \draw[thick] (0,0) -- (6,0);
    \draw[decorate,decoration={brace,amplitude=0.1cm}] (6,-0.5) --node[below]{$\tilde{\inter\mu}$} (0,-0.5);
    \draw[fill=black!20!white] (4,0) -- (4,3) -- (5,3) -- (5,0)
    --cycle;
    \draw[ultra thick] (4,0)
    --node[sloped]{$\blacktriangleright$}node[left]{$\tilde\beta$}
    (4,3);
      \fill (4,0) circle (0.1);
      \draw (4,0) node[below]{$q$};
      \fill (0,0) circle (0.1);
      \draw (0,0) node[above]{$\tilde p$};
      \draw[ultra thick] (0,0) -- node{$\blacktriangleright$}node[above]{$\tilde\alpha$} (4,0);
    \begin{scope}[xshift=3cm,yshift=3cm]
     \draw[thick] (0,0) -- (6,0);
     \fill (0,0) circle (0.1);
     \draw (0,0) node[left]{$\relloop\mu\lambda\cdot\tilde p$};
     \draw[ultra thick] (1,0) -- node{$\blacktriangleleft$} node[above]{$\tilde{\gamma}$}  (0,0);
     \fill (1,0) circle (0.1);
     \draw (1,0) node[above]{$\tilde r$};
    \end{scope}
    
  \end{tikzpicture}
  \caption{The deck transformation corresponding to
    $\relloop\mu\lambda$. The lifts of $\alpha,\beta$ and $\gamma$ of
    Figure \ref{fig:relloop} are shown.}
  \label{fig:relloop-deck}
\end{figure}
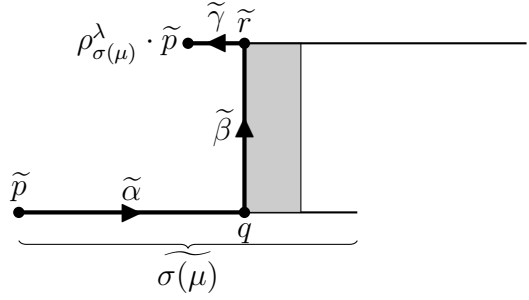

Again fixing a lift $\tilde{p}$ of $p$ and $\tilde{\mu}$ of $\mu$
in the universal cover $\tilde\mbc$ of $\mbc$, we can describe the
deck transformation given by $\relloop\mu\lambda$ (see Figure
\ref{fig:relloop-deck}.)
From this we immediately get:
\begin{lem}\label{lem:relative-loop-intersection}
  If $\mbc$ efficiently carries a track $t$, then in $\dualtree t
  \mbc$ the length of the arc\[
  \axis(\tube\mu) \cap \left(\relloop\mu\lambda \cdot \axis(\tube\mu)\right)
  \] is at least $|\lambda|_t$.
\end{lem}

\begin{defn}\label{defn:rel-tr}
  We define $\reltr{\mu}{\lambda}$ to be the measure of the arc
  between the leftmost point of $\lambda$ and the leftmost point of
  $\ol{\lambda}$, i.e. the length $|\alpha|_t-|\gamma|_t$ as shown in
  Figure \ref{fig:relloop}. We say the dual pair $\basepair \lambda$
  is \define{orientation preserving} if the holonomy $\lambda \to
  \dual\lambda$ extends to an orientation preserving homeomorphism of
  $\inter\mu$ (relative to its endpoints.)
\end{defn}

For the following three lemmas assume that $t \in \kappaTracks{\mbc}$.

\begin{defn}\label{defn:entangled1}
  If $\basepair\mu$ is an overlapping pair, $\basepair\lambda \subset \inter\mu$ with $\basepair\lambda$ orientation preserving,
  and the commutator $[\tube\mu,\relloop\mu\lambda]=1$, then we say
  that dual pairs $\basepair\mu$ and $\basepair\lambda$ are
  \define{entangled}.
\end{defn}

The following is obvious, but necessary for computational
considerations:

\begin{lem}\label{lem:compute-entanglement}
  Let $\mbc$ be a band complex with $\basepair\mu$ be an overlapping
  pair, and $\basepair\lambda \subset \inter\mu$ and such that
  $\basepair\lambda$ is orientation preserving. If we can solve the
  word problem in $\fungrp{\mbc}$ then we can decide if $\basepair\mu$
  and $\basepair\lambda$ are entangled.
\end{lem}

\begin{lem}\label{lem:commuting-rel-class}
  Let $(\mu,\ol{\mu})$ be an overlapping pair, suppose
  $\basepair\lambda \subset \inter{\mu}$ and that
  $(\lambda,\ol{\lambda})$ is orientation preserving and unmatched. If
  $[\tube{\mu},\relloop{\mu}{\lambda}]=1$, i.e. $\basepair\mu$ and
  $\basepair\lambda$ are entangled, then $\relloop{\mu}{\lambda}$ acts
  hyperbolically on $\dualtree t \mbc$,
\[\axis(\relloop{\mu}{\lambda}) = \axis(\tube{\mu}),\] and the
translation length of $\relloop{\mu}{\lambda}$ is
$\reltr{\mu}{\lambda}$.
\end{lem}

\begin{proof}
  Because $[\tube{\mu},\relloop{\mu}{\lambda}]=1$,
  \[\relloop\mu\lambda \cdot \axis(\tube\mu) =
  \axis\left(\relloop{\mu}{\lambda} \tube\mu
    (\relloop{\mu}{\lambda})^{-1}\right) = \axis(\tube\mu).\] It
  therefore follows that $\bk{\relloop\mu\lambda}$ fixes a bi-infinite
  arc in $\dualtree t \mbc$. From Figure \ref{fig:relloop-deck},
  $\relloop\mu\lambda$ translates this arc by
  $\tr{\relloop\mu\lambda}$ (as defined in Definition
  \ref{defn:rel-tr}), it therefore follows $\axis(\tube \mu)$ is the
  minimal invariant subtree for $\bk{\relloop\mu\lambda}$.
\end{proof}

\begin{lem}\label{lem:small-shifts-commute}
  Let $(\mu,\ol{\mu})$ be an overlapping pair, let $\basepair\lambda \subset \inter{\mu}$, and let $(\lambda,\ol{\lambda})$
  be orientation preserving. If $|\lambda|_t \geq \tr{\mu} + \kappa
  +1$ then $[\tube{\mu},\relloop{\mu}{\lambda}]=1$, i.e.
  $\basepair\mu$ and $\basepair\lambda$ are entangled.
\end{lem}

\begin{proof}
  Let \[ I = \axis(\tube{\mu }) \cap \axis \big(
  \relloop{\mu}{\lambda} \tube{\mu} (\relloop{\mu}{\lambda})^\mo
  \big). \] By Lemma \ref{lem:relative-loop-intersection} we have
  $|I|\geq \tr{\mu} + \kappa + 1$ let $J$ be a co-initial (or
  co-final) subsegment of $I$ with$|J| \geq \kappa + 1$. Without loss of generality (up
  to choosing if $J$ is co-initial or co-final) we may assume
  that \[\tube{\mu}\cdot J \subset I \subset \axis \big(
  \relloop{\mu}{\lambda} \tube{\mu} (\relloop{\mu}{\lambda})^\mo
  \big)\] and since $\lambda$ is orientation preserving and since
  translation length in invariant under conjugation we have :\[
  \big(\relloop{\mu}{\lambda} \tube{\mu}^\mo
  (\relloop{\mu}{\lambda})^\mo \big) \cdot (\tube{\mu} \cdot J) = J\]
  which by $\kappa$-acylindricity implies
  $[\tube{\mu},\relloop{\mu}{\lambda}]=1$.
\end{proof}

\subsection{The Bulitko trick}
\label{sec:bulitko-trick}

If we are given a band complex $\mbc$ with an overlapping pair
$\basepair\mu$ that is entangled with $\basepair\lambda$ then, by Lemma
\ref{lem:commuting-rel-class}, if $\kappaTracks \mbc \neq \emptyset$,
then $\bk{\tube\mu,\relloop\mu\lambda}$ must be cyclic. In particular
there is
$\bk{g} = \bk{\tube\mu,\relloop\mu\lambda}$ so that $g^n = \tube\mu$
and $g^m=\relloop\mu\lambda$ with $n,m$ minimal in absolute value. The
Bulitko trick will either enable us to find $n,m$ or will certify that
$\kappaTracks\mbc = \emptyset$.

Let $\phi:\fungrp\mbc \to \mathbb H$ be a homomorphism to the
fundamental group of a one-edged $\kappa$-acylindrical graph of
groups. Then as long as some element of $\bk{g}$ is sent to a
hyperbolic element $n,m$ will still be the minimal integers such that
$\phi(g)^n = \phi(\tube\mu)$ and $\phi(g)^m=\phi(\relloop\mu\lambda)$.

Since $[\relloop \mu \lambda,\tube\mu]=1$ then any $\phi:\fungrp\mbc
\to \mathbb H$ will send $\tube\mu$ and $\relloop\mu\lambda$ to
elements fixing a common axis, provided their images are
hyperbolic. In this case, the acylindrical Bulitko Lemma provides a
computable function depending only on $\tube\mu$ and
$\relloop\mu\lambda$ that bounds $n,m$ given above. Before giving the
lemma we must first present the necessary terminology.

Let $\mathbb{H}$ act $\kappa$-acylindrically on a based tree
$(T,t_0)$, let $p$ be some hyperbolic element of $\mathbb{H}$
(with respect to the given splitting of $\mathbb{H}$) and let $L\subset T$
denote its axis.  Consider the set of segments \[\big\{L' \subset L
\mid L = \bigcup_{n \in \Z} p^n L'\big\}. \] A minimal element of this
set with respect to inclusion is called a \define{fundamental domain} of $L$.

For an element $h \in \mathbb{H}$, let $[v_0,h\cdot v_0]$ denote the
geodesic between $v_0$ and $h\cdot v_0$. Let $g \in \mathbb{H}$. If a
segment
\[\sigma = [v_0,h\cdot v_0] \cap g\cdot L\] is non-empty then we
call it an $L$-periodic subsegment of $[v_0,h\cdot v_0]$. The
\define{$L$-periodicity} of $\sigma$ is the integer\[ \lfloor
\frac{|\sigma|}{|L_0|}\rfloor\] where $L_0$ is a
fundamental domain of $L$. We can now state

\begin{thm}[Acylindrical Bulitko Lemma, Theorem 1.3 of
  \cite{touikan-bulitko}]\label{lem:bulitko}
  There exists a computable function $\mathfrak{n}:\N \times \N \times
  \N \rightarrow \N$ such that for any nontrivial homomorphism $\phi:
  G \rightarrow \mathbb{H}$; where the group $G$ has a finite
  presentation $\bk{Y \mid S}$ and the group $\mathbb{H}$ has a
  $\kappa$-acylindrical splitting with based Bass-Serre tree
  $(T,t_0)$; and for any hyperbolic element in $p\in \mathbb{H}$
  (denote its axis $L \subset T)$, there exists a homomorphism
  $\phi^*: G \rightarrow \mathbb{H}$ such that for all $y \in Y$
  \begin{itemize}
  \item if $[t_0,\phi(y)\cdot t_o]$ has no $L$-periodic subsegments, then
    $\phi(y)=\phi^*(y)$, and
  \item if $[t_0,\phi(y)\cdot t_o]$ has $L$-periodic subsegments, then
    there is a bijective correspondence between the $L$-periodic
    subsegments of $[t_0,\phi(y)\cdot t_o]$ and $[t_0,\phi(y)^*\cdot
    t_o]$, but the $L$-periodicity of all the periodic subsegments of
    $[t_0,\phi(y)^*\cdot t_o]$ is at most
    $\mathfrak{n}(|Y|,|S|,\kappa)$.
  \end{itemize}
\end{thm}

\begin{prop}\label{prop:Bulitko-trick}
  Let $\basepair\mu$ be an overlapping pair in a band complex $\mbc$
  and suppose it is entangled with $\basepair\lambda$. There is an
  algorithm which terminates with one of the two following
  outputs:
  \begin{enumerate}[(i)]
  \item\label{it:bulitko-win} It gives an element $g$ such that $\bk{g} =
    \bk{\tube\mu,\relloop\mu\lambda}$.
  \item\label{it:bulitko-fail} It (correctly) certifies that
    $\kappaTracks \mbc = \emptyset$.
  \end{enumerate}
\end{prop}
\begin{proof}
  Let $G=\fungrp\mbc$. Since $\mbc$ is an explicitly given cell
  complex and since we can solve the word problem in $\fungrp\mbc$, it
  is possible to give a finite presentation $\bk{Y\mid S}$ of $G$
  where $\tube\mu$ and $\relloop\mu\lambda$ are included in the
  generating set $Y$.

  Let $\mathfrak n$ be the computable function given by Theorem
  \ref{lem:bulitko} and let $M=\mathfrak{n}(|Y|,|S|,\kappa)$. For
  every pair $n_i,m_j$ of absolute value less than $M$ let $(u_i,v_j)$
  be a pair such that $u_in_i + v_im_j = \gcd(n_i,m_i)$. Let
  $g_{ij} = (\tube\mu)^{u_i}(\relloop\mu\lambda)^{v_j}$ and check
  whether $(g_{ij})^{n_i} = \tube\mu$ and
  $(g_{ij})^{m_j}=\relloop\mu\lambda$. If we find some $g_{ij}$
  satisfying item (\ref{it:bulitko-win}) then we stop. Otherwise if
  all these verifications were negative we know that that
  $\kappaTracks{\mbc}=\emptyset.$

  Indeed suppose towards a contradiction that none of the $g_{ij}$
  were roots of $\tube\mu$ and $\relloop \mu\lambda$ but that there is
  some $t\in \kappaTracks \mbc$. By Lemma
  \ref{lem:commuting-rel-class} there is some $g \in G$ such that
  $\bk{g}$ is the maximal cyclic group stabilizing
  $L=\axis{(\tube\mu)} = \axis{\left(\relloop\mu\lambda\right)}$ let
  $n,m$ be the integers such that $g^n=\tube\mu$ and
  $g^m=\relloop\mu\lambda$. By hypothesis $|n|$ or $|m|$ is greater
  than $M$. Let $\mathbb H = \fungrp\mbc$ let $T$ be
  $\dualtree t \mbc$, let $\phi$ be the identity, let $p=g$ and let
  $L$ denote the axis of $g$. We can chose the basepoint $t_0$ of $T$
  so that $t_0 \in L$ which implies that that the geodesics
  $[t_0,\phi(\relloop\mu\lambda)\cdot t_0]$ and
  $[t_0,\phi(\tube\mu)\cdot t_0]$ consist of a single $L$-periodic
  segment. By Theorem \ref{lem:bulitko} there exists an endomorphism
  $\phi^*$ such that $[t_0,\phi^*(\relloop\mu\lambda)\cdot t_0]$ and
  $[t_0,\phi^*(\tube\mu)\cdot t_0]$ contain a single non-trivial
  $L$-periodic segment, this implies the hyperbolicity of
  $\phi^*(\tube\mu)$ and $\phi^*(\relloop\mu\lambda)$; thus the
  restriction $\phi^*|_\bk{g}$ is injective. On the other hand the
  bound on $L$ periodicity implies that $\phi^*(\tube\mu) = g^{n_0}$
  and $\phi^*(\relloop\mu\lambda) = g^{m_0}$ with
  $|n_0|, |m_0| \leq M$. Now we must have $\phi^*(g)=g^r$ which implies
  that $|n_0| = |r||n|$ and $|m_0|=|r||m|$ which contradicts the
  assumption that $|n|,|m| > M.$
\end{proof}

\subsection{Merging entangled pairs}\label{sec:periodic-mergers}
Let $\mbc$ be a measured band complex and let
$t \in \kappaTracks \mbc$. Suppose we have an overlapping pair
$\basepair\mu$ entangled with $\basepair\lambda \subset
\inter\mu$.
Then it will sometimes be possible to merge the bands $\band\lambda$ and
$\band\mu$ into a new band $\band\eta$.

Simply attaching a new band with overlapping bases
to $\mbc$ inside the segment $\inter\mu$ will add a cyclic free factor
to $\fungrp\mbc$. We must therefore also attach a 2-cell to encode that
$\tube\eta$ is a root of $\tube\mu$ and $\relloop\mu\lambda$. In order
to do so we may first have to widen $\band\lambda$.

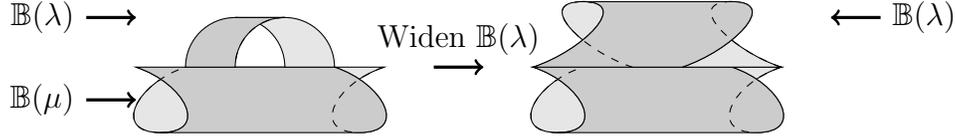
\begin{figure}[htb]
  \begin{center}
    \begin{tikzpicture}[scale=0.66]
      \begin{scope}[yscale=-0.66]
        \olband{0}{0}{4}{1} 
        \draw[very thick, ->] (-1,1)
        node[left]{$\band\mu$} -- (0,1);
      \end{scope}
      \opband{1}{0}{1}{2}
     
      \draw[very thick, ->]  (-1,1) node[left]{$\band\lambda$} --
      (0,1);
      \draw[very thick,->] (6,0)--node[above]{Widen $\band\lambda$}(7,0);
      \begin{scope}[xshift=8cm]
        \begin{scope}[yscale=-0.66]
        \olband{0}{0}{4}{1} 
      \end{scope}
      \begin{scope}[yscale=0.66]
        \olband{0}{0}{3}{2}    
      \end{scope}
      \draw[very thick, ->]  (7,1) node[right]{$\band\lambda$} --
      (6,1);
      \end{scope}
      
    \end{tikzpicture}
  \end{center}
  \caption{We widen the band $\band\lambda$ so that it becomes
    coinitial and coterminal with $\inter\mu$.}
  \label{fig:widen-downstairs}
\end{figure}

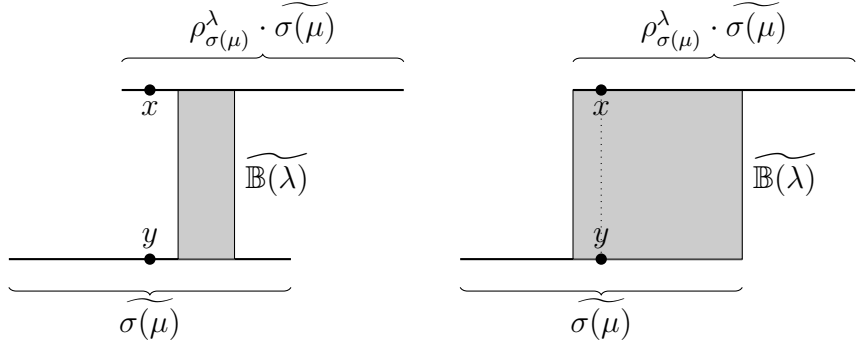
\begin{figure}[htb]
  \centering
    \begin{tikzpicture}[scale=0.75]
    \draw[thick] (1,0) -- (6,0);
    \draw[decorate,decoration={brace,amplitude=0.1cm}] (6,-0.5) --node[below]{$\tilde{\inter\mu}$} (1,-0.5);
    \draw[fill=black!20!white] (4,0) -- (4,3) -- (5,3) --node[right]{$\tilde{\band\lambda}$} (5,0)
    --cycle;
    \fill (3.5,0) circle (0.1) node[above]{$y$};

    \begin{scope}[xshift=2cm,yshift=3cm]
     \draw[thick] (1,0) -- (6,0);
     \draw[decorate,decoration={brace,amplitude=0.1cm}] (1,+0.5) --node[above]{$\relloop\mu\lambda\cdot\tilde{\inter\mu}$} (6,+0.5);
     \fill (1.5,0) circle (0.1) node[below]{$x$}; 
   \end{scope}
    
    \begin{scope}[xshift=8cm]
      \draw[thick] (1,0) -- (6,0);
      \draw[decorate,decoration={brace,amplitude=0.1cm}] (6,-0.5) --node[below]{$\tilde{\inter\mu}$} (1,-0.5);
      \draw[fill=black!20!white] (3,0) -- (3,3) -- (6,3) --node[right]{$\tilde{\band\lambda}$} (6,0)
      --cycle;
      \draw[dotted] (3.5,0) -- (3.5,3);
      \fill (3.5,0) circle (0.1) node[above]{$y$};
      \begin{scope}[xshift=2cm,yshift=3cm]
        \draw[thick] (1,0) -- (6,0);
     \draw[decorate,decoration={brace,amplitude=0.1cm}] (1,+0.5) --node[above]{$\relloop\mu\lambda\cdot\tilde{\inter\mu}$} (6,+0.5);
     \fill (1.5,0) circle (0.1) node[below]{$x$}; 
   \end{scope}
 \end{scope}

  \end{tikzpicture}
  \caption{The widening move corresponds to equivariantly widening
    every lift $\tilde{\band\lambda}$ of $\band\lambda$ in the
    universal cover. If $x,y$ lie in the same track $\tilde t$ in
    $\tilde{\mbc}$, then in $\tilde{\mbc'}$ we extend the track to
    pass through the widened $\tilde{\band\lambda}$.}
  \label{fig:widen-upstairs}
\end{figure}

Consider the operation of widening a band $\band\lambda$ illustrated
in Figure \ref{fig:widen-downstairs}.
We do this so that the resulting base $\lambda$ is coinitial with
$\inter\mu$ and $\dual\lambda$ is cofinal. The inverse of a widening
is a deformation retraction, so it preserves
$\fungrp\mbc$. Furthermore the element $\relloop\mu\lambda$ of the
fundamental group is unchanged. If $[\relloop\mu\lambda,\tube\mu]=1$
then for any track $t$ efficiently carried by $\mbc$ both
$\relloop\mu\lambda$ and $\tube\mu$ have the same axis by Lemma
\ref{lem:commuting-rel-class}. Figure \ref{fig:widen-upstairs} depicts
what happens when we pass to the universal cover.
Because $\relloop\mu\lambda$ has the same axis as $\tube\mu$ there is
a natural way to extend the pattern in $\tilde{\mbc}$ to a pattern of
$\tilde{\mbc'}$ so that the resulting dual trees $\dualtree t \mbc,
\dualtree{t'}{\mbc'}$ are equivariantly isomorphic. No 2-cells were
added so $t'$ is efficiently carried. We summarize in
the following lemma.

\begin{lem}\label{lem:widen}
  Let $\mu$ and $\lambda$ be as above and let
  $t \in \kappaTracks\mbc$. We can widen $\band{\lambda}$ so that, in
  the new band complex $\mbc'$, $\lambda$ is coinitial with
  $\inter\mu$ and $\ol\lambda$ is cofinal with
  $\inter\mu$. Furthermore the dual trees $\dualtree t \mbc$ and
  $\dualtree{t'}{\mbc'}$ are equivariantly isomorphic. In particular
  $\tr\lambda,\tr\mu$ remain invariant, and $t' \in \kappaTracks{\mbc'}$.
\end{lem}

\begin{lem}\label{lem:both-overlap-wide}
  Let $\basepair\mu$ be an overlapping pair entangled with
  $\basepair\lambda$ and suppose furthermore that $\basepair\lambda$
  is itself overlapping. Then after widening $\band\lambda$ as in
  Lemma \ref{lem:widen} we have $|\lambda|_t\geq \tr{\mu}.$ 
\end{lem}
\begin{proof}
  If $\lambda$ is coinitial and cofinal with $\inter\mu$, but
  $|\lambda|_t<\tr\mu$, then $\lambda$ can't overlap with its dual.
\end{proof}
Before continuing we need the following fact, which follows by
meditating on the Euclidean algorithm.

\begin{lem}\label{lem:number-theory}
  Let $n,m$ be positive integers and let $d = \gcd(n,m)$ then without
  loss of generality there are integers $u,v \in \Z_{\geq 0}$ such
  that $d = un-vm$, moreover we have non-decreasing sequences of
  integers $0=v_0 \leq v_1 \leq \ldots \leq v_{u+v} = v$ and
  $0=u_0 \leq u_1 \leq \ldots \leq u_{u+v} = u$ with
  \[u_i + v_i +1 = u_{i+1}+v_{i+1}\]
  such that the following inequalities
  hold \begin{equation}\label{eqn:partial-sums} 0 \leq u_in + v_im
    \leq m+n\end{equation}
\end{lem}

This fact motivates the following observation.

\begin{lem}\label{lem:premerge-I}
  Let $\basepair\mu$ be an overlapping pair entangled with
  $\basepair\lambda$. If $|\lambda|_t\geq\tr{\mu}$, then
  \[|\inter{\mu}| \geq \tr{\mu} + \reltr{\mu}{\lambda}\]
\end{lem}
\begin{proof}
  By hypothesis $|\lambda|_t \geq \tr{\mu}$ and \[\lambda \cup
  \ol{\lambda} \subset \inter{\mu} \Rightarrow
  |\lambda|_t+\reltr{\mu}{\lambda}| \leq |\inter{\mu}|\] which give the
  required inequality.
\end{proof}

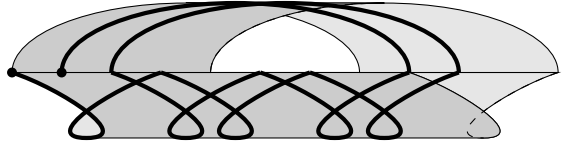
\begin{figure}[htb]
  \centering
  \begin{tikzpicture}[scale=0.66]
    \begin{scope}[yscale=-0.66]
      \olband{0}{0}{8}{3} 
    \end{scope}
    \begin{scope}[yscale=0.4]
      \opband{0}{0}{4}{7}    
    \end{scope}    

    \begin{scope}
      \begin{scope}[yscale=0.4]
        \opside{2}{0}{7}
        \opside{1}{0}{7}
      \end{scope}
      \begin{scope}[yscale=-0.66]
        \olside{0}{0}{8}{3}
        \olside{3}{0}{8}{3}
        \olside{6}{0}{8}{3}        
        \olside{2}{0}{8}{3}        
        \olside{5}{0}{8}{3}        
      \end{scope}
    \end{scope}
    \fill (0,0) circle (0.1);
    \fill (1,0) circle (0.1);
  \end{tikzpicture}
  \caption{The vertical path $\alpha_{\mu\lambda}$ contained in bands
    $\band\lambda,\band\mu$, corresponding to the series
    $3+3+3-7+3+3-7=1=\gcd(3,7)$.}
  \label{fig:merge-eg-1}
\end{figure}

Suppose we are in the situation of Lemma \ref{lem:premerge-I} and that
we have widened $\band\lambda$ as in Lemma \ref{lem:widen}. We will
illustrate the attachment of $\band\eta$ with a concrete
example. Suppose that $\tr\mu=3$, $\reltr\mu\lambda=7$, $|\mu|_t = 8$ and
$|\lambda|_t=4$, see Figure \ref{fig:merge-eg-1}.
$\gcd(3,7)=1$, we pick the linear combination $5*3-2*7=1$ which we
rewrite the series $3+3+3-7+3+3-7=1$, we do this because each initial
subsum is positive and at most $10=3+7$ as in Lemma
\ref{lem:number-theory}(\ref{eqn:partial-sums}). This Lemma implies
that such a series can be found for any pair of bands satisfying our
hypotheses. Now $\alpha_{\mu\lambda}$ in Figure \ref{fig:merge-eg-1}
is a simple path that is a concatenation of vertical sets prescribed
by the associated series. Its endpoints in $\inter\mu$ have distance
$\gcd(\tr\mu,\reltr\mu\lambda)$, which in our example is 1. By Lemmas
\ref{lem:number-theory} and \ref{lem:premerge-I}, we can always find
such an $\alpha_{\mu\lambda}$ for any pair of bands that satisfy the
hypotheses. This is why we needed to widen $\band\lambda$ in the first
place; if it were too narrow we wouldn't be able to construct
$\alpha_{\mu\lambda}$.

In $\fungrp\mbc$ the elements $\tube\mu$ and $\relloop\mu\lambda$
commute, so the product implied by the path
$\alpha_{\mu\lambda}$ is\[
(\tube\mu)^3(\relloop\mu\lambda)^\mo(\tube\mu)^2(\relloop\mu\lambda)^\mo
= (\tube\mu)^5(\relloop\mu\lambda)^{-2}.
\]
$(\tube\mu)^5(\relloop\mu\lambda)^{-2}$ is an element that translates
$\axis(\tube\mu)=\axis(\relloop\mu\lambda)$ by a distance of
$\gcd(\tr\mu,\reltr\mu\lambda)$. Noting that (signed) translation length
gives an embedding from the stabilizer of a bi-infinite line in a
$\kappa$-acylindrical tree to $\mathbb Z$, we conclude that
$g=(\tube\mu)^5(\relloop\mu\lambda)^{-2}$ is the element such that
$\bk{g} = \bk{\tube\mu,\relloop\mu\lambda}$ that will (up to sign
$\pm 1$) be produced by the algorithm given in Proposition
\ref{prop:Bulitko-trick}.

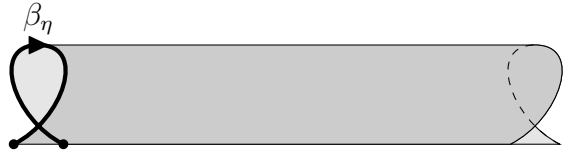
\begin{figure}[htb]
  \centering
  \begin{tikzpicture}[scale=0.66]
    \olband{0}{0}{10}{1}
    \olside{0}{0}{10}{1}
    \draw (0.5,2)node{$\blacktriangleright$} node[above]{$\beta_\eta$};
    \fill (0,0) circle (0.1);
    \fill (1,0) circle (0.1);
  \end{tikzpicture}
  \caption{Adding $\band\eta$ and the arc $\beta_\eta$.}
  \label{fig:merge-add-eta}
\end{figure}

We now attach a new band $\band\eta$ so that $\inter\eta = \inter\mu$
and that $\tr\eta = \gcd(\tr\mu,\reltr\mu\lambda)$. In our example this
forces $|\eta|_t = 10$. Let $\beta_\eta$ be the path shown in Figure
\ref{fig:merge-add-eta}.
Simply attaching $\band\eta$ to $\mbc$ gives the fundamental group
$\fungrp\mbc *\bk{\tube\eta}$. We also attach a 2-cell $B$ along
the simple closed path $\alpha_{\mu\lambda}*\beta_\eta$ to get the
resulting fundamental group
\[ \left(\fungrp\mbc *\bk{\tube\eta}\right)/\ncl{\tube\eta = g}
\approx \fungrp\mbc
\] by defining $\tube\eta = g$. First note that the new band complex
$\mbc'$ can be seen as containing $\mbc$ and that the track
$t \subset \mbc \subset \mbc'$ naturally extends to a track
$t' \subset \mbc'$ and that the trees
$\dualtree t \mbc \approx \dualtree{t'}{\mbc'}$ are equivariantly isomorphic. Further
note that, by the way the attaching map
$\alpha_{\mu\lambda}*\beta_\eta$ is defined, $\mbc'$ carries $t'$
efficiently.

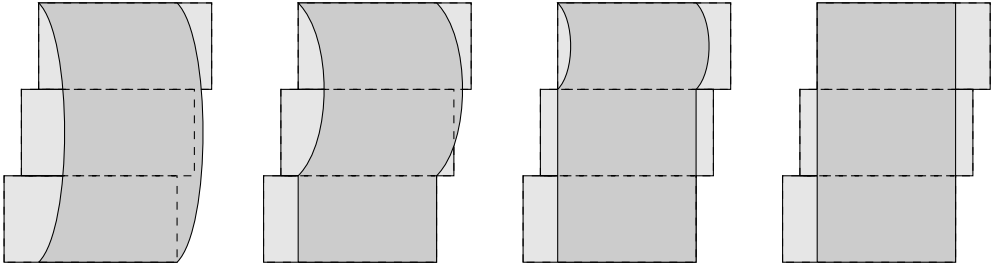
\begin{figure}[htb]
  \centering
  \begin{tikzpicture}[scale=0.23]
    \filldraw[black!10!white,draw=black] (0,0) -- (0,5) -- (10,5) --
    (10,0) -- cycle (1,5) -- (1,10) -- (11,10) -- (11,5) -- cycle
    (2,10) -- (2,15) -- (12,15) -- (12,10) -- cycle;
    \filldraw[black!20!white,draw=black] 
    (2,0) .. controls +(2,2) and +(2,-2) .. (2,15) -- (10,15) .. 
    controls +(2,-2) and +(2,2) .. (10,0) -- (2,0);
     \draw[dashed] (0,0) -- (0,5) -- (10,5) --
    (10,0) -- cycle (1,5) -- (1,10) -- (11,10) -- (11,5) -- cycle
    (2,10) -- (2,15) -- (12,15) -- (12,10) -- cycle;
    
    \begin{scope}[xshift=15cm]
       \filldraw[black!10!white,draw=black] (0,0) -- (0,5) -- (10,5) --
    (10,0) -- cycle (1,5) -- (1,10) -- (11,10) -- (11,5) -- cycle
    (2,10) -- (2,15) -- (12,15) -- (12,10) -- cycle;
    \filldraw[black!20!white,draw=black] 
    (2,0) -- (2,5) .. controls +(2,2) and +(2,-2) .. (2,15) -- (10,15) .. 
    controls +(2,-2) and +(2,2) .. (10,5) -- (10,0) -- (2,0);
     \draw[dashed] (0,0) -- (0,5) -- (10,5) --
    (10,0) -- cycle (1,5) -- (1,10) -- (11,10) -- (11,5) -- cycle
    (2,10) -- (2,15) -- (12,15) -- (12,10) -- cycle;
    \end{scope}

    \begin{scope}[xshift=30cm]
       \filldraw[black!10!white,draw=black] (0,0) -- (0,5) -- (10,5) --
       (10,0) -- cycle (1,5) -- (1,10) -- (11,10) -- (11,5) -- cycle
       (2,10) -- (2,15) -- (12,15) -- (12,10) -- cycle;
    \filldraw[black!20!white,draw=black] 
    (2,0) -- (2,10) .. controls +(1,1) and +(1,-1) .. (2,15) -- (10,15) .. 
    controls +(1,-1) and +(1,1) .. (10,10) -- (10,0) -- (2,0);
     \draw[dashed] (0,0) -- (0,5) -- (10,5) --
    (10,0) -- cycle (1,5) -- (1,10) -- (11,10) -- (11,5) -- cycle
    (2,10) -- (2,15) -- (12,15) -- (12,10) -- cycle;
    \end{scope}

    \begin{scope}[xshift=45cm]
      \filldraw[black!10!white,draw=black] (0,0) -- (0,5) -- (10,5) --
       (10,0) -- cycle (1,5) -- (1,10) -- (11,10) -- (11,5) -- cycle
       (2,10) -- (2,15) -- (12,15) -- (12,10) -- cycle;
    \filldraw[black!20!white,draw=black] 
    (2,0) --  (2,15) -- (10,15)-- (10,0) -- (2,0);
     \draw[dashed] (0,0) -- (0,5) -- (10,5) --
    (10,0) -- cycle (1,5) -- (1,10) -- (11,10) -- (11,5) -- cycle
    (2,10) -- (2,15) -- (12,15) -- (12,10) -- cycle;
    \end{scope}

  \end{tikzpicture}
  \caption{Zipping $\band\mu$ onto $\band\eta$, as seen from the
    universal cover. Here $\tube\mu = (\tube\eta)^3$.}  
  \label{fig:merge-zip}
\end{figure}

The next step is to zip the bands $\band{\mu}$ and $\band{\lambda}$
onto $\band\eta$. In our example, since $\tube\mu = (\tube\eta)^3$,
$\band\mu$ should ``wrap'' three times around $\band\eta$. To
accomplish this we \emph{horizontally} subdivide $\band{\mu}$ into 3
bands and successively perform three zipping operation (Definitions
\ref{move:zip} and \ref{move:squish}). This is best visualized in the
universal cover, see Figure \ref{fig:merge-zip}.
Thus, we have produced a new band complex $\mbc'$ efficiently carrying
a track $t'$ and by Proposition \ref{prop-preservation} we have
preserved the fundamental group and the dual Bass-Serre tree. Although we
considered a specific example, this discussion is sufficiently general
to make the following claim.

\begin{prop}[The periodic merger]\label{prop:periodic-merger}
  Let $\mbc$ be a band complex with an overlapping pair $\basepair\mu$
  and an overlapping pair $\basepair\lambda \subset \inter\mu$, such
  that $\basepair\mu$ and $\basepair\lambda$ are entangled. Suppose
  that $t\in \kappaTracks\mbc$. If, after widening $\band\lambda$, as
  in Lemma \ref{lem:widen} we have $|\lambda|_t\geq\tr\mu$ then
  \begin{itemize}
  \item there is a continuous
    map $m:\mbc \to \mbc'$ where $\mbc'$ with
    $m(t) = t'\subset \mbc'$ where $t'$ is a track efficiently
    carried by $\mbc'$. 
  \item The induced map $m_\sharp:\fungrp\mbc \to \fungrp{\mbc'}$ is
    an isomorphism and there is an $m_\sharp$-equivariant isomorphism
    of dual Bass-Serre trees $\dualtree t \mbc \to \dualtree{t'}{\mbc'}$.
  \item The resulting band complex $\mbc'$ has two fewer bases.
  \end{itemize}
\end{prop}
\begin{proof}
  $\mbc'$ is obtained by first perhaps widening $\band\lambda$ so that
  it is coinitial and coterminal with $\inter\mu$, then attaching a
  band $\band\eta$ so that $\inter\eta = \inter\mu$, attaching a
  2-cell, and finally horizontally subdividing and zipping $\band\mu$
  and $\band\lambda$ onto $\band\eta$. The resulting composition of
  operations preserves fundamental groups and dual trees.
\end{proof}

\subsection{A modification to $\ETone{\complex C}$: adding periodic
  mergers to the elimination process}
\label{sec:modification}
We now turn our attention to band complexes, viewed as combinatorial
objects. 

\begin{defn}\label{def:merging-inadmissible}
  Let $\mbc_v$ be a band complex in $\ETone{\complex C}$, then we say
  that $\mbc_v$ is \define{merging inadmissible} if it contains an
  overlapping pair $\basepair\mu$ that is entangled with
  $\basepair\lambda$, but the algorithm of Proposition
  \ref{prop:Bulitko-trick}, certifies that $\mbc_v$ can not
  efficiently carry a $\kappa$-track.
\end{defn}

Merging inadmissibility can be verified algorithmically. Indeed, given
$\basepair \lambda \subset \inter\mu$, Lemma
\ref{lem:compute-entanglement} states that entanglement can be
computed and merging inadmissibility is certified from the output of
the algorithm of Proposition \ref{prop:Bulitko-trick}. By definition,
this only depends on the underlying band complex, and not on the track
it carries.

Whether entangled base pairs can actually be merged, and the outcome
of this operation, depend on the track $t$ carried by $\mbc$. Merging
inadmissibility, however, guarantees that no matter the track being
carried by $\mbc$, such a merging is impossible.

Now we have proved that if $\basepair\mu$, $\basepair\lambda$ are a
pair of entangled overlapping pairs with
$\basepair\lambda \subset \inter\mu$, then for any
$t \in \kappaTracks{\mbc_v}$, after widening $\band\lambda$ as in
Lemma \ref{lem:widen} we will be able to apply a periodic merger by
Lemma \ref{lem:premerge-I}. It therefore follows that we can discard
merging inadmissible band complexes since they cannot carry
$\kappa$-tracks.

This next lemma simply follows from the fact that we can enumerate the
combinatorial outcomes all such periodic mergers since we widened some
band by a controlled amount, added a band, a 2-cell with an attaching map of length $M$, and
applied $N$ zipping operation, where $M,N$ are bounded by the output
of the algorithm of Proposition \ref{prop:Bulitko-trick}.

\begin{lem}\label{lem:bc-merge}
  Let $\mbc$ be a band complex containing overlapping pairs
  $\basepair\mu, \basepair\lambda$ that are entangled and such that
  $\basepair\lambda \subset \inter\mu$. Then we can effectively
  construct a finite set of band complexes\[
  \begin{tikzpicture}
    \node (1) at (0,0) {$\mbc$};
    \node (2) at (-1,-0.75) {$\mbc_1$};
    \node (3) at (1,-0.75) {$\mbc_{m_\mbc}$};
    \node (4) at (0,-0.75) {$\cdots$};
    \draw[->] (1) -- (2);
    \draw[->] (1) -- (3);
  \end{tikzpicture}
  \]
  containing all possible outcomes $m:\mbc \to \mbc'$ of merging
  $\band\mu$ and $\band\lambda$ as described by Proposition
  \ref{prop:periodic-merger}, with the track $t$ ranging over
  $\kappaTracks{\mbc}$.
\end{lem}

We note that although some of the band complexes produced by Lemma
\ref{lem:bc-merge} may not correspond to \emph{any} periodic mergers,
the resulting band complexes will have the same fundamental group and
a lower $\tau$-complexity. In particular if any of them admit a
$\kappa$-track, then so must the original $\mbc$. There is thus no
danger of introducing ``false positives'', by giving $\mbc$
illegitimate children. We now include periodic mergers in our
elimination tree.

\begin{defn}\label{defn:ETTwo}
  The elimination tree $\ETtwo{\mbc,J}$ is constructed inductively
  similarly as $\ETone{\mbc,J}$ in Section \ref{sec:restricted-ep} except
  with a new clause that takes precedence over item
  (\ref{it:ET-otherwise}) given in Section \ref{sec:ETone}.
  \begin{itemize}
  \item[(\ref{it:ET-otherwise}m)]\label{it:merge-clause} If $\mbc_v$ is a nonterminal band
    complex in $\ETtwo{\mbc J}$ that contains overlapping pairs
    $\basepair\mu$ and $\basepair\lambda$ satisfying the premises of
    Lemma \ref{lem:bc-merge} then define as its children the
    collection of band complexes given by Lemma \ref{lem:bc-merge}.
  \end{itemize}
\end{defn}

Obviously all the $\kappa$-inadmissibility and repetition
inadmissibility criteria on paths $\ETone{\mbc,J}$ also apply to
$\ETtwo{\mbc,J}$, as does the classification in Theorem
\ref{thm:branch-types}. For the remainder of this paper we will use
$\ETtwo{\mbc,J}$ as our elimination tree.

\section{Overlapping pairs must occur and stabilize.}
\label{sec:periodicity-reduction}

For this section let $\mbc$ be a band complex, and let
$t \subset \mbc$ be an automorphically minimal $\kappa$-track
efficiently carried by $\mbc$. The corresponding Rips sequence gives a
path in $\ETtwo{\mbc,J}$. Throughout this section we will fix a
superquadratic subpath
\begin{equation}
  \label{eqn:sq-rips-path}
  p(t_1):(\mbc_1,t_1) \to \cdots \to (\mbc_P,t_P)
\end{equation}
where, in particular, the track $t_1$ is an automorphically minimal
$\kappa$-track efficiently carried by $\mbc_1$. We will further assume
that the $J$-relative $\tau$-complexity remains constant throughout
$p$. The purpose of this section is to show that if $p$ is
sufficiently long, then it must have a tail in which some base $\mu$
is repeatedly the carrier. This is called $\mu$-periodicity (Definition
\ref{defn:mu-periodic}.)

We prove this by first defining a quantity called the excess that
remains constant throughout $p$. We then consider the quadratic part
of $\mbc_1$ and use this to show that a union of bases called the
participating segments has a length bounded above by some computable
multiple of the excess. This multiple is computed using a restricted
elimination process. It will be clear from the definition that the
excess is bounded above by some constant multiple of the length of the
longest base in each $\mbc_i$ occurring in $p$. Finally we will define
something called a C-T cycle with the property that whenever it
occurs, a considerable portion of the participating segments gets cut
out. Our bound on the total length of the participating segments will prevent these cycles from
occurring too often. The critical detail is that this bound does not
depends on the actual track $t.$ It only depends on the sequence of
underlying band complexes that occur along the path $p$ and the
assumption that $t$ is an automorphically minimal $\kappa$-track
efficiently carried by $\mbc$.

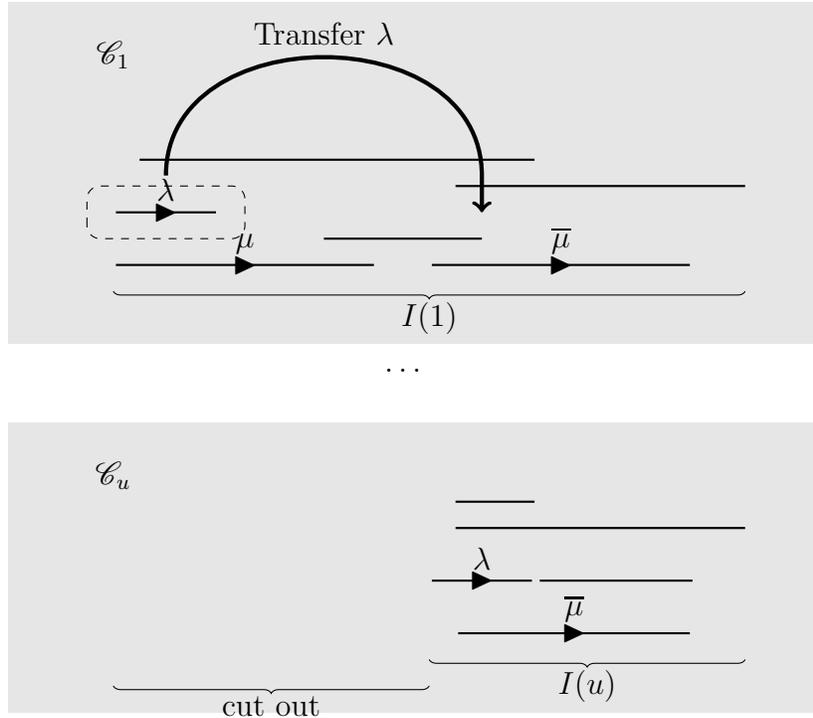
\begin{figure}[htb]
  \centering
  \begin{tikzpicture}[scale=0.7]
    \fill[black!10!white] (-1,-1.5) rectangle (14.5,5);
    \base{0}{1}{5}{\mu}
    \opdualbase{0}{7}{5}{\mu}
    \base{1}{1}{2}{\lambda}
    \draw[thick] (1.5,2) -- (9,2)
    (5,0.5) -- (8,0.5)
    (7.5,1.5) -- (13,1.5);
    \draw[dashed, rounded corners=5pt]
    (0.5,0.5) rectangle (3.5,1.5);
    \draw[ultra thick,->]
    (2,1.7) .. controls +(0,3) and +(0,3) .. node[above]{Transfer $\lambda$} (8,1.7) -- (8,1);
    \draw (1,4) node {$\mbc_1$};
    \draw[decorate,decoration={brace,amplitude=0.1cm}] (13,-0.5) --node[below]{$I(1)$} (1,-0.5);
    \draw (6.5,-2) node {$\ldots$};
    \begin{scope}[yshift=-7cm]
      \fill[black!10!white] (-1,-1.5) rectangle (14.5,4);
      \draw (1,3) node {$\mbc_u$};
      \opdualbase{0}{7.5}{4.5}{\mu}
      \base{1}{7}{2}{\lambda}
      \draw[thick] (7.5,2.5) -- (9,2.5) (7.5,2) -- (13,2) (9.1,1) --
      (12,1);
      \draw[decorate,decoration={brace,amplitude=0.1cm}] (13,-0.5) --node[below]{$I(u)$} (7,-0.5);
    
      \draw[decorate,decoration={brace,amplitude=0.1cm}] (7,-1)
      --node[below]{cut out} (1,-1);    

    \end{scope}

  \end{tikzpicture}
  \caption{In $\mbc_1$, $\mu$ is a carrier and $\lambda$ gets
    transferred. In $\mbc_u$, $\lambda$ is again leading. At least
    $|\mu|_{t_1}$ got cut out of $I(1)$.}
  \label{fig:move-mu}
\end{figure}

We will give an example of what is meant by ``cut out''. Suppose that
the carrier $\mu$ in $\mbc_1$ does not overlap with its dual, that
some base $\lambda$ gets transferred, and that eventually in $\mbc_u$
$\lambda$ is again a leading base, see Figure \ref{fig:move-mu}.
Let the interval $I(u)$ be the natural image of $I(1)$ in $\mbc_u$. On
one hand we have $|I(u)|_{t_u} < |I(1)|_{t_1}$, on the other hand we
have $|I(1)|_{t_1} - |I(u)|_{t_u} \geq |\mu|_{t_1}$, i.e we cut out at
least $|\mu|_{t_1}$ from $I(1)$.

\subsection{The excess invariant}\label{sec:excess}

The excess invariant given in Definition \ref{defn:excess} is
originally due to Makanin \cite{Makanin-1982}. It also occurs in
\cite{Bestvina-Feighn-1995,KM-IrredII}.

\begin{defn}\label{defn:participating-segment}
  Let $p$ be the path given in (\ref{eqn:sq-rips-path}), by $C(p)$ we
  denote \define{the set of bases that are carriers at some point
    along $p$} and by $T(p)$ we denote \define{the set of bases that
    are transferred at some point along $p$}.  We denote the
  \define{participating segments:}
  \[\inter p =\bigcup_{\mu\in C(p)\cup T(p)}
  \left(\mu\cup\dual\mu\right).\]
\end{defn}

Recall that we are reusing symbols, so that if $\lambda$ denotes a
base in $\mbc_1$ then it also naturally denotes a base in every
$\mbc_u$ that occurs along $p$. With this in mind it makes sense to
reuse notation to denote the corresponding subset $\inter p$ of each
$\mbc_u$ in $p$ the same way, if there is a danger of confusion we
will explicitly write $\inter p \subset \mbc_u$. Noting that entire
transformations are compositions of zipping moves and collapses
(c.f. Section \ref{sec:basic-moves}), the continuous map
$\mbc_{u}\to\mbc_{u'}$ actually induces a \emph{retraction} of
$\inter p \subset \mbc_u$ onto $\inter p \subset \mbc_{u'}$. It
follows that $\inter p \subset \mbc_{u'}$ naturally embeds into
$\inter p \subset \mbc_u$. Along the path $p(t_1)$ given in
(\ref{eqn:sq-rips-path}) we will write $|\inter p|_{t_u}$ to denote
the hitting measure of $\inter p \subset \mbc_{u}$ with respect to
$t_u$. By the embedding above we have a chain of proper
inequalities\begin{equation}\label{eqn:SQ-ineqs} |\inter p|_{t_1} >
  \cdots > |\inter p|_{t_P} > 0.
\end{equation}

\begin{defn}\label{defn:excess}
Suppose now a subset $\sigma$ is a union of bases and a track $t$ is carried by $\mbc$. We call the
following quantity \emph{excess}:\[ \excess t \sigma = \left(\sum_{\mu
    \subset \sigma} |\mu|_t\right) - 2|\sigma|_t.\] 
\end{defn}

Excess measures how far the  $\sigma$ is from being quadratic
as quantified by the hitting measure. From the definition of the
entire transformation, a straightforward counting argument gives the
following.

\begin{lem}[Excess is invariant, c.f. {\cite[(7.6.1)]{Bestvina-Feighn-1995}}]\label{lem:excess-invariant}
  Let $p(t_1)$ be the Rips process given in (\ref{eqn:sq-rips-path}),
  i.e. all moves are entire transformations and
  $J$-relative $\tau$-complexity is constant. Then we have
  equalities\[
  \excess{t_1}{\inter p} = \cdots = \excess{t_w}{\inter p}.\]  
\end{lem}

\subsection{Bounding the quadratic part}
\label{sec:frep}

We can decompose  $\inter p \subset \mbc_1$ into\[
\inter p = Q(p) \cup SQ(p) \subset \mbc_1
\] where $Q(p)$, the \define{quadratic part of $\inter p$}, is the
closure of the set of points that are contained in exactly two
bases. We define $SQ(p)$, the \define{superquadratic part of
  $\inter p$}, to be the closure of $\inter p \setminus Q(p)$ (all
points in $SQ(p)$ lie in at most three bases.) It is worth noting that
this decomposition is almost disjoint (intersection consists of
finitely many points) and that some bases may lie partially in $Q(p)$
and partially in $SQ(p)$.

If we forget the tracks $t_i$ carried by $\mbc_i$ in $p(t_1)$ then we
have an induced superquadratic path
\begin{equation}\label{eq:sq-path}
  p: \mbc_1 \to \cdots \to \mbc_P
\end{equation}
with constant $J$-relative $\tau$ complexity that lies in
$\ETone{\mbc,J}$. Note that we are not requiring $\mbc_p$ to be
terminal.

\begin{lem}\label{lem:frep-ineq}
  We can construct a computable function $\frep$, depending only on
  the band complex $\mbc$ and a path $p$ as in (\ref{eq:sq-path}),
  with positive integer values such that for any automorphically
  minimal track $t$ efficiently carried by $\mbc_1$ we have the
  inequality:
  \begin{equation}
    \label{eq:frep-ineq}
    |Q(p)|_t \leq \frep(\mbc_v,p)\cdot|SQ(p)|_t.
  \end{equation}
\end{lem}
\begin{proof}
  Subdivide each band (Definition \ref{move:subdivide}) so that each
  base either entirely lies in $Q(p)$ or intersects $Q(p)$ with empty
  interior. $Q(p)$ is now a union of maximal sections (Definition
  \ref{defn:maximal-section}). Let $J_Q$ be the complement of $Q(p)$.
  Consider the $J_Q$-restricted admissible elimination tree
  $\ADone{\mbc_1,J_Q}$. Since $\inter p$ is a union of bases that is
  closed under taking duals, the leaves of $\ADone{\mbc_1,J_Q}$ will be
  band complexes in which every base in $Q(p)$ is eventually moved
  onto $SQ(p)$. On the other hand, since $\ADone{\mbc_1,J_Q}$ is
  $J_Q$-relatively quadratic, by Propositions \ref{lem:quadratic-rep}
  and \ref{cor:must-be-superquadratic} the subtree $\ADone{\mbc_1,J_Q}$
  is finite and therefore algorithmically constructible.

  Going backwards in $\ADone{\mbc_1,J_Q}$ from every admissible leaf
  to the root $\mbc_1$ (recall Section
  \ref{sec:set-of-tracks}) and repeatedly applying the upper bound of
  Lemma \ref{lem:going-backwards} gives us a finite (algorithmically
  constructible) set of upper bounds for lengths of the bases in
  $Q(P)$ in $\mbc_1$ in terms of $|SQ(p)|_t$. The maximum over this
  set can be used to compute an upper bound of $|Q(p)|_t$ in terms of
  $|SQ(p)|_t$ for any automorphically minimal track $t$ efficiently
  carried by $\mbc_1$.
\end{proof}

Whenever a base $\lambda_u$ is the carrier in
$(\mbc_u,t_u) \to (\mbc_{u'},t_{u'})$ in (\ref{eqn:sq-rips-path}),
after transferring other bases, it gets shortened. The bound given by
Lemma \ref{lem:frep-ineq} and the invariance of excess tell us that,
although bases get shorter, their lengths remain bounded below
throughout $p$. This will force them to overlap.

\begin{lem}\label{lem:lambdamax-bound}
  Let $p(t_1)$ be as in (\ref{eqn:sq-rips-path}) and let
  $\lambda_u \in C(p) \cup T(p)$ be the base in $\mbc_u$ such that
  $|\lambda_u|_{t_u}$ is maximal. The initial length
  $|\inter p|_{t_1}$ is always bounded by
  \[ |\inter p|_{t_1} \leq
  \left(\frep(\mbc_1,p)+1\right)N^2|\lambda_u|_{t_u}
  \] where $N$ is the number of bases in $\mbc_u$.
\end{lem}
\begin{proof}
  In $(\mbc_u,t_u)$, every $\mu \in C(p)\cup T(p)$ has length at most
  $|\lambda_u|_{t_u}$. It therefore follows that $|\inter p|_{t_u}
  \leq N\cdot|\lambda_u|_{t_u}$. Now since there are no more than
  $N$ bases, each point in $\inter p$ is contained by at most $N$
  bases. The strict upper bound for the invariant excess\[
  \excess{t_1}{\inter p} = \excess{t_u}{\inter p} < N^2\cdot|\lambda_u|_{t_u}
  \] as well as the bound $|SQ(p)|_{t_1} \leq \excess{t_1}{\inter p}$ combine with Lemma \ref{lem:frep-ineq} to
  give the desired
  inequality since $|\inter p|_{t_1} = |Q(p)|_{t_1}+|SQ(p)|_{t_1}$.
\end{proof}

\subsection{Orientation reversing overlaps: two lemmas}
The proof in the next section requires two additional lemmas.

\begin{lem}\label{lem:OR-small-overlap}
  Suppose $\fungrp\mbc$ has no elements of order 2, and
  $t \in \kappaTracks\mbc$. Then for any dual pair $\basepair\mu$ such
  that $\mu \cap \dual \mu \neq \emptyset$ that is orientation
  reversing we have 
  \[ |\mu \cap \dual{\mu}|_t < \kappa+1.
  \]
\end{lem}
\begin{proof}
  Since $\mu \cap \dual\mu \neq \emptyset$ and $\basepair\mu$ is
  orientation reversing, the image of $\band\mu$ in $\mbc$ contains a
  Möbius band. Let the simple closed curve $\gamma$ be homotopic to
  the core of this band. We may view $\gamma$ an element of
  $\fungrp\mbc$. If $\mbc$ efficiently carries a track $t$ then
  $\gamma$ must invert some segment in $\dualtree t \mbc$; thus
  $\gamma \neq 1$.

  By assumption $\gamma^2\neq 1$ , but an analysis of the action of
  $\gamma^2$ on $\tilde\mbc$ and $\dualtree t \mbc$ similar to the one
  shown in Figure \ref{fig:tubular-action} shows that $\gamma^2$ fixes
  an arc of length $|\mu \cap \dual{\mu}|_t$ in $\dualtree t
  \mbc$. $\kappa$-acylindricity therefore ensures the required bound.
\end{proof}

\begin{lem}\label{lem:short-strips-II}
  Let $(\mu,\ol{\mu})$ be an overlapping pair and suppose that
  $(\lambda,\ol{\lambda})$ is orientation reversing and $\lambda \cup
  \ol{\lambda} \subset \inter{\mu}$. If $\fungrp{\mbc}$ has no 2-torsion
  then for any $t \in \kappaTracks\mbc$ \[|\lambda|_t < \tr{\mu} + \kappa + 1. \]
\end{lem}
\begin{proof}
  Suppose towards a contradiction that $|\lambda|_t \geq
  \tr{\mu} + \kappa +1$. Assume that $\lambda$ is to the left of
  $\ol{\lambda}$ and let $q$ be the leftmost point of $\lambda$. By
  successively transferring $\ol{\lambda}$ through $\band{\mu}$ we can
  arrange so that the leftmost point of $\ol{\lambda}$ is moved to a
  distance of less than $\tr{\mu}$ to the right of $q$ it therefore
  follows that $|\ol{\lambda} \cap \lambda|_t \geq \kappa +1$, which
  contradicts Lemma \ref{lem:OR-small-overlap}.
\end{proof}

\subsection{Cutting off too much}\label{sec:cutting-off-too-much}
We now introduce C-T cycles which are guaranteed to ``shorten''
$\inter p \subset \mbc_v$ by some fixed amount each time they occur. The idea of a C-T
cycle is inspired from Case 3 in the proof of \cite[Theorem
12]{TAH-2006}.

\begin{defn}\label{defn:CT-cycle}
  Let $p$ be as in (\ref{eqn:sq-rips-path}). A \define{C-T cycle} is a
  sub-path $\mbc_j \to \cdots \to \mbc_k$ such that for each
  $\lambda \in C(p)$ there is some $j\leq l \leq k$ such that $\lambda$ is the
  carrier in some $\mbc_l$ and one of the following occurs:
  \begin{itemize}
  \item If $\lambda$ either doesn't overlap with its dual or is
    orientation reversing, then for some $l<l'\leq k$, some base
    $\delta$ that was transferred by $\lambda$ in $\mbc_{l}
    \rightarrow \mbc_{l+1}$ is a leading base again in $\mbc_{l'}$.
  \item If $(\lambda,\ol{\lambda})$ form an overlapping pair, then
    there are $l<l'<l''\leq k$ such that $\lambda$ gets carried in
    $\mbc_{l'}$ (i.e. it ceases to be a carrier) and is a leading base
    again in $\mbc_{l''}$. Also there is some some base $\delta$ that
    was transferred in $\mbc_{l} \rightarrow \mbc_{l+1}$ is again a
    leading base in $\mbc_{l'''}$ for some $l<l'''\leq k$.
  \end{itemize}
\end{defn}

\begin{lem}\label{lem:ct-cycle-cut}
  Let $p(t_1)$ be a path as in (\ref{eqn:sq-rips-path}) and let
  $\mbc_j \to \cdots \to \mbc_k$ be a C-T cycle in $p(t_1)$, then
  either every base in $C(p)\cup T(p)$ has length at most
  $\max\{2\kappa,1\}$, or
  \begin{equation}\label{eqn:cut-bound}
    |\inter p|_{t_j} - |\inter p|_{t_k} \geq \frac{|\lambda|_{t_j}}{2(2+\kappa)}
  \end{equation} 
where $\lambda$ is the longest base in $\mbc_j$.
\end{lem}

\begin{proof}
  First note that until $\lambda$ is the carrier in some $\mbc_l$
  where $j \leq l \leq k$ we have
  \[|\lambda|_{t_j} = |\lambda|_{t_j+1} = \cdots = |\lambda|_{t_l}.\]
  By definition of a C-T cycle, there is such a $j \leq l \leq k$
  where $\lambda$ is the carrier in $\mbc_{l}$.  Assume that
  $|\lambda|_{t_j} > \max\{2\kappa,1\}$.  We will now show
  (\ref{eqn:cut-bound}), the proof divides into cases; decimals denote
  subcases.
  \\~\\
  \noindent {\bf Case 1:} \textit{$\lambda$ does not overlap with its
    dual.} In this case let $\delta$ be some base that is carried by
  $\lambda$ it is moved by at least $|\lambda|_{t_j}$ to the right. By
  definition of a C-T cycle there some $l< l' \leq k$ where $\delta$
  is a leading base in $\mbc_{l'}$ again which means that
  $|\inter p|_{t_l} - |\inter p|_{t_l'} \geq |\lambda|_{t_j}$ so the
  result holds. (See Figure \ref{fig:move-mu}.)
  \\~\\
  \noindent {\bf Case 2:} \textit{$\lambda$ has nontrivial
    intersection
    with $\dual{\lambda}$.}
  \\~\\
  \noindent {\bf Case 2.1:} \textit{$\basepair \lambda$ is orientation
    reversing.} Let $\delta$ be some base that carried by
  $\lambda$. Then $\delta$ is moved
  $2|\lambda|_{t_l}-|\lambda \cap \ol{\lambda}|_{t_l} -
  |\delta|_{t_l}$
  to the right. On one hand $|\delta|_{t_l} \leq |\lambda|_{t_l}$, on
  the other hand by Lemma \ref{lem:OR-small-overlap}
  $|\lambda \cap \ol{\lambda}|_{t_l} \leq \kappa$ we therefore
  conclude that $\delta$ is moved at least $|\lambda|_{t_l}-\kappa$ to
  the right. When $\delta$ is again a leading base in
  $\mbc_{l'}, l < l' \leq k$ then
  \[|\inter p|_{t_l} - |\inter p|_{t_{l'}} \geq |\lambda|_{t_l}-\kappa-1
  \geq \frac{|\lambda|_{t_l}}{2}\] so (\ref{eqn:cut-bound}) holds.
  \\~\\
  \noindent {\bf Case 2.2:} \textit{$\basepair \lambda$ form an
    overlapping pair.} There is some maximal $l'$ such that in $\mbc_{l}
  \rightarrow \ldots \rightarrow \mbc_{l'-1}$, $\lambda$ is always the
  carrier base. Note moreover that $\tr{\lambda}$ remains constant.
  \\~\\
  \noindent {\bf Case 2.2.1:} \textit{$|\lambda|_{t_l} \leq
    (2\kappa+2)\tr{\lambda}$.} Let $\delta$ be some base that is
  carried by $\lambda$. Then it is moved to the right by
  $\tr{\lambda}$ and by hypothesis there is some $l < l''' \leq k$
  such that $\delta$ is leading again, hence \[|\inter p|_{t_l} -
  |\inter p|_{t_{l'}} \geq \frac{|\lambda|_{t_l}}{2\kappa+2}\] and
  (\ref{eqn:cut-bound}) holds.
  \\~\\
  \noindent {\bf Case 2.2.2:} \textit{$|\lambda|_{t_l} >
    (2 \kappa+2)\tr{\lambda}$.}
  \\~\\
  \noindent {\bf Case 2.2.2.1:} \textit{$\displaystyle |\inter p|_{t_l} - |\inter
    p|_{t_{l'}} \geq \frac{|\lambda|_{t_l}}{2(2\kappa+2)}$.}
  (\ref{eqn:cut-bound})
  immediately holds.
  \\~\\
  \noindent {\bf Case 2.2.2.2:}
  \textit{$\displaystyle |\inter p|_{t_l} - |\inter p|_{t_{l'}} <
    \frac{|\lambda|_{t_l}}{2(2\kappa+2)}$.}
  In particular, in $\mbc_{l'}$
  \begin{equation}\label{eqn:2222}|\lambda|_{t_{l'}} > |\lambda|_{t_l} \left( 1 -
      \frac{1}{2(2\kappa+2)} \right) =
    |\lambda|_{t_l}\left(\frac{4\kappa+3}{2(2\kappa+2)}\right)>
    (2\kappa+3/2)\tr{\lambda}
  \end{equation} by the Case 2.2.2 assumption on $|\lambda|_{t_l}$.
  \\~\\
  \noindent {\bf Case 2.2.2.2.1:} \textit{The carrier $\eta$ in
    $\mbc_{l'}$ does not overlap with its dual.} In this case
  $\lambda$ is moved to the right by at least $|\lambda|_{t_{l'}}$ and
  we note that for all $\kappa \geq 0$ we have:\[ 1-
  \frac{1}{2(2\kappa + 2)} > 1/2 > \frac{1}{2(2\kappa+2)};
  \]
  thus when $\lambda$ is a leading base again in $\mbc_{l''}$, by
  (\ref{eqn:2222}) we have
  \[|\inter p|_{t_{l'}}- |\inter p|_{t_{l''}} >
  \frac{|\lambda|_{t_l}}{2}\] so (\ref{eqn:cut-bound}) holds.
  \\~\\
  \noindent {\bf Case 2.2.2.2.2} \textit{The carrier $\eta$ in
    $\mbc_{l'}$ has nontrivial intersection with $\dual{\eta}$ and is
    orientation reversing.} Note that
  \[|\eta|_{t_{l'}}> |\lambda|_{t_{l'}} \geq |\lambda|_{t_l} \left( 1
    - \frac{1}{2(2\kappa+2)} \right).\]
  As in Case 2.1 we deduce that $\eta$ carries $\lambda$ more than
  $|\lambda|_{t_{l'}}-\kappa$ to the right. So that in $\mbc_{l''}$
  when $\lambda$ is leading again at least $|\lambda|_{t_{l'}}-\kappa$
  was cut from $\inter p$. Suppose towards a contradiction that
  \[|\lambda|_{t_{l'}}-\kappa <
  \frac{|\lambda|_{t_l}}{2(2\kappa+2)}.\]
  Then since we are in Case 2.2.2.2 we have
  \[ |\lambda|_{t_l} \left( 1 - \frac{1}{2(2\kappa+2)} \right)-
  \kappa< |\lambda|_{t_{l'}}-\kappa \]
  Combining these gives
  \begin{eqnarray*}
    & \displaystyle |\lambda|_{t_l} \left( 1 - \frac{1}{2(2\kappa+2)} \right)- \kappa 
    & <  \frac{|\lambda|_{t_l}}{2(2\kappa+2)}\\
    \Rightarrow 
    &  \displaystyle |\lambda|_{t_l} \left( 1 - \frac{2}{2(2\kappa+2)}
      \right)  
    & < \kappa\\
    \Rightarrow 
    &
      \displaystyle |\lambda|_{t_l} \left( \frac{4\kappa+2}{4\kappa+4}
      \right) 
    & < \kappa\\
    \Rightarrow
    & \displaystyle |\lambda|_{t_l} < \left(1+\frac{2}{(4\kappa+4)}\right)\kappa 
    & < \frac 3 2 \kappa\\
  \end{eqnarray*}
  which contradicts our assumption that $|\lambda|_{t_l}>2\kappa$
  \\~\\
  \noindent {\bf Case 2.2.2.2.3:} \textit{$\basepair \eta$ form an
    overlapping pair where $\eta$ in $\mbc_{l'}$ is the carrier.} We
  may assume that, after repeatedly getting shortened in
  $\mbc_{l} \rightarrow \ldots \rightarrow \mbc_{l'-1}$, $\lambda$ is
  still long enough in $\mbc_{l'}$ for $\basepair \lambda$ to
  an overlapping pair. Indeed suppose this was not the case. By not
  overlapping we have $|\lambda|_{t_{l'}} < \tr{\lambda}$, and
  substituting into the the Case 2.2.2 assumption yields
  \[|\lambda|_{t_{l'}} < \frac{|\lambda|_{t_l}}{(2\kappa+2)}.\]
  This contradicts (\ref{eqn:2222}) since
  \[ 1-\frac{1}{2(2\kappa+2)} > \frac{1}{2\kappa+2}. \]
  \\~\\
  \noindent {\bf Case 2.2.2.2.3.1:} \textit{(Recall Definition
    \ref{defn:overlapping-interval}) In $\mbc_{l'}$,
    $\inter{\eta} \subset \inter{\lambda}$.} This means, since we are
  in Case 2.2.2.2, that
  \[|\eta|_{t_{l'}} > (2\kappa+3/2)\tr{\lambda} > \tr{\lambda} +
  \kappa + 1\]
  so by Lemma \ref{lem:small-shifts-commute} $\tube{\lambda}$ and
  $\relloop{\lambda}{\eta}$ commute, so the overlapping pairs
  $\basepair\mu$ and $\basepair\lambda$ are entangled. By Lemma
  \ref{lem:both-overlap-wide} we can apply the periodic merger of
  Proposition \ref{prop:periodic-merger} to merge the bands
  $\band\eta$ and $\band\mu$. Since we are working in the elimination
  tree $\ETtwo{\complex{C}}$ by (4m) of Definition
  \ref{defn:ETTwo} \emph{we must merge the bands $\band{\eta}$ and
    $\band{\lambda}$ which decreases the complexity.} This contradicts
  the assumption that the  $\tau$-complexity remains constant
  throughout $p$.
  \\~\\
  \noindent {\bf Case 2.2.2.2.3.2:} \textit{In $\mbc_{l'}$,
    $\inter{\eta} \supset \inter{\lambda}$.} We finally distinguish
  two sub-cases:
  \\~\\
  \noindent {\bf Case
    2.2.2.2.3.2.1:}\textit{$(\kappa+2)\tr{\eta}\leq
    |\lambda|_{t_{l'}}.$}
  Again as in Case 2.2.2.2.3.1 we can perform a periodic merger.
  \\~\\
  \noindent {\bf Case 2.2.2.2.3.2.2:}\textit{$(\kappa+2)\tr{\eta} >
    |\lambda|_{t_{l'}}.$} In this case $\lambda$ gets moved by
  $\tr{\eta}$ to the right. $\lambda$ is again a leading base in
  $\mbc_{l''}$ then we will have cut at least $\tr{\eta}$ from
  $\inter p$ in passing from $\mbc_l'$ to $\mbc_{l''}$. On one hand
  since this is a sub-case of Case 2.2.2.2 we have
  \[\tr{\eta} > \frac{|\lambda|_{t_{l'}}}{\kappa+2} > \frac{
    |\lambda|_{t_l} }{ \kappa + 2 } \left( 1 - \frac{1}{2(2\kappa+2)}
  \right),\]
  where the last inequality is from (\ref{eqn:2222}), and we can
  estimate
  \begin{eqnarray*}
    &\displaystyle \frac{2(2\kappa+2)}{\kappa+2}\frac{1}{2} & \geq 1\\
  \Rightarrow &\displaystyle \frac{2(2\kappa+2)}{\kappa+2}\left( 1 - \frac{1}{2
    (2\kappa + 2)}  \right) & >  1\\
  \Leftrightarrow & \displaystyle\frac{1}{\kappa + 2}\left( 1 - \frac{1}{2
    (2\kappa + 2)}\right) & > \frac{1}{2(2\kappa+2)}\\
  \Leftrightarrow &\displaystyle \frac{|\lambda|_{t_l}}{\kappa + 2}\left( 1 -
  \frac{1}{2
    (2\kappa + 2)}\right) & > \frac{|\lambda|_{t_l}}{2(2\kappa+2)}\\
  \Rightarrow & \tr{\eta} & > \frac{|\lambda|_{t_l}}{2(2\kappa+2)}.
  \end{eqnarray*} (\ref{eqn:cut-bound}) therefore holds and, furthermore,
  all possibilities have been exhausted. 
\end{proof}

We can now combine lemmas \ref{lem:lambdamax-bound} and
\ref{lem:ct-cycle-cut}. 

\begin{cor}\label{cor:not-too-many-cycles}
  Let $p$ be a path as in (\ref{eqn:sq-rips-path}), i.e. a path
  induced by an automorphically minimal $\kappa$-track efficiently
  carried by $\mbc_1$. Then at most
  \[N^2(\frep(\mbc_1,p)+1)( 2 \kappa+2 )\] disjoint C-T cycles can
  occur.
\end{cor}
\begin{proof}
  By Lemmas \ref{lem:ct-cycle-cut} and \ref{lem:lambdamax-bound},
  whenever a C-T cycle $\mbc_j \to \cdots \to \mbc_k$ occurs we can
  bound from below the difference
  \[|\inter p|_{t_j} - |\inter p|_{t_k} \geq \frac{|\inter
    p|_{t_1}}{N^2(\frep(\mbc_1,p)+1)( 2 \kappa+2 )}.\] Since $|\inter
  p|_{t_1} > |\inter p|_{t_2} > \cdots>0$ the desired bound on the
  number of C-T cycles follows.

\end{proof}

\subsection{C-T-inadmissibility and a reduction to $\mu$-periodicity}
\label{sec:C-T-inadmissble}
Given a path $p$ in a (restricted) elimination tree $\ETtwo{\mbc,J}$
we can define the sets $C(p)$ and $T(p)$ (Definition
\ref{defn:participating-segment}) and therefore corresponding C-T
cycles (Definition \ref{defn:CT-cycle}).

\begin{defn}\label{defn:CT-inadmissible}
  If a path $p:\mbc_v \to \cdots \to \mbc_u$ in $\ETtwo{\mbc,J}$
  contains more disjoint C-T cycles than the computable bound given by Corollary
  \ref{cor:not-too-many-cycles}, then it is called
  \define{C-T-inadmissible.}
\end{defn}

\begin{defn}[Admissible]\label{defn:admissible}
  A subtree of $\ETtwo{\mbc,J}$ is said to be
  \define{admissible} if it doesn't contain any leaves that are
  inadmissible (recall (\ref{it:inadmissible}) in Section
  \ref{sec:EP}), $\kappa$-inadmissible paths (Definition
  \ref{defn:kappa-inadmissible}), repetition inadmissible paths
  (Definition \ref{defn:rep-inadmissible}), or C-T inadmissible
  paths. We denote by $\ADtwo{\mbc,J}\subset \ETtwo{\mbc,J}$ the
  maximal admissible subtree, and call it the \define{admissible elimination tree}.
\end{defn}

\begin{prop}\label{prop:admissible-contains-II}
  Any path $p:\mbc\to\cdots$ in $\ETtwo{\mbc,J}$ induced by an
  automorphically minimal $\kappa$-track efficiently carried by $\mbc$
  must also lie in $\ADtwo{\mbc,J}$.
\end{prop}
\begin{proof}
  This follows immediately from Proposition
  \ref{prop:admissible-contains} and Corollary
  \ref{cor:not-too-many-cycles}.
\end{proof}

\begin{defn}\label{defn:mu-periodic}
  Let $\mu$ be the carrier base in $\mbc$. A path\[
  p:\mbc_u \to \cdots
  \]
  in $\ETtwo{\mbc,J}$ is called \define{$\mu$-periodic} if, throughout
  $p$, $\mu$ is the carrier base and $\basepair\mu$ forms an
  overlapping pair.
\end{defn}

\begin{prop}\label{prop:periodicity-reduction}
  Any infinite path of $\ADtwo{\mbc,J}$ has a tail
  \[\mbc_u \to \cdots\]
  that is $\lambda$-periodic for some base $\lambda$.
\end{prop}
\begin{proof}
  Let $p$ be some infinite path of $\ADtwo{\mbc, J}$. We can
  form the sets $C(p)$ and $T(p)$, since every tail $p'$ of $p$ gives
  $C(p') \subset C(p)$ and $T(p')\subset T(p)$, passing to a tail of
  $p$ we may assume that each base in $C(p)$ is carrier infinitely
  often and every base in $T(p)$ is carried infinitely often.

  We may further assume that the (relative) $\tau$-complexity remains
  constant, that no annulus subdivisions or Möbius moves occur, and
  that $p$ is superquadratic. If $C(p)$ consists of more than one
  element then infinitely many C-T cycles (Definition
  \ref{defn:CT-cycle}) occur so $p$ is not contained in an admissible
  subtree.

  It therefore follows that some base $\lambda$ must repeatedly be the
  carrier throughout $p$. Now if $\lambda$ doesn't overlap with $\dual
  \lambda$ after $N$ (the total number of bases) entire transformations
  $\lambda$ can no longer be a maximal leading base. It follows that
  $\basepair \lambda$ form an overlapping pair.
\end{proof}

Therefore, if we can find a computable bound on the number of times in
a row the same base $\lambda$ can be a carrier base in some
$\lambda$-periodic path $\mbc_v\to\cdots$ induced by a minimal
$\kappa$-track efficiently carried by $\mbc_v$, then we will be able to
effectively constructed a finite subtree of $\ADtwo{\complex C}$,
whose leaves give a set of tracks containing all automorphically
minimal $\kappa$-tracks efficiently carried by $\mbc$.

\section{Bounding the periodicity of overlapping pairs}
\label{sec:bound-periodicity}
\begin{defn}\label{defn:po-carrier}
  We say that a maximal leading base $\lambda$ is a \define{principal
    overlapping carrier} if $\basepair\lambda$ form an
  overlapping pair and there are no other overlapping pairs
  $\basepair\mu$ that can be merged with $\basepair\lambda$.
\end{defn}

\begin{figure}[htb]
  \centering
  \begin{tikzpicture}[scale=0.7]
    \fill[black!10!white] (-1,-1) rectangle (8,3);
    \base{0}{0}{5}{\lambda}
    \opdualbase{1}{2}{5}{\lambda}
    \base{2}{0}{3}{\delta}
    \base{1}{1}{1}{\mu}
    \draw[ultra thick,->] (3.5,-1.1) -- (3.5,-1.9);
    \begin{scope}[yshift=-5cm]
      \fill[black!10!white] (-1,-1) rectangle (8,3);
      \base{0}{1}{4}{\lambda}
      \opdualbase{1}{3}{4}{\lambda}
      \base{2}{2}{3}{\delta}
      \base{1}{1}{1}{\mu}
    \end{scope}

    \draw[gray,dashed](0,3) -- (0,-6) (2,3) -- (2,-6);
    \draw[decorate,decoration={brace,amplitude=0.1cm}] (2,-6.2) --node[below]{$\tr{\lambda}$} (0,-6.2);
  \end{tikzpicture}
  \caption{An entire transformation where the leading base $\lambda$
    forms a principal overlapping carrier. The transfer base $\delta$ is
    moved to the right by $\tr \lambda$. If $\delta$ is eventually a
    leading base again, then $|\lambda|$ will have decreased by
    $\tr\lambda$.}
  \label{fig:periodicity-motivation}
\end{figure}
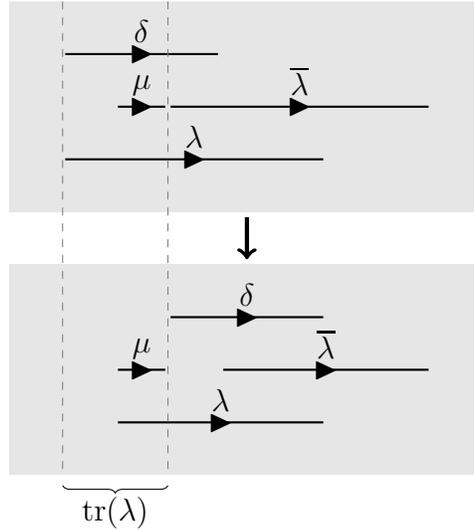

By Proposition \ref{prop:periodicity-reduction}, in any sufficiently
long admissible branch in $\ETtwo{\complex C}$, the situation depicted
in in Figure \ref{fig:periodicity-motivation} will occur.
At each such entire transformation although $|\lambda|$ decreases,
$\tr\lambda$ is invariant. Taking inspiration from words, where a
large initial segment of a word overlaps with a terminal segment,
e.g. $abcabcabcab$, we have the following:

\begin{defn}
  Let $\mbc$ be a band complex efficiently carrying a track $t$ and
  suppose that $\basepair\lambda$ form an overlapping pair. We define
  the \define{periodicity of $\lambda$} to be the positive integer
  \[ \period{t}{\lambda} = \Big\lfloor
  \frac{|\lambda|_t}{\tr{\lambda}}\Big\rfloor.\]
\end{defn}

\begin{lem}\label{lem:period-bound}
  Let $\lambda$ be a principal overlapping carrier in an band complex
  $\mbc$ efficiently carrying a track $t$ with $N$ bases. Then in the
  Rips process induced by $t$ starting at $(\mbc,t)$, $\lambda$ can be the
  maximal leading base at most $N\period{t}{\lambda}$ times in a row.
\end{lem}
\begin{proof}
  Each time a base $\delta$ is carried by $\lambda$ it is moved to the
  right by $\tr{\lambda}$. The next time that base is carried the base
  $\lambda$ will be shortened by $\tr{\lambda}$ so the advertised
  bound holds.
\end{proof}

It therefore follows that the periodicity of $\lambda$ bounds the
number of consecutive times $\lambda$ can be a carrier in a Rips
process. Equivalently, this is the maximal length of a
$\lambda$-periodic path. Given an
band complex $\mbc$ that has a principal overlapping carrier
$\lambda$, we will compute an upper bound for $\period t \lambda$
that holds for every automorphically minimal $\kappa$-track $t$. If $\period t \lambda$ exceeds this bound,
we prove the existence of another automorphically equivalent track $t'$ efficiently
carried by $\mbc$ such that $\size{t'} < \size{t}$. It is worth noting
that, since we are requiring \emph{automorphic} equivalence, simply
applying Bulitko's lemma is not good enough, as the latter only
bounds minimal periodicity in possibly non-injective homomorphic
images.

To compute this bound we will construct \emph{auxiliary elimination
  trees}. This will require all the machinery developed up until now
as well as a few new ideas. The auxiliary tree will in fact be a
rooted tree of trees. Off the leaves of this tree we will be able to
read an upper bound for the periodicity.

Proposition \ref{prop:periodicity-reduction} combined with this
periodicity bound for every occurring principal overlapping carrier in
$\ADtwo{\complex C}$ will finally enable us to construct a finite
subtree guaranteed to give us all the tracks we need for Theorem
\ref{thm:main}.

\subsection{Periodic block form, and the tree
  $\tpbf{\mbc}{J}$}\label{sec:pbf}

\begin{defn}\label{defn:periodic-block}
  Suppose that for some overlapping pair $(\lambda,\ol{\lambda})$,
  $\inter{\lambda}$ is a maximal section (recall Definition
  \ref{defn:maximal-section}). Then we call $\inter{\lambda}$ a
  \define{periodic block.}  A band complex such that every base lies
  in some periodic block is said to be in \define{periodic block
    form}.
\end{defn}

\begin{figure}[htb]
  \centering
  \begin{tikzpicture}[yscale=0.6,xscale=0.45]
    \draw[ultra thick,->] (9,1.5) --node[above]{(\ref{it:add-delta})} (13,1.5);
    \draw[ultra thick,->] (13,0) -- node[above]{(\ref{it:move-through-delta})} (9,-2);
    \draw[ultra thick,->] (7,-9) -- (7,-9.5) --node[below]{(\ref{it:truncate-delta})} (15,-9.5);

    \fill[black!10!white] (-0.5,-1) rectangle (8.5,4);
    \base{1}{0}{4}{\lambda}
    \opdualbase{0}{1}{4}{\lambda}
    \base{2}{0}{3}{\mu}
    \base{3}{1}{2}{\rho}
    \base{2.5}{2}{5}{\nu}
    \draw[thick] (3,1.5) -- (8,1.5);
    \draw[thick] (4,0.5) -- (7,0.5);
    \draw[thick] (5,1) -- (8,1);
    \begin{scope}[xshift=14cm]
      \fill[black!10!white] (-0.5,-1.5) rectangle (14.5,4);
      \base{1}{0}{4}{\lambda}
      \opdualbase{0}{1}{4}{\lambda}
      \base{2}{0}{3}{\mu}
      \base{3}{1}{2}{\rho}
      \base{2.5}{2}{5}{\nu}
      \draw[thick] (3,1.5) -- (8,1.5);
      \draw[thick] (4,0.5) -- (7,0.5);
      \draw[thick] (5,1) -- (8,1);
      \base{-1}{0}{5}{\delta}
      \opdualbase{-1}{9}{5}{\delta}
      \draw[gray,dashed](0,-1.6) -- (0,4.5) (5,-1.6) -- (5,4.5);
    \end{scope}
    \begin{scope}[yshift=-7cm]
      \fill[black!10!white] (-0.5,-1.5) rectangle (14.5,4);
      \base{1}{9}{4}{\lambda}
      \opdualbase{0}{10}{4}{\lambda}
      \base{2}{9}{3}{\mu}
      \base{3}{11}{2}{\rho}
      \base{2.5}{2}{5}{\nu}
      \draw[thick] (3,1.5) -- (8,1.5);
      \draw[thick] (4,0.5) -- (7,0.5);
      \draw[thick] (5,1) -- (8,1);
      \base{-1}{0}{5}{\delta}
      \opdualbase{-1}{9}{5}{\delta}
      \draw[gray,dashed](9,-1.6) -- (9,4.5) (14,-1.6) -- (14,4.5);
    \end{scope}
    \begin{scope}[yshift=-7cm,xshift=14cm]
      \fill[black!10!white] (1.5,-3) rectangle (14.5,4);
      \base{1}{9}{4}{\lambda}
      \opdualbase{0}{10}{4}{\lambda}
      \base{2}{9}{3}{\mu}
      \base{3}{11}{2}{\rho}
      \base{2.5}{2}{5}{\nu}
      \draw[thick] (3,1.5) -- (8,1.5);
      \draw[thick] (4,0.5) -- (7,0.5);
      \draw[thick] (5,1) -- (8,1);
      \base{-1}{2}{3}{\delta}
      \opdualbase{-1}{11}{3}{\delta}
      \draw[decorate,decoration={brace,amplitude=0.1cm}](14,-1.5) --node[below]{$\inter\lambda$} (9,-1.5);
    \end{scope}
  \end{tikzpicture}
  \caption{Creating a new periodic block for the principal overlapping
    carrier $\lambda$. Depicted is the sequence
   (\ref{it:add-delta}) - (\ref{it:truncate-delta}) in the construction
    of $\tpbf{\mbc}{J}$.}
  \label{fig:tpbf}
\end{figure}
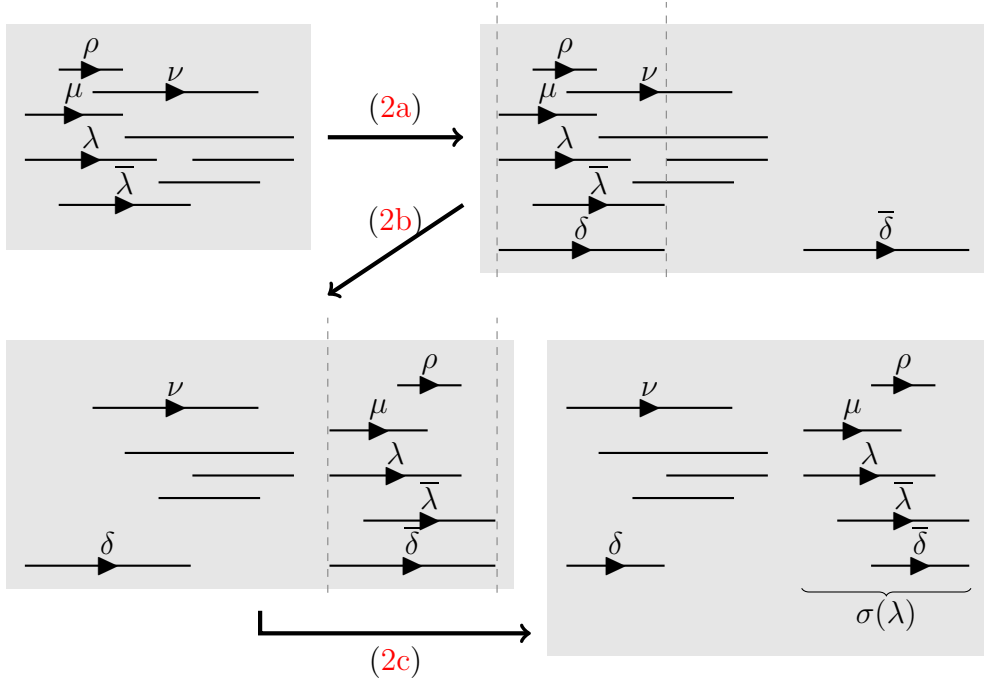

We now describe another $J$-restricted elimination process which
constructs the tree $\tpbf{\mbc}{J}$. This process brings a band
complex $\mbc$ into periodic block form and is a variation of the
construction of $\ADtwo{\mbc,J}$; only this time, whenever we
encounter a principal overlapping carrier, we add a new band
$\band\delta$, transfer $\basepair\lambda$ all the way to the right
using entire transformations, and enlarge $J$ to
$J \cup \inter{\lambda}$, see Figure \ref{fig:tpbf}. By Proposition
\ref{prop:periodicity-reduction}, we will have constructed a finite
 admissible tree in which every leaf is a band complex whose
bases either lie in $J$ or in some periodic block.
We build $\tpbf{\mbc}{J}$ as follows:
\begin{enumerate}[(1)]
\item We perform a $J$-restricted elimination process, adding only
  admissible band complexes (recall Definition
  \ref{defn:admissible}).
\item\label{it:new-block} If at some point as we grow our tree there is a $\mbc_i$ with a principal overlapping carrier $\lambda$ then instead
  as doing an entire transformation we do the following:
  \begin{enumerate}[(a)]
  \item\label{it:add-delta} We attach a band $\band{\delta}$ to
    $\mbc_i$ by identifying $\delta \bandeq \inter{\lambda}$, so that
    $\dual \delta$ doesn't meet any other bases. We  extend the
    ordering $<_i$ so that
    the maximal section corresponding to $\dual\delta$ is terminal (Definition
    \ref{defn-ordering}.)
  \item\label{it:move-through-delta} We transfer all the bases contained in $\inter{\lambda}$ through
    $\band{\delta}$ onto $\ol{\delta}$.
  \item\label{it:truncate-delta} We collapse the naked initial segment
    of $\delta$. This gives the $\mbc'_i$.
  \item\label{it:enlarge-J} We set $J' = J \cup \inter{\lambda}$ and
    we continue growing our tree at $\mbc'_i$ by returning to step 1
    but with $J'$ in place of $J$.
  \end{enumerate}
  
\item If after step (\ref{it:enlarge-J}) we have moved all the bases onto $J'$ then
  we stop. 
\end{enumerate}

By direct inspection we verify:
\begin{lem}\label{lem:stays-sq}
  If $(\mbc_i,J)$ is quadratic or superquadratic, then after
  performing steps (\ref{it:add-delta}) - (\ref{it:enlarge-J}), the
  resulting band complex is still quadratic or superquadratic (respectively).
\end{lem}

\begin{lem}\label{lem:pb-complexity}
  If we add a band $\band{\delta}$ to $\mbc_1$ as in
  (\ref{it:add-delta}), then transfer all the bases contained in
  $\inter{\lambda}$ onto $\ol{\delta}$ as in
  (\ref{it:move-through-delta}) to produce $\mbc_i'$, then $\tau(\mbc_i)
  \geq \tau(\mbc_i')$. Moreover if $J$ and $J'$ are as in
  (\ref{it:enlarge-J}) above then we have a strict
  inequality \[\tau(\mbc_i, J) > \tau(\mbc'_i,J').\]
\end{lem}
\begin{proof}
  We keep track of the $\tau$-complexity. Adding the band
  $\band{\delta}$ adds a base to a maximal section $\sigma$ with
  $b(\sigma)\geq 2$ (recall Definitions \ref{defn:maximal-section} and
  \ref{defn:tau-complexity}), which increases the $\tau$-complexity by
  1, and creates another maximal section $\sigma(\ol{\delta})$ with
  only one base, which doesn't contribute to the $\tau$-complexity. We
  then move the base $\lambda$ onto $\sigma(\ol{\delta})$ which
  decreases $b(\sigma)$ by 1. Now $b(\sigma(\ol{\delta}))=2$, so it
  still contributes 0 to the $\tau$-complexity. So far we have added 1
  and removed 1 from the $\tau$-complexity. Since $\lambda$ was
  assumed to be a leading base and $\delta \supset \lambda$, doing the
  rest of (\ref{it:move-through-delta}) and (\ref{it:truncate-delta})
  amounts to sequence of transformations, which do not increase the
  $\tau$-complexity.

  After all this $\inter{\lambda}$ is a maximal section with
  $\tau(\inter\lambda) > 0$ and $\inter{\lambda} \cap J = \emptyset$
  which implies that $\tau(\mbc_i,J) > \tau(\mbc'_i,J')$.
\end{proof}

\begin{cor}\label{cor:tpbf-finite}
  $\tpbf{\mbc}{J}$ is finite. Furthermore, if $t \subset
  \mbc$ is an efficiently carried automorphically minimal
  $\kappa$-track, then the path $\mbc \to \cdots \to \mbc_p$ to a
  band complex in periodic block form induced by $t$ is contained in
  $\tpbf{\mbc}{J}$.
\end{cor}
\begin{proof}
  We first prove the first statement. Suppose towards a contradiction
  that this was not the case, then $\tpbf{\mbc}{J}$ has an infinite
  branch $\branch{}$. Seeing as we are constructing an admissible
  elimination tree in the sense of Definition \ref{defn:admissible},
  by Proposition \ref{prop:periodicity-reduction} this infinite branch
  can be assumed to start with some $\mbc_u$ with a principal
  overlapping carrier $\lambda$. By the definition of
  $\tpbf{\mbc}{J}$, this means we must construct a new periodic block,
  item (\ref{it:new-block}), which by Lemma \ref{lem:pb-complexity}
  strictly decreases the relative $\tau$-complexity, so this event can
  only happen finitely many times, contradicting the the fact that
  $\branch{}$ is infinite.

  The second claim follows immediately from Proposition
  \ref{prop:admissible-contains-II}.
\end{proof}

\subsection{Normalized periodic block form}\label{sec:second-periodic-merger}

Once a band complex $\mbc$ is in periodic block form, it will be
possible to perform periodic block mergers, which decreases the number
of periodic blocks and the $\tau$-complexity. Furthermore it will enable
us to put a partial order on the periodic blocks, called a periodic
hierarchy. First we give another version of entanglement.

\begin{defn}\label{defn:block-overlapping-pair}
  If $\basepair\mu$ is an overlapping pair such that $\inter\mu$
  is a maximal section (Definition \ref{defn:maximal-section}) then
  $\basepair\mu$ is called a \define{block overlapping pair}.
\end{defn}

In particular, if $\mbc$ is in periodic block from, then every maximal
section is in fact a block overlapping pair.

\begin{figure}[htb]
  \centering
  \begin{tikzpicture}[scale=0.7]
    \olband{-1}{0}{5}{1}
    \filldraw[fill=black!10!white] (1,0) rectangle (3,-2);
    \begin{scope}[yshift=-2cm,yscale=-1]
      \olband{0}{0}{7}{1.5}
    \end{scope}
    \draw[very thick] (2,0) -- (2,-2);
    \fill (2,-1) circle (0.1) 
    (2,-1) node[right]{$p$};
    \draw (2.5,2) node{$\blacktriangleleft$} node[above]{$\beta$}; 
    \fill (3,0) circle (0.1); 
    \fill (2,0) circle (0.1); 
    \tubeside{2}{0}{1}
    \draw[very thick] (3,0) --node[sloped]{$\blacktriangleleft$}
    node[below]{$\gamma$} (2,0);
    \draw[very thick] (2,-1) --node[sloped]{$\blacktriangleright$} node[left]{$\alpha$} (2,0);
    \begin{scope}[yshift=-2cm,yscale=-1]
      \tubeside{2}{0}{1.5}
      \draw[very thick]  (2,0) -- (3.5,0);
    \end{scope}
  \end{tikzpicture}
  \caption{$\delta$-relative tubular elements. The relative tubular
    element $\orelloop \delta \lambda$ is the loop based at $p$ given
    by $\alpha * \beta *\gamma *\alpha^\mo$.}
  \label{fig:other-tubes}
\end{figure}
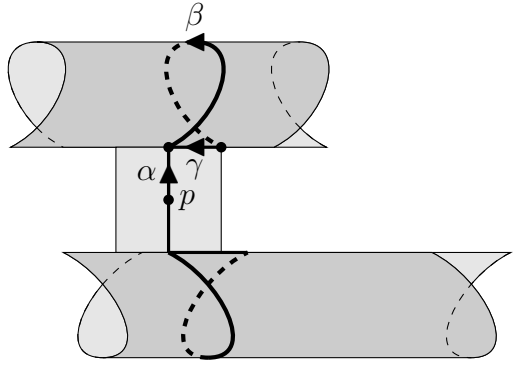

\begin{defn}\label{defn:tp2-entanglement}
  Let $\basepair \lambda$ and $\basepair \mu$ disjoint block
  overlapping pairs. Suppose there is a band $\band\delta$ connecting
  $\inter\lambda$ and $\inter\mu$, i.e. $\delta \subset \inter\lambda$
  and $\dual\delta \subset \inter\mu$. For any $p \in \band\delta$ we
  can define $\delta$-relative tubular elements
  $\orelloop{\delta}{\lambda},\orelloop \delta \mu$ as in Figure
  \ref{fig:other-tubes}. \define{$\basepair \lambda$ and
    $\basepair \mu$ are entangled by $\band\delta$} if
  \[ [\orelloop{\delta}{\lambda},\orelloop \delta \mu]=1.
  \]
\end{defn}

\begin{lem}\label{lem:commuting-delta-loops}
  Let $\basepair\mu$ and $(\lambda,\ol{\lambda})$ be entangled by
  $\band{\delta}$.  If $[\orelloop{\delta}{\mu},
  \orelloop{\delta}{\lambda}] = 1$ then $\orelloop{\delta}{\mu}$ and
  $\orelloop{\delta}{\lambda}$ have the same axis in the dual
  Bass-Serre tree $\dualtree{t}{\mbc}$ and
  lie in a common cyclic subgroup.
\end{lem}
\begin{proof}[sketch]
  This is proved  the same way as Lemma \ref{lem:commuting-rel-class}.
\end{proof}

\begin{lem}\label{lem:wide-delta}
  In the band complex $\mbc$, let $\basepair\mu$ and
  $\basepair\lambda$ be block overlapping pairs let $\band\delta$
  satisfy $\delta \subset \inter\mu$ and
  $\dual\delta \subset \inter\lambda$. If for some
  $t \in \kappaTracks\mbc$
  \[ |\delta|_t > \max\left\{\tr\mu +\kappa, \tr\lambda+\kappa
  \right\},
  \] then $\basepair\mu$ and $\basepair\lambda$ are entangled by
  $\band\delta$. Furthermore the tubular elements $\orelloop \delta
  \lambda$ and $\orelloop \delta \mu$ have the same axis in $\dualtree
  t \mbc$.
\end{lem}
\begin{proof}[sketch]
  This is proved  the same way as Lemma
  \ref{lem:small-shifts-commute}.
\end{proof}

We can't simply widen some band $\band\lambda$
and still have a well formed band complex and preserve the dual
tree. We could do this in Lemma \ref{lem:widen}, by studying what
happened in the dual tree and by noting that the result was still a
well formed band complex. A similar analysis for block overlapping
pairs gives the following:

\begin{lem}\label{lem:widen-bop}
  Let $\basepair\mu$ be a block overlapping pair in a band complex
  $\mbc$ and let $t\in\kappaTracks{\mbc}$. Then we can obtain a new
  band complex $\mbc'\supset\mbc$, equipped with an efficiently
  carried $\kappa$-track $t'$, by widening $\band\mu$. After widening,
  $\mu$ has a initial naked segment and we recover $\mbc' \to \mbc$ by
  collapsing (Definition \ref{move:collapse}) this initial segment. The
  dual Bass-Serre trees $\dualtree{t'}{\mbc'}$ and $\dualtree{t}{\mbc}$ are
  therefore equivariantly isomorphic. Furthermore $\tr{\mu}$ remains
  invariant.
\end{lem}

\begin{lem}\label{lem:delta-widen}
  Let $(\mu,\ol{\mu})$ and $(\lambda,\ol{\lambda})$ be block
  overlapping pairs entangled by $\band{\delta}$ in a band complex
  $\mbc$.  Then for any $\kappa$-track $t$ efficiently carried by
  $\mbc$ we can obtain a new band complex $\mbc'\supset \mbc$ by first
  widening $\band\lambda, \band\mu$ by at most $|\lambda|_t+|\mu|_t$
  and then widening $\band\delta$ so that
  \begin{enumerate}[(i)]
  \item either $\inter\lambda \subset \delta$ or $\inter\mu \subset
    \dual\delta$, and
  \item the track $t$ extends to a track $t\subset t'$ efficiently
    carried by $\mbc'$ such that there is a $\fungrp{\mbc}$-equivariant
    isomorphism of dual Bass-Serre trees \[ \dualtree t \mbc \stackrel{\sim}{\to}
    \dualtree{t'}{\mbc'}.
    \]
  \end{enumerate}
\end{lem}

We note that $\lambda,\mu$ only need to be increased to length at most
$|\lambda|_t+|\mu|_t$ for this to work.

\subsubsection{The periodic block merger}
\begin{figure}[htb]
\centering
\begin{tikzpicture}[scale=0.12]
  
  \filldraw[black!10!white,draw=black,thick] (4,0) -- (7,0) -- (7,-5)
  -- (4,-5) --cycle;

  \fill[black!20!white] (2,0) -- (1,1) .. controls (3,3) and (2.8,6)
  .. (1,6) -- (20,6) .. controls (21.8,6) and (22,3) .. (20,1) --
  (19,0) -- cycle;
  \draw[thick] (1,1) .. controls (3,3) and (2.8,6)
  .. (1,6) -- (20,6) .. controls (21.8,6) and (22,3) .. (20,1) --
  (19,0) -- (2,0);
  \draw[dashed,thick] (20,6) .. controls (18.2,6) and (18,3) .. (20,1);
  
  \filldraw[black!10!white,draw=black,thick] (19,0) -- (20,1) -- (21,0)
  --cycle;
  \filldraw[black!10!white,draw=black,thick](1,1) .. controls (3,3) and
  (2.8,6) .. (1,6) (1,6) .. controls (-0.8,6) and (-1,3) .. (1,1);
  \fill[black!20!white] (0,0) -- (1,1) -- (2,0);
  \draw[dashed,thick] (1,1) -- (2,0);
  \draw[thick] (1,1) -- (0,0) -- (2,0);

  \begin{scope}[shift={(12,-5)},rotate=180]
    \fill[black!20!white] (2,0) -- (1,1) .. controls (3,3) and (2.8,6)
  .. (1,6) -- (10,6) .. controls (11.8,6) and (12,3) .. (10,1) --
  (9,0) -- cycle;
  \draw[thick] (1,1) .. controls (3,3) and (2.8,6)
  .. (1,6) -- (10,6) .. controls (11.8,6) and (12,3) .. (10,1) --
  (9,0) -- (2,0);
  \draw[dashed,thick] (10,6) .. controls (8.2,6) and (8,3) .. (10,1);
  
  \filldraw[black!10!white,draw=black,thick] (9,0) -- (10,1) -- (11,0)
  --cycle;
  \filldraw[black!10!white,draw=black,thick](1,1) .. controls (3,3) and
  (2.8,6) .. (1,6) (1,6) .. controls (-0.8,6) and (-1,3) .. (1,1);
  \fill[black!20!white] (0,0) -- (1,1) -- (2,0);
  \draw[dashed,thick] (1,1) -- (2,0);
  \draw[thick] (1,1) -- (0,0) -- (2,0);
  \end{scope}

  \begin{scope}[shift={(30,0)}]
     \filldraw[black!10!white,draw=black,thick] (1,0) -- (12,0) -- (12,-5)
  -- (1,-5) --cycle;

  \fill[black!20!white] (2,0) -- (1,1) .. controls (3,3) and (2.8,6)
  .. (1,6) -- (20,6) .. controls (21.8,6) and (22,3) .. (20,1) --
  (19,0) -- cycle;
  \draw[thick] (1,1) .. controls (3,3) and (2.8,6)
  .. (1,6) -- (20,6) .. controls (21.8,6) and (22,3) .. (20,1) --
  (19,0) -- (2,0);
  \draw[dashed,thick] (20,6) .. controls (18.2,6) and (18,3) .. (20,1);
  
  \filldraw[black!10!white,draw=black,thick] (19,0) -- (20,1) -- (21,0)
  --cycle;
  \filldraw[black!10!white,draw=black,thick](1,1) .. controls (3,3) and
  (2.8,6) .. (1,6) (1,6) .. controls (-0.8,6) and (-1,3) .. (1,1);
  \fill[black!20!white] (0,0) -- (1,1) -- (2,0);
  \draw[dashed,thick] (1,1) -- (2,0);
  \draw[thick] (1,1) -- (0,0) -- (2,0);

  \begin{scope}[shift={(12,-5)},rotate=180]
    \fill[black!20!white] (2,0) -- (1,1) .. controls (3,3) and (2.8,6)
    .. (1,6) -- (10,6) .. controls (11.8,6) and (12,3) .. (10,1) --
    (9,0) -- cycle;
  \draw[thick] (1,1) .. controls (3,3) and (2.8,6)
  .. (1,6) -- (10,6) .. controls (11.8,6) and (12,3) .. (10,1) --
  (9,0) -- (2,0);
  \draw[dashed,thick] (10,6) .. controls (8.2,6) and (8,3) .. (10,1);
  
  \filldraw[black!10!white,draw=black,thick] (9,0) -- (10,1) -- (11,0)
  --cycle;
  \filldraw[black!10!white,draw=black,thick](1,1) .. controls (3,3) and
  (2.8,6) .. (1,6) (1,6) .. controls (-0.8,6) and (-1,3) .. (1,1);
  \fill[black!20!white] (0,0) -- (1,1) -- (2,0);
  \draw[dashed,thick] (1,1) -- (2,0);
  \draw[thick] (1,1) -- (0,0) -- (2,0);
\end{scope} 
\end{scope}

  \begin{scope}[shift={(60,0)}]
    \fill[black!20!white] (2,0) -- (1,1) .. controls (3,3) and (2.8,6)
    .. (1,6) -- (20,6) .. controls (21.8,6) and (22,3) .. (20,1) --
    (19,0) -- cycle;
  \draw[thick] (1,1) .. controls (3,3) and (2.8,6)
  .. (1,6) -- (20,6) .. controls (21.8,6) and (22,3) .. (20,1) --
  (19,0) -- (2,0);
  \draw[dashed,thick] (20,6) .. controls (18.2,6) and (18,3) .. (20,1);
  
  \filldraw[black!10!white,draw=black,thick] (19,0) -- (20,1) -- (21,0)
  --cycle;
  \filldraw[black!10!white,draw=black,thick](1,1) .. controls (3,3) and
  (2.8,6) .. (1,6) (1,6) .. controls (-0.8,6) and (-1,3) .. (1,1);
  \fill[black!20!white] (0,0) -- (1,1) -- (2,0);
  \draw[dashed,thick] (1,1) -- (2,0);
  \draw[thick] (1,1) -- (0,0) -- (2,0);
  
  \begin{scope}[shift={(12,0)},rotate=180]
    \fill[black!20!white] (2,0) -- (1,1) .. controls (3,3) and (2.8,6)
    .. (1,6) -- (10,6) .. controls (11.8,6) and (12,3) .. (10,1) --
    (9,0) -- cycle;
    \draw[thick] (1,1) .. controls (3,3) and (2.8,6)
    .. (1,6) -- (10,6) .. controls (11.8,6) and (12,3) .. (10,1) --
    (9,0) -- (2,0);
    \draw[dashed,thick] (10,6) .. controls (8.2,6) and (8,3) .. (10,1);
    
    \filldraw[black!10!white,draw=black,thick] (9,0) -- (10,1) -- (11,0)
    --cycle;
  \filldraw[black!10!white,draw=black,thick](1,1) .. controls (3,3) and
  (2.8,6) .. (1,6) (1,6) .. controls (-0.8,6) and (-1,3) .. (1,1);
  \fill[black!20!white] (0,0) -- (1,1) -- (2,0);
  \draw[dashed,thick] (1,1) -- (2,0);
  \draw[thick] (1,1) -- (0,0) -- (2,0);
\end{scope}
\end{scope}
\end{tikzpicture}
\caption{Steps (1) and (2) of type periodic block
  merger.}\label{fig:tp2-merge}
\end{figure}
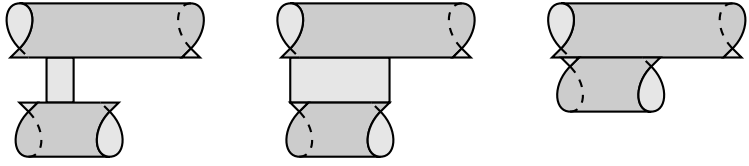

We now describe the \define{periodic block merger}, see Figure
\ref{fig:tp2-merge}. Let $\mbc$ have block overlapping pairs
$\basepair\mu$, $\basepair\lambda$ entangled by $\band{\delta}$. We do
the following:
\begin{enumerate}[(1)]
\item\label{it:tp1-widen} Widen bands $\band\lambda,\band\mu$ and
  then $\band{\delta}$ to obtain a band complex $\mbc'$ as given in
  Lemma \ref{lem:delta-widen}.
\item\label{it:tp2-widen} Assume that
  $\ol{\delta} \supset \inter{\lambda}$. We treat $\ol{\delta}$ as
  carrier base and move all the bases contained in $\inter{\lambda}$
  onto $\delta$ via entire transformations.
\item The dual pairs $\basepair{\lambda}, \basepair{\mu}$ are now
  entangled in the sense of Definition \ref{defn:entangled1} and both
  form overlapping pairs, so by Lemma \ref{lem:both-overlap-wide} we
  can apply the periodic merger given in Proposition
  \ref{prop:periodic-merger}, merging $\band{\lambda}, \band{\mu}$
  into some new $\band{\eta}$.
\end{enumerate}

We now record the following observation.

\begin{lem}\label{lem:um-tp2-merge}
  Let $\mbc$ be a band complex containing block overlapping pairs
  $\basepair\lambda, \basepair\mu$ entangled by $\band\delta$. Then we
  can effectively construct a finite set of band complexes\[
  \begin{tikzpicture}
    \node (1) at (0,0) {$\mbc$};
    \node (2) at (-1,-0.75) {$\mbc_1$};
    \node (3) at (1,-0.75) {$\mbc_n$};
    \node (4) at (0,-0.75) {$\cdots$};
    \draw[->] (1) -- (2);
    \draw[->] (1) -- (3);
  \end{tikzpicture}
  \]
  containing all possible combinatorial outcomes
  $m:(\mbc,t) \to (\mbc',t')$ of applying a periodic block merger as
  $t$ ranges over $\kappaTracks{\mbc}$.
\end{lem}

\subsubsection{Normalization: merging away entanglement}
It may be that $\basepair\mu$ is a block overlapping pair that is
entangled with $\basepair\lambda$ but $|\lambda|_t < \tr\mu$, for
example if $\lambda,\dual\lambda$ are very short and near the
extremities of $\basepair\mu$. In this situation it is not possible to
directly apply a periodic merger. However since block overlapping
pairs can always be widened we have the following:

\begin{lem}\label{lem:merge-entangled}
  Let $\mbc$ have a block overlapping pair $\basepair\mu$ overlapping
  pair entangled with $\basepair\lambda$. For any
  $t\in\kappaTracks{t}$ we can widen $\band{\mu}$ by $\tr\mu$ so that
  after widening $\band\lambda$ we have $|\lambda|_t\geq \tr\mu$.
\end{lem}
\begin{proof}
  By Lemma\ref{lem:widen-bop}, we can widen $\basepair\mu$ (by at most
  $\tr\mu$) so that $\lambda$ is at distance more than $\tr\mu$ from
  the endpoints of $\inter\mu$. It follows that after widening as in
  Lemma \ref{lem:widen}, $\lambda$ is sufficiently long.
\end{proof}

\begin{cor}\label{cor:merge-entangled}
  Let $\mbc$ have a block overlapping pair $\basepair\mu$ overlapping
  pair entangled with $\basepair\lambda$. For any
  $t\in\kappaTracks{t}$, after perhaps widening $\band\mu$ by at most
  $\tr{\mu}$ and $\band\lambda$ as in Lemma \ref{lem:widen}, we can
  perform a periodic merger (Proposition \ref{prop:periodic-merger})
  of $\band\mu$ and $\band\lambda$.
\end{cor}

\begin{defn}\label{defn:Normalized}
  Let $\mbc$ be a band complex in periodic block form and let
  $t \in \kappaTracks\mbc$, then we can repeatedly apply periodic
  block mergers and periodic mergers as described in Corollary
  \ref{cor:merge-entangled} so no block overlapping pair is entangled
  with any other base pair. If such is the case then we call
  corresponding $m:(\mbc,t) \to (\mbc',t')$ the \define{normalization
    of $(\mbc,t)$} and we say that $\mbc'$ is in \define{in normalized
    periodic block form}.
\end{defn}

We now give the combinatorial equivalent. Recall that entanglement is
algorithmically decidable and depends only on the band complex $\mbc$.

\begin{defn}\label{defn:normalized-desc}
  Let $\mbc$ be a band complex in periodic block form. Then its
  \define{normalized children} is the collection of band complexes
  
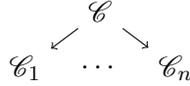
\[
  \begin{tikzpicture}
    \node (1) at (0,0) {$\mbc$};
    \node (2) at (-1,-0.75) {$\mbc_1$};
    \node (3) at (1,-0.75) {$\mbc_n$};
    \node (4) at (0,-0.75) {$\cdots$};
    \draw[->] (1) -- (2);
    \draw[->] (1) -- (3);
  \end{tikzpicture}
  \] obtained by enumerating all possible combinatorial outcomes of
  normalizations $m:(\mbc,t) \to (\mbc',t')$ where $t$ ranges over
  $\kappaTracks{\mbc}$.
\end{defn}

\subsection{Periodic hierarchies and maximal periodic
  blocks}\label{sec:periodic-hierarchies}
Throughout this section $\mbc$ will be a band complex in normalized
periodic block form.

\begin{defn}\label{defn:long}
  Let $\mbc$ be in normalized periodic block form, let
  $t \in \kappaTracks\mbc$, and let $\inter\lambda$ be a periodic
  block.  A base $\delta \subset \inter{\lambda}$ is \define{$t$-long in
    $\inter{\lambda}$} if $|\delta|_t \geq
  \tr{\lambda}+\kappa+1$. Otherwise it is called \define{$t$-short}.
\end{defn}

\begin{lem}\label{lem:no-long-bases}
  Let $\mbc, t,\lambda$, and $\delta$ be as in Definition
  \ref{defn:long}. Then if $\dual\delta$ also lies in $\inter\lambda$,
  $\delta$ must be $t$-short in $\inter\lambda$.
\end{lem}
\begin{proof}
  If $\basepair\delta$ is orientation reversing, then this follows
  immediately from Lemma \ref{lem:short-strips-II}. If $\basepair\delta$ is
  orientation preserving, but $t$-long in $\inter\lambda$, then it is
  entangled with $\basepair\lambda$, but not merged. By Corollary \ref{cor:merge-entangled}  this contradicts
  the assumption that $\mbc$ is in normalized periodic block form.
\end{proof}

It is possible for a base to be long in one periodic block, but its
dual must lie in another periodic block and it must be short in that
periodic block. If a base and its dual are both long in their
respective periodic blocks then we can perform a periodic block
merger.

\begin{defn}[Periodic hierarchies]\label{defn:periodic-hierarchy}
  Let the band complex $\mbc$ be in normalized periodic block form. A
  \define{periodic hierarchy} $\calH$ is a partial order $<_\calH$ on the set
  of periodic blocks that is generated as follows:
  \begin{enumerate}[(i)]
  \item\label{it:periodic-choice} If $\band\delta$ has bases lying in
    in periodic blocks $\inter{\lambda_1}$ and $\inter{\lambda_2}$ then we
    may either declare $\inter{\lambda_1} <_\calH \inter{\lambda_2}$,
    $ \inter{\lambda_2}<_\calH\inter{\lambda_1}$, or that
    $\inter{\lambda_1}$ and $\inter{\lambda_2}$ are incomparable.
  \item\label{it:partial-order} We extend (\ref{it:periodic-choice}) to a partial order, if
    possible.  
  \end{enumerate}
  A periodic block $\inter\lambda$ is \define{$\calH$-maximal} if it
  is maximal with respect to the partial order.  If
  $t \in \kappaTracks{\mbc}$ then we define the \define{induced
    periodic hierarchy} $\calH(t)$ to be generated by setting
  $\inter{\lambda_1} <_{\calH(t)} \inter{\lambda_2}$ if and only if
  $\delta$ is $t$-long in $\inter{\lambda_1}$ in
  (\ref{it:periodic-choice}).
\end{defn}

It is obvious that periodic hierarchies, being finite combinatorial
objects, can be effectively listed. What is less obvious is whether
the definition of an induced periodic hierarchy actually gives a
periodic hierarchy.

\begin{lem}
  If $\mbc, t$ and $\calH(t)$ are as in Definition
  \ref{defn:periodic-hierarchy} then $<_{\calH(T)}$ gives a partial
  order on the set of periodic blocks; thus $\calH(t)$ is a periodic
  hierarchy.
\end{lem}
\begin{proof}
  Suppose that for some $\band\delta$, both $\delta$ and $\dual\delta$
  were long in $\inter{\lambda_1}$ and $\inter{\lambda_2}$,
  respectively. Then by Lemma \ref{lem:wide-delta} the block
  overlapping pairs are entangled and can be merged, contradicting the
  assumption that $\mbc$ is normalized, so (\ref{it:periodic-choice})
  of Definition \ref{defn:periodic-hierarchy} is satisfied. Thus
  (\ref{it:periodic-choice}) gives a directed graph $\Gamma$ without
  loops of length 2 with periodic blocks as vertices. If $\Gamma$ has
  a directed cycle, then this would imply that for some
  $\basepair\lambda$, $\tr\lambda < \tr\lambda$ which is
  absurd. We can therefore extend $\cal_H$ to a
  partial order giving (\ref{it:partial-order}) of Definition \ref{defn:periodic-hierarchy}.
\end{proof}

\begin{cor}\label{cor:max-short}
  If $\basepair\lambda$ is an $\calH(t)$-maximal block overlapping
  pair, then every base $\delta \subset \inter\lambda$ is $t$-short.
\end{cor}

\subsection{Bounding the periodicity of maximal periodic blocks}
\begin{defn}
  Suppose $\inter{\lambda}$ is a periodic block and suppose that we
  can vertically subdivide $\band{\lambda}$ into three
  bands \[\band{\lambda_1},\band{\lambda_2}, \band{\lambda_3}\] such
  that: \begin{itemize}
  \item $(\lambda_2,\ol{\lambda_2})$ form an overlapping pair.
  \item No bases other than $\ol{\lambda_1}$ and $\lambda_3$ intersect
    $\inter{\lambda_2}$.
  \item $\band{\lambda_2}$ contains no connections (recall Definition
    \ref{defn:MBC}(\ref{it:conn-a}).)
  \end{itemize}
  Then we call $\band{\lambda_2}$ a \define{clean
    tube.}
\end{defn}

\begin{figure}[htb]
  \centering
  \begin{tikzpicture}[scale=0.8, every node/.style={scale=0.8}]
      \node [draw, cloud, cloud puffs=7, minimum width = 5cm, minimum height=5cm, fill=gray!60]
      at (-2,0) {};
      \node [draw, cloud, cloud puffs=4, minimum width = 5cm, minimum height=5cm, fill=gray!60]
      at (6,0) {};
      \draw[fill=white] (0,-1) .. controls +(-0.5,0) and +(-0.5,0) .. (0,1) -- 
      (4,1) ..controls +(0.5,0) and +(0.5,0) ..(4,-1) -- (0,-1);
      \begin{scope}
        \clip[draw] (0,-1) .. controls +(-0.5,0) and +(-0.5,0) .. (0,1) -- 
        (4,1) ..controls +(0.5,0) and +(0.5,0) ..(4,-1) -- cycle;
        \draw[dashed] (-1,-1) ..controls +(1,0.5) and +(-1,-0.5).. (1,1);
        \draw[dashed] (1,-1) ..controls +(1,0.5) and +(-1,-0.5).. (3,1);
        \draw[dashed] (3,-1) ..controls +(1,0.5) and +(-1,-0.5).. (5,1);
        \begin{scope}[xshift=1cm]
          \draw[dashed] (-1,-1) ..controls +(1,0.5) and +(-1,-0.5).. (1,1);
          \draw[dashed] (1,-1) ..controls +(1,0.5) and +(-1,-0.5).. (3,1);
        \end{scope}
      \end{scope}
  \end{tikzpicture}
  \caption{A clean tube in a band complex is literally an embedded
    $S^1\times [-1,1]$. The track $t$ spirals around the clean
    tube. Unwinding it by a Dehn twist will decrease the size.}
  \label{fig:clean-tube}
\end{figure}
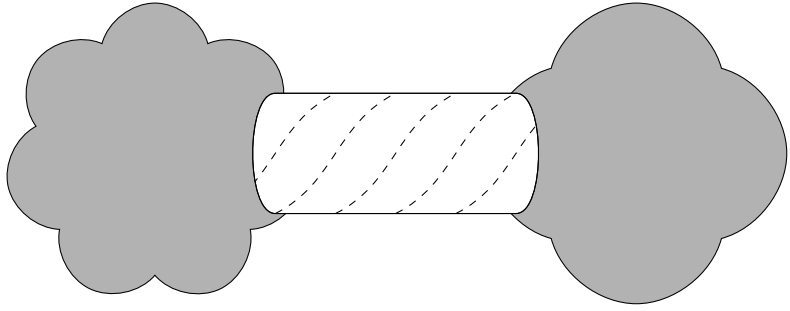

The significance of clean tubes is illustrated in Figure
\ref{fig:clean-tube}.
Although the next result is obvious from this picture. It is important
to state it carefully to get explicit bounds.
\begin{lem}\label{lem:dehn-twist}
  Let $t$ be an automorphically minimal $\kappa$-track efficiently
  carried by $\mbc$, and let $\inter{\lambda}$ be a periodic block.
  Then for any clean tube $\band{\lambda_2} \subset \band{\lambda}$,
  \[|\lambda_2|_t\leq 2\tr{\lambda}.\] 
\end{lem}
\begin{proof}
  Suppose towards a contradiction that for some clean tube
  $\band{\lambda_2}$ we have: \[|\lambda_2|_t> 2\tr{\lambda}.\]
  Let $I \subset \lambda_2 \cap \ol{\lambda_2}$ be interval of length
  $\tr{\lambda}$ that is of distance $\tr{\lambda}$ the leftmost
  endpoint of of $\inter{\lambda_2}$. Let $p \in t \cap I$ be the
  rightmost point of $t\cap I$, if we follow the the connected
  component of $t\cap \band{\lambda_2$} that contains $p$ and
  intersects $I$ again in $p'$ we see that the distance between $p$
  and $p'$ in $\lambda_2$ is exactly $\tr{\lambda}$. Recall the
  notation of Definition \ref{defn:measured-band} and
  consider the map
  \[ J_{\band{\lambda_2}} \times [-1,1] \rightarrow
  \mbc \]
  with $J_{\band{\lambda_2}} \times\{1\} \rightarrow \lambda_2$ and
  $J_{\band{\lambda_2}}\times\{-1\} \rightarrow \ol{\lambda_2}$. The
  preimage of $I$ has two connected components
  $I_\pmo \subset J_{\band{\lambda_2}} \times\{\pmo\}$. Let $\alpha$
  be the straight line in $J_{\band{\lambda_2}} \times [-1,1]$ between
  the rightmost point of $I_1$ and the rightmost point of $I_{-1}$ and
  let $\beta$ be the line in $J_{\band{\lambda_2}} \times [-1,1]$
  between the leftmost point of $I_1$ and the leftmost point of
  $I_{-1}$.  $\alpha$ and $\beta$ are chosen to be transverse to the
  preimage of $t$.  Let $Q$ be the quadrilateral in
  $J_{\band{\lambda_2}} \times [-1,1]$ enclosed by $I_{\pmo},\alpha,\beta$
  (see Figure \ref{fig:quadrilateral}.)  Then via
  $J_{\band{\lambda_2}} \times [-1,1] \rightarrow \mbc$, $Q$ is
  mapped to an annulus $A \subset \mbc$
  such that $Q \subset \band{\lambda_2}$ with $\alpha,\beta$ mapping
  onto each component of $\partial A$.
  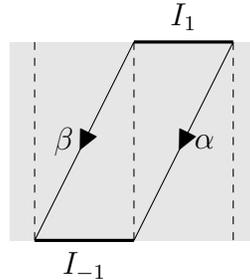
\begin{figure}[htb]
    \centering
    \begin{tikzpicture}[scale={0.66}]
      \fill[black!10!white] (-0.5,0) -- (4.5,0) -- (4.5,-4) -- (-0.5,-4)
      --cycle; 
      \draw[dashed] (2,0) -- (2,-4); 
      \draw[dashed] (0,0) -- (0,-4); 
      \draw[dashed] (4,0) -- (4,-4); 
      \draw (4,0) -- node[sloped]{$\blacktriangleleft$} node[right]{$\alpha$} (2,-4);
      \draw (2,0) -- node[sloped]{$\blacktriangleleft$}
      node[left]{$\beta$}(0,-4); \draw[very thick] (2,0) --
      node[above]{$I_1$} (4,0); \draw[very thick] (0,-4) --
      node[below]{$I_\mo$} (2,-4);
    \end{tikzpicture}
    \caption{The quadrilateral $Q$ inside $\band{\lambda_2}$. The
      annulus $A$ is obtained by identifying $I_1$ and $I_{-1}$. The
      track is drawn as dashed lines.}\label{fig:quadrilateral}
  \end{figure}

  We parameterize this annulus $A$
  as \begin{equation}\label{eqn:parametrize-annulus} \{r \exp(i\theta)
    | 1\leq r \leq 2, 0 \leq \theta \leq 2\pi\} \subset
    \mathbb{C}.\end{equation} By construction, $t_A=t \cap A$ corresponds
  to the curve:
  \begin{eqnarray*}
    t_A: [0,1] & \rightarrow & A\\
    s & \mapsto & (2-s)\exp(i(2\pi s)).\\
  \end{eqnarray*}
  If we make a Dehn $\tau_A$ twist around $A$ then post-composing
  gives:
  \begin{eqnarray*}
    \tau_A\circ t_A: [0,1] & \rightarrow & A\\
    s & \mapsto & (2-s).\\
  \end{eqnarray*}  
  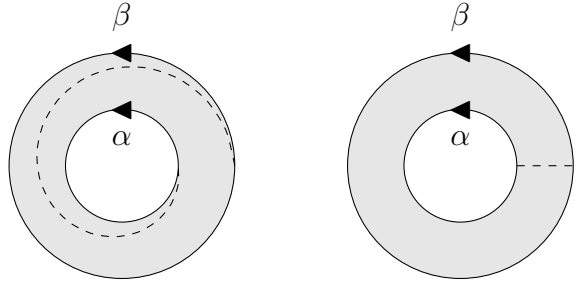
\begin{figure}[htb]
    \centering
    \begin{tikzpicture}[scale={0.75}]
      \fill[black!10!white] (0,0) circle (2);
      \draw[domain=0:1,smooth,samples=20,variable=\t,dashed] plot
      ({(2-\t)*cos(2*pi*\t r)},  {(2-\t)*sin(2*pi*\t r)})
      coordinate(end);
      \fill[white] (0,0) circle (1);
      \draw (0,0) circle (1) (0,0) circle (2) (0,1) node{$\blacktriangleleft$} (0,1)++(0,-0.2)
      node[below]{$\alpha$} (0,2) node{$\blacktriangleleft$} (0,2)++(0,0.2)
      node[above]{$\beta$};
      \fill[shift={(6,0)},black!10!white] (0,0) circle (2);
      \draw[dashed,shift={(6,0)}] (0:1) -- (0:2);
      \fill[shift={(6,0)},white] (0,0) circle (1);
      \draw[shift={(6,0)}] (0,0) circle (1) (0,0) circle (2) (0,1) node{$\blacktriangleleft$} (0,1)++(0,-0.2)
      node[below]{$\alpha$} (0,2) node{$\blacktriangleleft$} (0,2)++(0,0.2)
      node[above]{$\beta$};
    \end{tikzpicture}
    \caption{On the left the image $A$ of $Q$ after identifying $I_1$
      and $I_2$ parameterized as (\ref{eqn:parametrize-annulus})
      (drawn to scale). On the right the result of the Dehn twist
      $\tau_A$. The the dashed line represents $t_m \cap
      A$. Note that $\tau_A$ restricts to the identity on $\partial
      A$}\label{fig:twist}
  \end{figure}

  Such a Dehn twist is illustrated in Figure \ref{fig:twist}. Since
  the homeomorphism $\tau_A$ restricts to the identity on
  $\partial A$, it extends to a homeomorphism of $\mbc$ which we will
  also denote $\tau_A$. Consider the new track $\tau_A(t)$. On one
  hand, $\tau_A(t)$ is still efficiently carried by $\mbc$. On the
  other hand, after perturbing by an isotopy we have,
  $\tau_A(t) \cap \inter{\lambda} = t \cap \inter\lambda -1$, but
  otherwise for every subset
  $S \subset \mbc \setminus \band{\lambda_2}$ we have
  $S\cap t = S \cap \tau_A(t)$. It follows that
  \begin{equation}
    \label{eq:size-matters}
    \size{\tau_A(t)} < \size{t}
  \end{equation}
  The homeomorphism $\tau_A:\mbc \to \mbc$ lifts to a homeomorphism of
  $\tilde{\tau_A}:\tilde{\mbc} \to \tilde{\mbc}$ sending the lift
  $\tilde{t}$ to the lift $\tilde{\tau_A(t)}$. Furthermore $\tau_A$
  maps naturally to an element of $\aut{(\fungrp{\mbc})}$. It follows
  that the trees $\dualtree{t}{\mbc}$ and $\dualtree{\tau_A(t)}{\mbc}$
  are $(\tau_A)_\sharp$-equivariantly-isomorphic so that
  $t\sim_{\aut{(\fungrp{\mbc})}}\tau_A(t)$. $\tau_A(t)$ is therefore
  obviously a $\kappa$-track, and with (\ref{eq:size-matters}) we see
  that $t$ is not automorphically minimal, which is a contradiction.
 \end{proof}

 Informally, if clean tubes are longer than twice the translation
 length, there would be enough room to make an annulus, as shown in
 Figures \ref{fig:quadrilateral}, \ref{fig:twist}. This would enable
 us to shorten the track $t$ by a Dehn twist, contradicting
 automorphic minimality. Clean tubes therefore must be short. We now
 use this observation to give a combinatorial bound for periodicity.

\begin{cor}[Periodicity bound for maximal periodic blocks]\label{cor:periodicity-bound}
  Let $t\subset \mbc$ be an efficiently carried automorphically
  minimal $\kappa$-track and let $\inter\lambda$ be a
  $\calH(t)$-maximal periodic block. Let $|C|$ denote the number of
  connected components of connection preimages in $\band\lambda =
  J_{\band{\lambda}}\times [-1,1]$ and
  let $B$ denote the number of bases, other than $\lambda,
  \dual\lambda$ contained in $\inter\lambda$. Then \[
  \period{t}{\lambda} \leq (6+2\kappa)B+2|C|+2.
  \]
\end{cor}
\begin{proof}
  Parameterize $\band\lambda$ as $J_{\lambda} \times [-1,1]$ so that
  the preimage of $t$ and of every connection is contained in a union
  of vertical lines. Every base $\delta \subset \band \lambda$ has a
  preimage with connected components $\delta^+ \subset J_\lambda
  \times \{1\}$ and $\delta^- \subset J_\lambda \times \{-1\}$. Let
  $C$ denote the preimage of the connections. Consider the
  complement \[
  Y = J_{\band{\lambda}} \times [-1,1] \setminus \left( C \cup \left(\bigcup_\delta \left(\delta^\pm
  \times [-1,1]\right)\right)\right),
\] where $\delta$ runs over the bases contained in $\inter\lambda$. On
one hand $Y$ has at most $|C| + 2B+1$ connected components, on the
other hand every maximal clean tube contained in $\band{\lambda}$ is
the preimage of one of these components, thus every connected
component of $Y$ has length at most $2\tr\lambda$ be Lemma
\ref{lem:dehn-twist}. Furthermore, by Corollary \ref{cor:max-short}
and Definition \ref{defn:long}, every
$\delta^\pm$ has width at most $\tr{\lambda}+\kappa$. We therefore
have the bound:\begin{eqnarray*}
  |\lambda|_t &\leq& 2B\left(\tr{\lambda}+\kappa\right) + 2\tr{\lambda}\left(|C| + 2B+1\right)\\
  & \leq & 2B\left(\tr{\lambda}+\tr{\lambda}\kappa\right) +
  2\tr{\lambda}\left(|C| + 2B+1\right)\\
& \leq & \tr{\lambda}\left((6+2\kappa)B+2|C|+2\right);
  \end{eqnarray*}
  from which we immediately obtain the advertised bound.
\end{proof}

The significance of this bound is that it only depends on
the combinatorial band complex $\mbc$ and the combinatorial periodic
hierarchy $\calH(t)$. We  obtain the following computable
function.

\begin{defn}\label{defn:period-lambda}
  Let $\mbc$ be in normalized periodic block form, let $\calH$ be a
  periodic hierarchy on $\mbc$, and let $\inter\lambda$ be
  $\calH$-maximal periodic block.  We
  define\[\period{\calH}{\lambda} = (6+2\kappa)B+2|C|+2\]
  where $B,|C|$ are as in Corollary \ref{cor:periodicity-bound}.
\end{defn}

Unfortunately it may be that the periodic block whose periodicity we
are interested in is not $\calH$-maximal. We deal with this in the
next section.

\subsection{Bounding the periodicity of principal overlapping
  carriers: auxiliary trees}\label{sec:taux}
Throughout this section we will use the following notation. If $\mbc$
is a band complex, then the term $J$ in the pair $(\mbc,J)$ will
always denote a union of block overlapping pairs. If $\mbc$
happens to be in normalized periodic block form and is equipped with a
periodic hierarchy $\calH$, then we will denote the corresponding
triple $(\mbc,J;\calH)$.

We will try to bound the periodicity of the block overlapping pairs
that constitute $J$. It may happen however that two block overlapping
pairs get merged at some point. Because of this we will use the
following naming convention.

\begin{conv}[Renaming merged bases in auxiliary trees]\label{conv:renaming}
  Suppose that two bands $\band{\lambda},\band\mu$ get merged onto
  some band $\band \eta$ in a periodic merger. Then, as far as naming
  bases is concerned, we will consider $\eta = \lambda = \mu$, i.e. we
  will allow a base to have the multiple names.
\end{conv}

This renaming convention is justified since on one hand, we want to
bound the periodicity of an overlapping pair that gets merged, so we
must to keep track of what it got merged with. On the other hand by
the following result, which is an immediate consequence of Position
\ref{prop:periodic-merger}, we are guaranteed that any periodicities
we compute will be overestimates. 

\begin{lem}\label{lem:periodicity-non-decr}
  Let $\basepair\lambda$ be an overlapping pair such that
  $\band\lambda$ gets zipped onto $\band\eta$ via a periodic merger
  $\mbc \to \mbc'$ mapping an efficiently carried $\kappa$-track
  $t \subset \mbc$ into the efficiently carried $\kappa$-track
  $t'\subset \mbc'$.
  Then\[ \period{t}{\lambda} \leq \period{t'}{\eta}.  \]
  In particular periodicities of block overlapping pairs are
  nondecreasing when passing to normalizations.
  \end{lem}

\begin{defn}\label{defn:aux-descendants}
  Let $(\mbc,J)$ be a band complex with a principal overlapping pair
  $\basepair\lambda$. Then we define the \define{auxiliary
    children} of $(\mbc,J)$ to collection:\[
  \begin{tikzpicture}
    \node (mbc) at (0,0) {$(\mbc,J)$};
    \node (ad1) at (-2,-1) {$(\mbc_1,J\cup\inter\lambda;\calH_1)$};
    \node (ad2) at (2,-1) {$(\mbc_N,J\cup\inter\lambda;\calH_N)$};
    \node (dots) at (0,-1) {$\cdots$};
    \draw[->] (mbc) -- (ad1);
    \draw[->] (mbc) -- (ad2);
  \end{tikzpicture}
  \]
  where
  \[\left\{(\mbc_1,J\cup\inter\lambda;\calH_1), \ldots,
    (\mbc_N,J\cup\inter\lambda;\calH_N)\right\}\]
  is obtained by first taking the leaves of $\tpbf{\mbc}{J}$, then
  taking their normalized children (Definition
  \ref{defn:normalized-desc}), and finally by taking all combinatorial
  possibilities for periodic hierarchies (Section
  \ref{sec:periodic-hierarchies}).
\end{defn}

In Section \ref{sec:periodic-hierarchies} we bounded the periodicity
of a maximal block overlapping pair, but it may be that
$\inter\lambda$, where $\lambda$ is a principal overlapping pair in
$(\mbc,J)$, is not a maximal periodic block in some of the auxiliary
descendant $(\mbc_i,J\cup\inter\lambda,\calH_i)$, due to the choice of
periodic hierarchy $\calH_i$. It follows that simply passing to
auxiliary children isn't sufficient to bound $\period{t}{\lambda}$
where $t$ ranges over the automorphically minimal $\kappa$-tracks
efficiently carried by $\mbc$.

In Section \ref{sec:aux-tree} we will construct from $(\mbc,J)$ a
finite auxiliary tree $\Taux{\mbc,J\cup\inter\lambda}$. If $\mbc$
efficiently carries an automorphically minimal $\kappa$-track $t$,
then we will construct an induced tree
$\Taux{\mbc,t,J\cup\inter\lambda}$ in Section \ref{sec:ind-tree}. This
induced tree will be proved to contain a combinatorial witness for
an upper bound of $\period{t}{\mu}$ for some base $\mu$ such that
$\inter\mu \subset J$ is a periodic block. We will also have a
containment\[ \Taux{\mbc,t,J\cup\inter\lambda} \subset
\Taux{\mbc,J\cup\inter\lambda}.
\]
From this it will follow that the construction of
$\Taux{\mbc,\inter\lambda}$ will give a way to bound
$\period{t}{\lambda}$ in $\mbc$ where $t$ ranges over the
automorphically minimal $\kappa$-tracks.

\subsubsection{The auxiliary tree $\Taux{\mbc,J\cup\inter\lambda}$}\label{sec:aux-tree}

We construct $\Taux{\mbc,J\cup\inter\lambda}$ with the following
recursive algorithm. The reader may skip ahead to Figure
\ref{fig:ind-aux-tree} to get an idea of what this tree is supposed to
look like.

\begin{enumerate}[(1)]
  \setcounter{enumi}{-1}
  
\item If $\basepair\lambda$ is a principal overlapping pair in
  $(\mbc,J)$ then we declare $(\mbc,J)$ to be the \define{root of}
  $\Taux{\mbc,J\cup\inter\lambda}$. Write $J' = J'\cup\inter\lambda$.
\item\label{it:aux-descendants} The $\Taux{\mbc, J'}$-children of
  the root $(\mbc,J)$ are the auxiliary children of $\mbc$
  (Definition \ref{defn:aux-descendants}). These are connected to the
  root by \define{auxiliary edges}. Further descendants are added as
  follows:
  \begin{enumerate}
  \item\label{it:taux-terminal} If some periodic block
    $\inter\mu \subset J'$ in an auxiliary child
    $(\mbc',J';\calH')$ of $(\mbc,J)$ is $\calH'$-maximal, then
    $(\mbc',J';\calH')$ is called a \define{witnessing terminal}. We
    stop growing $\Taux{\mbc,J'}$ at $(\mbc',J')$.
  
  \item\label{it:taux-reorder-and-go} Otherwise we modify the order
    $<'$ on $\mbc'$ (Definition \ref{defn-ordering}) so that some
    $\calH'$-maximal periodic block $\inter\delta$ is initial
    with $\delta$ the  carrier base.

    We start building $\ETtwo{\mbc,J'}$ rooted at $(\mbc',J')$. For
    every path originating from the root we forbid the base $\delta$
    from being the carrier base more than
    {$B\cdot\period{\calH'}{\delta}$} times in a row, where $B$
    denotes the number of bases in $\mbc'$. Once $\delta$ ceases to be
    the carrier base we forget about $\calH'$, and go to
    (\ref{it:horizontal-growth}) below.
  \end{enumerate}
\item\label{it:horizontal-growth}
  We continue growing $\Taux{\mbc,J'}$ using the
  following rules:
  \begin{enumerate}
  \item\label{it:std-growth} If the leading base $\mu$ in some
    $(\mbc'',J')$ does not form a principal overlapping pair, then
    continue growing $\ADtwo{\mbc',J'}$ at $(\mbc'',J')$ by adding
    admissible descendants (Definition \ref{defn:admissible}), if
    there are any.
    
    If $(\mbc'',J')$ has no admissible children and some of its
    bases are not in $J$, then $(\mbc'',J')$ is \define{halted by
      inadmissibility}.
    
    Otherwise, if the elimination process stops because all the bases
    of $\mbc''$ were moved into $J'$, then we take the auxiliary
    children
    \[\left\{(\mbc''_r,J';\calH_r)\right\}\]
    of $(\mbc'',J')$ (which is already in periodic block from, but
    may not be normalized), equipped with periodic hierarchies. All
    these children are \define{witnessing terminals} as in
    (\ref{it:taux-terminal}).
    
  \item\label{it:overlapping-case} If the carrier base base $\mu$ in
    some $(\mbc'',J')$ is a principal overlapping carrier we first
    construct $\Taux{\mbc'',J''}$ rooted at $(\mbc'',J')$, where
    $J'' = J'\cup\inter\mu$. This is the recursion.

    Next we take
    \[W = W(\mbc'',J'',\mu)=\left\{(\mbc_i,J_i;\calH_i)\right\}\]
    to be the set of all witnessing terminals of $\Taux{\mbc'',J''}$
    in which $\inter\mu$ is an $\calH_i$-maximal block overlapping
    pair (recall Convention \ref{conv:renaming}).
    
    If this set is empty then $(\mbc'',J')$ is declared to be a
    \define{halted terminal} and no further descendants are
    added. Otherwise the following number is defined and computable:
    \begin{equation}\label{eq-period-bound}
      \period{(\mbc'',J')}{\mu}=\max_{(\mbc_i,J_i;\calH_i)\in
        W}\period{\calH_i}{\mu}.      
    \end{equation}
    We continue growing $\ADtwo{\mbc',J'}$ at $(\mbc'',J')$ but we
    forbid $\mu$ from being the carrier base more than
    $B\cdot\period{(\mbc'',J')}{\mu}$ times in a row in every path
    originating at $(\mbc'',J')$. Once $\mu$ ceases to be a
    carrier base we go back to (\ref{it:std-growth}) or
    (\ref{it:overlapping-case}) as appropriate.
  \end{enumerate}
  
\end{enumerate} 

\begin{defn}\label{defn:depth}
  Let $\mbc'$ be a band complex occurring in $\Taux{\mbc,J}$ the
  \define{depth of $\mbc'$ in $\Taux{\mbc,J}$} is the number of
  auxiliary edges in $\Taux{\mbc,J}$ connecting $\mbc$ and $\mbc'$.
\end{defn}

\begin{lem}\label{lem:depth-bound}
  The maximal depth of a descendant of $\mbc$ in $\Taux{\mbc,J}$ is at
  most the relative $\tau$-complexity $\tau(\mbc,J)$.
\end{lem}
\begin{proof}
  By Lemma \ref{lem:pb-complexity} all the leaves of $\tpbf{\mbc,J}$
  have strictly smaller $\tau$ complexity relative to
  $J \cup \inter\lambda$. Furthermore, periodic mergers never increase
  $\tau$-complexity. Finally if the relative $\tau$ complexity
  $\tau(\mbc',J')=0$, then by (\ref{it:std-growth}) all its children
  are witnessing terminals.
\end{proof}

We think of auxiliary edges as being vertical; thus

\begin{defn}\label{defn:horizontal}
  A subtree of $\Taux{\mbc,J}$ sitting inside some $\ADtwo{\mbc',J'}$
  or, equivalently, without auxiliary edges is called
  \define{horizontal}.
\end{defn}

\begin{prop}\label{prop:taux-finite}
  Let $\lambda$ be a principal overlapping pair in a band complex
  $\mbc$ occurring in $\ADtwo{\mbc}$, then
  $\Taux{\mbc,\inter{\lambda}}$ is finite.
\end{prop}
\begin{proof}
  We first to show that for any $\mbc''$ occurring in some horizontal
  subtree $\ADtwo{\mbc',J'}\subset \Taux{\mbc,\inter{\lambda}}$ with a
  principal overlapping carrier $\mu$, the auxiliary tree
  $\Taux{\mbc'',J'\cup\inter\mu}$ is finite.

  We prove this by induction on $\tau(\mbc'',J')$. If
  $\tau(\mbc'',J')=0$ then if it is not itself yet a witnessing
  terminal, then its auxiliary children are halted terminals and
  the result follows. Otherwise if
  $\tau(\mbc'',J')=1$, by Lemma \ref{lem:pb-complexity}, all auxiliary
  children must have $(J' \cup \inter\lambda)$-relative
  $\tau$-complexity equal to 0, so they are witnessing terminals.

  Now we suppose that all auxiliary trees are finite for all relative
  $\tau$ complexities less than $n$, and that $\tau(\mbc'',J')=n$. Any
  auxiliary descendant $(\mbc''',J'';\calH''')$ will have
  $\tau(\mbc''',J'') < n$ by Lemma \ref{lem:pb-complexity}. We
  construct $\ADtwo{\mbc''',J''}$ for each auxiliary descendant of
  $(\mbc'',J')$ according to rules (\ref{it:std-growth}) and
  (\ref{it:overlapping-case}).  Whenever a principal overlapping
  carrier $\eta$ occurs in some $\mbc^{(4)}$, the corresponding
  auxiliary tree built in (\ref{it:overlapping-case}) is finite by the
  induction hypothesis. Either $\mbc^{(4)}$ is a halted terminal or we
  can compute the finite $\period{(\mbc^{(4)},J'')}{\mu}$ (as given in
  (\ref{eq-period-bound}) of step (\ref{it:overlapping-case})). By Proposition
  \ref{prop:periodicity-reduction}, this prevents horizontal subtrees
  from having infinite branches; thus by König's Lemma they are
  finite. It follows that $\Taux{\mbc'',J'\cup\inter\lambda}$ is
  finite.  The result now follows by induction.
\end{proof}

Having established that $\Taux{\mbc,J}$ is finite, and therefore
effectively constructible, we can now define the following computable
function.

\begin{defn}\label{def:periodicity-bound}
  For a band complex $\mbc$ with a principal overlapping pair
  $\basepair\lambda$ we
  define\[ \period{\mbc}{\lambda} = \max_{(\mbc_i,J_i,\calH_i)\in
    S}\period{\calH_i}{\lambda} \]
  where $S$ is the set of witnessing terminals $(\mbc_i,J_i,\calH_i)$
  of $\Taux{\mbc,\inter\lambda}$ where, following the renaming
  Convention \ref{conv:renaming}, $\inter\lambda$ is
  $\calH_i$-maximal.
\end{defn}

It remains to show that this $\period{\mbc}{\lambda}$ gives an upper
bound for $\period{t}{\lambda}$ where $t$ is an automorphically
minimal $\kappa$-track efficiently carried by $\mbc$. This will be
done by studying the induced tree.

\subsubsection{The induced tree $\Taux{\mbc,t,J\cup\inter\lambda}$} \label{sec:ind-tree}
A triple $(\mbc,t,J)$ will denote a band complex $\mbc$, an
automorphically minimal $\kappa$-track $t$ efficiently carried by
$\mbc$, and a union $J$ of periodic blocks. If $\mbc$ is in normalized
periodic block form,  in the 3+1 tuple
$\left(\mbc,t,J;\calH(t)\right)$, $\calH(t)$ will denote the periodic
hierarchy induced by $t$ (Definition \ref{defn:periodic-hierarchy}.) 

In any (restricted) elimination tree $\mathfrak{T}$ rooted at
$(\mbc,J)$ a track $t\subset \mbc$ efficiently carried by $\mbc$
induces a directed path in $\mathfrak{T}$. By Corollary
\ref{cor:tpbf-finite} and by Definition \ref{defn:Normalized} the
following makes sense:

\begin{defn}\label{defn:ind-aux-descendants}
  Let $(\mbc,t,J)$ be a band complex with a principal overlapping pair
  $\basepair\lambda$ then its \define{induced auxiliary child} is
  given by the labelled graph\[
  \begin{tikzpicture}
    \node (mbc) at (0,0) {$(\mbc,J,t)$};
    \node (ad1) at (0,-1) {$(\mbc',t',J\cup\inter\lambda;\calH(t'))$};
    \draw[->] (mbc) -- (ad1);

  \end{tikzpicture}
  \]
  where $(\mbc',t,J\cup\inter\lambda;\calH(t'))$ is the auxiliary
  descendant of $(\mbc,J)$ (Definition \ref{defn:aux-descendants})
  induced by $t \subset \mbc$.
\end{defn}

The algorithm to construct $\Taux{\mbc,t,J\cup\inter\lambda}$ is
analogous to the algorithm to construct
$\Taux{\mbc,J\cup\inter\lambda}$ given in Section
\ref{sec:aux-tree}. The numbering of the clauses is intended to
coincide.

\begin{enumerate}[(1)]
  \setcounter{enumi}{-1}
\item If $\basepair\lambda$ is a principal overlapping pair in
  $(\mbc,t,J)$ then we declare $(\mbc,t,J)$ to be the \define{root of}
  $\Taux{\mbc,t,J\cup\inter\lambda}$. Write $J' = J'\cup\inter\lambda$.

\item\label{it:ind-aux-descendants} The $\Taux{\mbc,t,J'}$-child
  of the root $(\mbc,t,J)$ is the induced auxiliary child of
  $\mbc$. These are
  connected to the root by an \define{auxiliary edge}. Further
  descendants are added as follows:

  \begin{enumerate}
  \item\label{it:ind-taux-terminal} If some periodic block
    $\inter\mu \subset J'$ in an auxiliary descendant
    $\left(\mbc',t',J';\calH(t')\right)$ of $(\mbc,t,J)$ is
    $\calH(t')$-maximal, then $\left((\mbc',J';\calH(t')\right)$ is
    called a \define{witnessing terminal}. We stop growing
    $\Taux{\mbc,t,J'}$.
    
  \item\label{it:taux-reorder-and-go} Otherwise we change the order
    $<$ on $\left(\mbc',t',J';\calH(t')\right)$ so that
    some $\calH(t')$-maximal periodic block $\inter\delta$
    is initial with $\delta$ the carrier base.
  
    We start our Rips process, building the path in $\ADtwo{\mbc,J'}$
    rooted at $(\mbc',J')$ induced by the track $t' \subset \mbc'$. By
    Corollary \ref{cor:periodicity-bound}, $\delta$ is not the carrier
    base more than {$B\cdot\period{\calH'}{\delta}$} times in a row,
    where $B$ denotes the number of bases in $\mbc'$. Once $\delta$
    ceases to be the carrier base we go to
    (\ref{it:ind-horizontal-growth}) below.
  \end{enumerate}

\item\label{it:ind-horizontal-growth} We continue growing the path
  induced by $t$ in $\ADtwo{\mbc,J'}$ using the following rules:
  \begin{enumerate}
  \item\label{it:ind-std-growth} If the carrier base $\mu$ in some
    $(\mbc'',t'',J')$ is not a principal overlapping carrier, then
    we add its descendant in $\ADtwo{\mbc,J'}$ as usual.
    
    If all the bases $\mbc''$ are moved into $J'$, then we take the
    induced auxiliary child
    \[\left( \mbc''',t''',J';\calH(t''') \right) \]
    of $(\mbc'',J')$ (which is already in periodic block from, but may
    not be normalized). Again we call
    $ \left(\mbc''',t''',J';\calH(t''')\right)$ a \define{witnessing
      terminal} as in (\ref{it:ind-taux-terminal}).
    
  \item\label{it:ind-overlapping-case} If the carrier base $\mu$ in some
    $(\mbc'',t'',J')$ forms a principal overlapping pair with its dual
    we first construct the induced auxiliary tree
    $\Taux{\mbc'',t'',J''}$ rooted at $(\mbc'',J')$, where
    $J'' = J'\cup\inter\mu$.

    Next we take
    \[R = R(\mbc'',t'',J'',\mu) = \left\{ \left( \mbc_i,t_i,J_i;
        \calH(t_i)\right) \right \}\]
    to be the set of all witnessing terminals of $\Taux{\mbc'',J''}$
    in which $\inter\mu$ is an $\calH(t_i)$-maximal block overlapping
    pair (recall Convention \ref{conv:renaming}).
    
    If this set is empty then $(\mbc'',J')$ is declared to be a
    \define{halted terminal} and no further descendants are
    added. Otherwise the following number is defined and computable
    (recall Definition \ref{defn:period-lambda}):
    \begin{equation}\label{eq:ind-period-bound}
      \period{(\mbc'',t'',J')}{\mu}=\max_{\left(\mbc_i,t_i,J_i;\calH(t_i)\right)\in
        R}\period{\calH(t_i)}{\mu}.       
    \end{equation}
    We now continue growing the path at $(\mbc'',t'',J'')$ in
    $\ADtwo{\mbc',J'}$ induced by $t'\subset \mbc'$. By Corollary
    \ref{cor:periodicity-bound} and Lemma
    \ref{lem:periodicity-non-decr}, $\mu$ will not be the carrier base
    more than $B\cdot\period{(\mbc'',t'',J')}{\mu}$ times in a row.

    Once $\mu$ ceases to be a carrier base we go back to
    (\ref{it:ind-std-growth}) or (\ref{it:ind-overlapping-case}) as
    appropriate.
  \end{enumerate}
\end{enumerate}
  
\begin{figure}[htb]
  \centering
  \begin{tikzpicture}[xscale=1, every node/.style={scale=0.8},inner sep=0.5mm]
    \tikzstyle{descr} = [
    rectangle,draw=black!50,
    execute at begin node={\begin{varwidth}{10em}},
      execute at end node={\end{varwidth}}]
    \tikzstyle{dot}=[circle,fill=black]
    
    \node (CTv) at (0,-1) [dot]{}; 
    \node[above] (CT) at (CTv.north) {$(\mbc,t)$};
    \node[above]  (CTd) at (CT.north) [descr] {$\basepair\lambda$ overlapping.};
    \node (c1t1v) at (0,-2) [dot]{};
    \node[below] (c1t1) at (c1t1v.south) {$(\mbc_1,t_1,J)$};
    \node[below] (c1t1d) at (c1t1.south) [descr]{$\inter\lambda$ not $\calH(t_1)$-maximal.};
    
    \node (c2t2v) at (3,-2) [dot]{};
    
    \node (c3t3v) at (3,-5) [dot]{};
    
    \node (c4t4v) at (6,-5) [dot]{};

    \node (c5t5v) at (6,-6) [dot]{};
    
    \node (c6t6v) at (9,-2) [dot]{};
    
    \node (c7t7v) at (9,-3) [dot]{};
    
    \node[above] (c2t2) at (c2t2v.north) {$(\mbc_2,t_2,J)$};
    \node[above]  (c2t2d) at (c2t2.north) [descr] {$\basepair\mu$ overlapping.};

    \node[above] (c4t4) at (c4t4v.north) {$(\mbc_4,t_4,J')$};
    \node[above]  (c4t4d) at (c4t4.north west) [descr] {$\basepair\delta$ overlapping.};
    
    \node[above] (c6t6) at (c6t6v.north) {$(\mbc_6,t_6,J)$};
    \node[above]  (c6t6d) at (c6t6.north) [descr] {$\basepair\eta$ overlapping.};
    
    \node[below] (c3t3) at (c3t3v.south) {$(\mbc_3,t_3,J')$};
    \node[below] (c3t3d) at (c3t3.south) [descr]{$\inter\mu$ and
      $\inter\lambda$ not $\calH(t_3)$-maximal.};
    
    \node[below] (c5t5) at (c5t5v.south) {$(\mbc_5,t_5,J'')$}; 
    \node[below] (c5t5d) at (c5t5.south) [descr]{$\inter\mu$ is
      $\calH(t_5)$-maximal. Continue growing at $\mbc_2$.}; 
    \node[below] (c7t7) at (c7t7v.south) {$(\mbc_7,t_7,J\cup\inter\eta)$}; 
    \node[below] (c7t7d) at
    (c7t7.south) [descr]{$\inter\lambda$ is $\calH(t_y)$-maximal. In
      $\mbc$ $\period{t}{\lambda}$ is at most
      $\period{\calH(t_5)}{\lambda}$.};
    
    \node (endmu) at (6,-2) [dot]{};

    \node (d1) at ($(c1t1v)!0.5!(c2t2v)$) {$\cdots$};
    \node (d2) at ($(c3t3v)!0.5!(c4t4v)$) {$\cdots$};
    \node (d3) at ($(c2t2v)!0.5!(endmu)$) {$\cdots$};
    \node (d4) at ($(endmu)!0.5!(c6t6v)$) {$\cdots$};
    
    \draw[very thick,->] (endmu)+(-1,0) -- (endmu);
    \draw[very thick,->] (endmu) -- +(1,0);
    \draw[very thick,->] (c1t1v) -- +(1,0);
    \draw[very thick,->] (c2t2v) -- +(1,0);
    \draw[very thick,->] (c2t2v)+(-1,0) -- (c2t2v);
    \draw[very thick,->] (c6t6v)+(-1,0) -- (c6t6v);
    \draw[very thick,->] (c3t3v) -- +(1,0);
    \draw[very thick,->] (c4t4v)+(-1,0) -- (c4t4v);
    \draw[very thick,->] (CTv) -- (c1t1v);
    \draw[very thick,->] (c2t2v) -- (c3t3v);
    \draw[very thick,->] (c4t4v) -- (c5t5v);
    \draw[very thick,->] (c6t6v) -- (c7t7v);
    
    \draw[decorate,decoration={brace,amplitude=0.3cm}] (5.9,-2.1) --
    (3.1,-2.1);
    \node (periodicity-bound) at (5,-3) [descr] {$\mu$ is carrier base
      at most $N\period{\calH(t_5)}{\mu}$ times.};
  \end{tikzpicture}
  \caption{The tree $\Taux{\mbc,t,\inter\lambda}$, here
    $J=\inter\lambda, J' = \inter\lambda\cup\inter\mu$, and
    $J'' = J' \cup \inter\delta$. The band complexes are numbered in
    order of appearance. Whenever an overlapping pair occurs, an
    auxiliary edge is constructed. $\mbc_4$ and $\mbc_6$ are halted
    terminals. $\mbc_5$ and $\mbc_7$ are witnessing terminals.}
\label{fig:ind-aux-tree}
\end{figure}
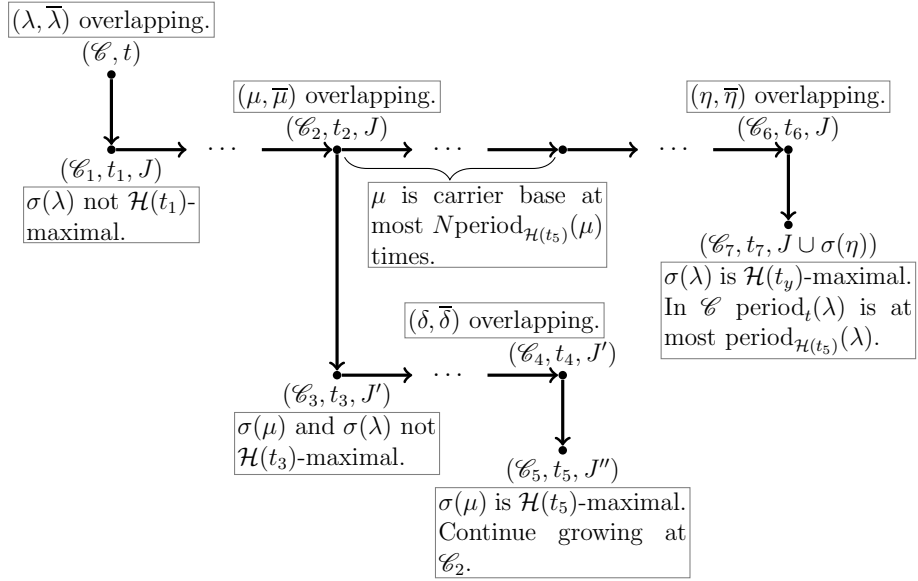

The induced auxiliary tree can be thought of as being constructed one
vertex at a time, see Figure \ref{fig:ind-aux-tree}, as opposed to
a branching process. Another important distinction is that induced
auxiliary trees do not have terminals that are halted by
inadmissibility.  

\begin{lem}\label{lem:ind-tree-finite}
  $\Taux{\mbc,t,J}$ is finite.
\end{lem}
\begin{proof}
  By Lemma \ref{lem:depth-bound} there is a bound on the number of
  auxiliary edges in any path. The finiteness of horizontal paths
  follows from the definition of a Rips process, i.e. $\size{t}$
  (Definition \ref{defn:track-size}) effectively bounds the length of
  such paths.
\end{proof}

\begin{lem}\label{lem:ind-has-a-bound}
  If $\basepair\lambda$ is a principal overlapping pair in $\mbc$ then
  there is a witnessing terminal
  $\left(\mbc_T,t_T,J_T;\calH(t_T)\right)$ in
  $\Taux{\mbc,t,\inter\lambda}$ in which $\inter\lambda$ (following
  renaming convention \ref{conv:renaming}) is a
  $\calH(t_T)$-maximal periodic block such that \[
  \period{t}{\lambda} \leq \period{\calH(t_T)}{\lambda}.
  \]
\end{lem}
\begin{proof}
  We first must show that such a witnessing terminal
  $\left(\mbc_T,t_T,J_t;\calH(t_T)\right)$ exists. Suppose towards a
  contradiction that this was not the case, then the top level
  elimination tree in $\Taux{\mbc,t,\inter\lambda}$ must end in a
  halted terminal $(\mbc_1,t_1,\inter\lambda)$, otherwise all the
  bases are moved onto $\inter\lambda$, and this periodic block will be maximal (since it's unique.)

  Let $\basepair{\lambda_1}$ be the overlapping pair in $\mbc_1$ and
  let $J_2 = \inter\lambda\cup\inter{\lambda_1}$. Let
  $(\mbc_2,t_2,J_2)$ be the auxiliary descendant, as in
  (\ref{it:ind-aux-descendants}), of $\mbc_1$. The horizontal path in
  $\Taux{\mbc,t,\inter\lambda}$ starting at $(\mbc_2,t_2,J_2)$ must
  end with a halted terminal, otherwise it ends in a witnessing
  terminal $\left(\mbc_T,t_T,J_2;\calH(t_T)\right)$ with either
  $\inter\lambda$ or $\inter\lambda_1$ $\calH(t_T)$-maximal. This
  either contradicts the hypothesis that $\mbc_1$ is a halted terminal
  or the hypothesis that $\Taux{\mbc,t,\inter\lambda}$ has no
  witnessing terminals with $\inter\lambda$ maximal.

  Continuing in this fashion we obtain a sequence of halted terminals
  of increasing depth
  \[ \left(\mbc_1,t_1,J_1\right), \left(\mbc_2,t_2,J_2\right), \ldots,
  \left(\mbc_F,t_F,J_F\right)
  \]
  where $\basepair\lambda_i$ is the principal overlapping pair in
  $\mbc_i$ and $J_{i+1}= J_i \cup\lambda_i$. Since
  $\tau(\mbc_i,J_i) > \tau(\mbc_{i+1},J_{i+1})$ this sequence is
  finite which forces some principal overlapping pair $\lambda_j$ in
  $\mbc_F$ to be $\calH(t_F)$-maximal, contradicting the fact $\mbc_j$
  is a halted terminal.

  It therefore follows that there is a witnessing terminal
  $\left(\mbc_T,t_T,J_t;\calH(t_T)\right)$ in
  $\Taux{\mbc,t,\inter\lambda}$ in which $\inter\lambda$ is
  maximal. The desired upper bound on $\period{t}{\lambda}$ now
  follows immediately from Definition \ref{defn:period-lambda}, Lemma
  \ref{lem:periodicity-non-decr}, and Corollary
  \ref{cor:periodicity-bound}.
\end{proof}

So far we have shown that the induced auxiliary tree contains a
witnessing terminal whose combinatorial periodicity bounds the actual
periodicity. We now bound the periodicity, for all tracks.

\begin{prop}\label{prop:the-periodicity-bound-works!}
  Let $\lambda$ be a principal overlapping pair in a band complex
  $\mbc$ in $\ADtwo{\mbc}$, then for all automorphically
  minimal $\kappa$-tracks $t$ efficiently carried by $\mbc$ the
  following holds\[
  \period t \lambda \leq \period \mbc \lambda,
  \] where $\period\mbc\lambda$ is the computable function given by
  Definition \ref{def:periodicity-bound}.
\end{prop}
\begin{proof}
  We first show that for any automorphically minimal $\kappa$-track
  $t$ efficiently carried by $\mbc$ we have a natural containment:
  \begin{equation}
    \label{eq:tree-containment}
    \Taux{\mbc,t,\inter\lambda} \subset \Taux{\mbc,\inter\lambda}.
  \end{equation}
  We will show this by analyzing how $\Taux{\mbc,t,\inter\lambda}$ is
  constructed by
  adding one band complex at a time.

  Going through the construction algorithms point-by-point, by
  Proposition \ref{prop:admissible-contains-II}, the only
  problem that could arise is in step (\ref{it:ind-overlapping-case})
  of the construction of the auxiliary trees. It could be that for
  some $(\mbc',t',J')$ with principal overlapping pair $\basepair\mu$,
  $(\mbc',t',J')$ is not a halted terminal and
  $\period{(\mbc',t',J')}{\mu}$ from (\ref{eq:ind-period-bound}) is
  greater than $\period{(\mbc',J')}{\mu}$ from (\ref{eq-period-bound}).

  Note however that in the recursive construction of
  $\Taux{\mbc,t,\inter\lambda}$, we must first construct
  $\Taux{\mbc',t',J'\cup\lambda}$ before adding a ``horizontal''
  child of $(\mbc',t',J')$. It follows, by the definition and
  properties of auxiliary children, that the next vertex added to
  $\Taux{\mbc,t,\inter\lambda}$ is still contained in
  $\Taux{\mbc,\inter\lambda}$.

  Continuing in this manner it is obvious (the reader is, of course,
  free to supply their own argument by induction on relative
  $\tau$-complexity) that the set $R$ appearing in
  (\ref{eq:ind-period-bound}) of Case (\ref{it:ind-overlapping-case})
  in the construction of the induced auxiliary tree is a subset of $W$
  appearing in (\ref{eq-period-bound}) of Case
  (\ref{it:overlapping-case}) of the construction of the auxiliary
  tree. We conclude that
  \[ \period{(\mbc',t',J')}{\mu}\leq\period{(\mbc',J')}{\mu}.\]
  (\ref{eq:tree-containment}) now follows; thus by Lemma
  \ref{lem:ind-has-a-bound} a witnessing terminal of
  $\Taux{\mbc,\inter\lambda}$ bounds $\period{t}{\lambda}$ from above.
\end{proof}

\subsection{The proof of Theorem \ref{thm:main}: a description of the
  main algorithm}\label{sec:main-proof}
Suppose we are given a finite 2-complex $\complex{C}$ such that
$\pi_1(\complex{C})$ has no elements of order 2, a solution to the
word problem for $\pi_1(\complex{C})$, an acylindricity constant
$\kappa$, and a finite collection
\[ S=\big\{ \{h_i\}_{i\in I_n} \mid n = 1,\ldots, m \big\}\]
of finite generating sets of subgroups
$\elliptics = \big\{ \bk{h_i}_{i\in I_n}\big\}_{n=1}^m$.

We start by replacing $\complex{C}$ by $\complex{C}_S$ given in
Section \ref{sec:relative-splittings}. This can be done
algorithmically. Using the construction of Section
\ref{sec:correspondences}, we see that for any track
$t\subset \complex C$ there is a corresponding band complex that
carries it efficiently. These band complexes can be effectively
enumerated; thus it is possible to construct the finite set
$\{\mbc_1,\ldots,\mbc_{n_C}\}$ of band complexes given by Proposition
\ref{prop:finite-ubcs} (See also Section \ref{sec:ETone}
(\ref{it:et-first-children})). This gives the first level of our
elimination tree. We will now define the ultimate elimination tree
$\ETP{\complex C}$ in this paper. Here is final inadmissibility
criterion:

\begin{defn}\label{defn:period-inadmissible}
  Let $\mbc$ be a band complex. A path \[
  p: \mbc_u \to \cdots \to \mbc_v
  \] in $\ADtwo{\mbc}$ (recall section \ref{sec:C-T-inadmissble}) is
  called \define{periodicity-inadmissible} if 
  \begin{enumerate}
  \item $p$ is a $\mu$-periodic path (Definition
    \ref{defn:mu-periodic}) for some base $\mu$ in $\mbc_u$, and
  \item the length of $p$ is greater than $N\cdot\period{\mbc_u}{\mu}$
    (Definition \ref{def:periodicity-bound}), where $N$ is the number
    of bases in $\mbc_u$.
  \end{enumerate}
\end{defn}

$\ETP{\complex C}$ is constructed identically to $\ADtwo{\complex C}$,
but we also forbid periodicity inadmissible paths. By Proposition
\ref{prop:taux-finite} we can compute the periodicity bound and thus
effectively decide whether a path is periodicity-inadmissible, so the
resulting tree is finite.

To help the reader, however, we will give here a more explicit
construction of $\ETP{\complex C}$ that will summarize the important
results of this paper. We start with our root, the polygonal 2-complex
$\complex C$. We add the descendants $\mbc_1,\ldots,\mbc_{n_C}$. We
then build $\ETP{\complex C}$ ``one generation at a
time'' as follows:
\begin{enumerate}[(1)]
\item For every admissible vertex without descendants we add the
  descendants as described in Definition \ref{defn:ETTwo}.
\item If a freshly added descendant is merging inadmissible
  (Definition \ref{def:merging-inadmissible}), we declare it
  inadmissible and stop adding its descendants.
\item We now consider every directed path constructed so far in our
  elimination tree. If a path $\mbc_u \to \ldots \to \mbc_v$ is
  either,
 \begin{itemize}
 \item $\kappa$-inadmissible (Definition
   \ref{defn:kappa-inadmissible}),
 \item repetition inadmissible (Definition
   \ref{defn:rep-inadmissible}),
 \item C-T-inadmissible (Definition \ref{defn:CT-inadmissible}), or 
 \item periodicity inadmissible,
 \end{itemize}
 then we declare the last vertex $\mbc_v$ to be inadmissible and stop
 adding its descendants. 
\end{enumerate}

\begin{figure}[htb]
  \centering
  \begin{tikzpicture}
    \draw[fill=gray!20] (0,1) -- (4,-5.5) -- (-4,-5.5)--cycle;
    \node (C) at (0,0) {$\complex C$};
    \node (C1) at (-1,-1.5) {$\mbc_1$};
    \node (Ci) at (0,-1.5) {$\mbc_i$};
    \node (Cn) at (1,-1.5) {$\mbc_n$};
    \node (Cu) at (2,-3) {$\mbc_u$};
    \node (Cv) at (-1,-4.5) {$\mbc_v$};
    \draw (-0.5,-1.5) node{$\ldots$} (0.5,-1.5) node{$\ldots$};
    \draw[->](C)--(C1);
    \draw[->](C)--(Ci);
    \draw[->](C)--(Cn);
    \draw[thin, decoration = {zigzag},decorate](Ci) -- (Cu);
    \draw[ultra thick, decoration = {zigzag},decorate](Cu) -- (Cv);
    \draw[ultra thin] (Cv)+(-0.25,0.5) -- +(0.25,-0.5);
    \draw[ultra thin] (Cv)+(0.25,0.5) -- +(-0.25,-0.5);
    
  \end{tikzpicture}
  \caption{As we build $\ETP{\complex C}$ we have a path (drawn thick)
    in which
    $\mbc_u$ and $\mbc_v$ are equal. This path is therefore repetition
    inadmissible and $\mbc_v$ has no further descendants.}
  \label{fig:rep-inadmissible-eg}
\end{figure}
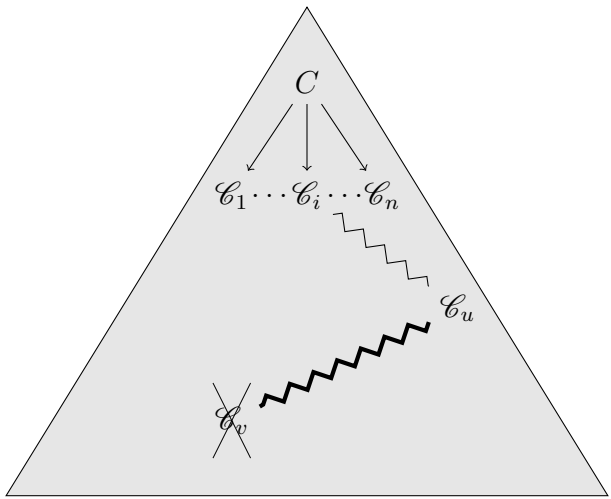

We will now argue that tree $\ETP{\complex C}$ is finite and can be
algorithmically constructed. First note that at every step we can
algorithmically construct the set of descendants of a band complex
(Definition \ref{defn:Derived-BC}), in particular $\ETP{\complex C}$
has finite branching. By König's lemma it is therefore enough to show
that $\ETP{\complex C}$ has no directed infinite paths.

Suppose towards a contradiction that this was the case. By Theorem
\ref{thm:branch-types} any such infinite branch must be either
thinning, quadratic or superquadratic. In the thinning and quadratic
cases, such a branch must either have a $\kappa$-inadmissible subpath
or must contain a repetition (see Section \ref{sec:repetition}.) In
the superquadratic case any infinite branch must either have a tail
with infinitely many C-T cycles or the tail must be $\mu$-periodic for
some base $\mu$ (see the proof of Proposition
\ref{prop:periodicity-reduction}.) In both of these cases such a tail
will have either a C-T-inadmissible or a periodicity inadmissible
initial segment. It therefore follows that $\ETP{\complex C}$ has no
infinite branches, furthermore the four inadmissibility criteria are
algorithmically verifiable; thus
$\ETP{\complex C} \subset \ETone{\complex C}$ is algorithmically
constructible.

The leaves of $\ETP{\complex C}$ give a subset of all possible tracks
of $\complex C$ (recall Section \ref{sec:set-of-tracks}.) We will now
show that this subset contains a representative of every
automorphically minimal $\kappa$-track.

Suppose that there was some automorphically minimal $\kappa$-track
$t \subset \complex C$ that induced a path
$p:\complex C \to \mbc_1 \to \ldots \to \mbc_l$ in
$\ETone{\complex C}$ with $\mbc_l$ terminal, which isn't contained in
$\ETP{\complex C}$. Then, by definition of $\ETP{\complex C}$, $p$
must either contain a $\kappa$-inadmissible, a repetition
inadmissible, a C-T-inadmissible, or a periodicity inadmissible
subpath. Propositions \ref{prop:kappa-inadmissible},
\ref{prop:repetition}, Corollary \ref{cor:not-too-many-cycles} and
Proposition \ref{prop:the-periodicity-bound-works!} cover each of
these cases and contradict the assumption that $t$ is an
automorphically minimal $\kappa$-track. It follows that the admissible
leaves of $\ETP{\complex C}$ give a set of
tracks \[\left\{t_1,\ldots,t_{n(\complex C, \kappa,S)}\right\}\] in
$\complex C$ that satisfy the requirements of Theorem
\ref{thm:main}.\hfill\qedsymbol

\bibliographystyle{alpha} \bibliography{biblio.bib}

\end{document}